\documentclass[11pt, letterpaper]{amsart}
\usepackage{amsmath}
\usepackage{amssymb}
\usepackage[left=1in,right=1in,bottom=0.8in,top=0.8in]{geometry}
\usepackage{amsfonts}
\usepackage{graphicx}
\usepackage[font=small,labelfont=bf]{caption}
\usepackage{epstopdf}
\usepackage{hyperref}
\usepackage[dvipsnames]{xcolor}
\usepackage{amsthm}
\usepackage[T1]{fontenc}
\usepackage{float}
\usepackage{pgfplots}
\usepackage{listings}
\usepackage{longtable}
\usepackage{mathrsfs}
\usepackage{bbm}
\usepackage{euscript}
\usepackage{mathtools}
\usepackage{tikz}
\usepackage{comment}
\usepackage{tikz-cd}
\usepackage{amsxtra}
\usepackage{pifont}
\usepackage{amsbsy}
\usepackage{wasysym}
\usepackage{calligra}
\usepackage{url}
\usepackage{enumerate}
\usepackage[dvipsnames]{xcolor}
\usepackage[
backend=biber,
style=alphabetic,
]{biblatex}

\makeatletter
\renewcommand{\tocsection}[3]{%
  \indentlabel{\@ifnotempty{#2}{\bfseries\ignorespaces#1 #2.\quad}}\bfseries#3}
  \renewcommand{\tocsubsection}[3]{%
  \indentlabel{\hspace{2.5em}\@ifnotempty{#2}{\ignorespaces#1 #2.\quad}}#3}
\makeatother

\newtheorem{theorem}{Theorem}[section]
\newtheorem*{backgroundtheorem}{Theorem}
\newtheorem{corollary}[theorem]{Corollary}
\newtheorem{proposition}[theorem]{Proposition}
\newtheorem{lemma}[theorem]{Lemma}

\theoremstyle{definition}
\newtheorem{definition}[theorem]{Definition}

\newtheorem{setup}[theorem]{Setup}
\newtheorem{example}[theorem]{Example}
\newtheorem{notation}[theorem]{Notation}

\theoremstyle{remark}
\newtheorem{remark}[theorem]{Remark}

\definecolor{energy}{RGB}{114,0,172}
\definecolor{refblue}{RGB}{30,90,180}
\definecolor{freq}{RGB}{45,177,93}
\definecolor{spin}{RGB}{251,0,29}
\definecolor{signal}{RGB}{203,23,206}
\definecolor{circle}{RGB}{217,86,16}
\definecolor{average}{RGB}{203,23,206}

\colorlet{shadecolor}{gray!20}
\pgfplotsset{compat=1.9}

\usepgflibrary{fpu}

\numberwithin{equation}{section}
\newcommand{\QQ}{\mathbb Q}
\newcommand{\PP}{\mathbb P}
\newcommand{\RR}{\mathbb R}
\newcommand{\CC}{\mathbb C}
\newcommand{\ZZ}{\mathbb Z}
\newcommand{\OO}{\mathcal O}
\newcommand{\cA}{\mathcal A}
\newcommand{\cB}{\mathcal B}

\newcommand{\cD}{\mathcal D}
\newcommand{\cE}{\mathcal E}
\newcommand{\cF}{\mathcal F}
\newcommand{\cG}{\mathcal G}
\newcommand{\cH}{\mathcal H}
\newcommand{\cI}{\mathcal I}
\newcommand{\cK}{\mathcal K}
\newcommand{\cL}{\mathcal L}
\newcommand{\cM}{\mathcal M}

\newcommand{\cP}{\mathcal P}
\newcommand{\cQ}{\mathcal Q}

\newcommand{\cS}{\mathcal S}
\newcommand{\cT}{\mathcal T}
\newcommand{\cU}{\mathcal U}
\newcommand{\cV}{\mathcal V}
\newcommand{\cW}{\mathcal W}
\newcommand{\bS}{\mathbf S}
\newcommand{\rk}{\operatorname{rk}}
\newcommand{\Hom}{\operatorname{Hom}}
\newcommand{\Pic}{\operatorname{Pic}}
\newcommand{\Gr}{\operatorname{Gr}}
\newcommand{\pr}{\operatorname{pr}}
\newcommand{\Sym}{\operatorname{Sym}}
\newcommand{\length}{\ell}

\newcommand{\id}{\operatorname{id}}
\newcommand{\im}{\operatorname{Im}}
\newcommand{\End}{\operatorname{End}}
\newcommand{\GL}{\operatorname{GL}}

\newcommand{\SO}{\operatorname{SO}}
\newcommand{\Sp}{\operatorname{Sp}}

\hypersetup{
    colorlinks=true,
    linkcolor=refblue,
    citecolor=signal,
    urlcolor=signal
}

\author[Soham Ghosh]{Soham Ghosh}
\email{soham13@uw.edu}
\thanks{SG was supported in part by a Simons Dissertation Fellowship SFI-MPS-SDF-00014771, NSF CAREER DMS-2044564 and NSF FRG DMS-2151718 grants.}
\author[S\'andor J Kov\'acs]{S\'andor J Kov\'acs}
\email{skovacs@uw.edu}
\address{University of Washington, Department of Mathematics, Seattle,WA 98195,USA}
\thanks{SJK was supported in part by NSF Grant DMS-2100389, and the Craig McKibben
and Sarah Merner Endowed Professorship in Mathematics.}
\title{Characterizations of products of projective spaces and quadrics}
\begin{document}

\renewcommand\thesubsection{\thesection.\Alph{subsection}}
\numberwithin{equation}{theorem}
\renewcommand\theequation{\thetheorem.\arabic{equation}}

\begin{abstract}
We characterize products of projective spaces and smooth quadrics using Schur functors and positivity of vector bundles. First, we obtain characterizations of projective spaces and smooth quadrics using positivity of arbitrary Schur functors of extensions of a unipotent bundle by the tangent bundle, generalizing several existing results. We then prove that almost nef regular foliations whose direct sum is a subbundle of the tangent bundle determine a product decomposition, if the corresponding Schur functors admit suitable positive Schur submodules. For smooth weak Fano varieties similar results are obtained by replacing almost nefness and individual integrability of the foliations by the integrability of their total image. As special cases, we obtain product analogues of well-known characterizations of projective spaces and quadrics.
\end{abstract}

\maketitle
\tableofcontents

\section{Introduction and Main results}

Positivity of the tangent bundle imposes strong restrictions on the
geometry of a projective variety, leading to several characterizations of projective spaces. Some of the most well-known results are due to Kobayashi-Ochiai \cite{KO73}, Mori \cite{Mori79},  Wahl \cite{Wahl83}, Andreatta-Wi\'sniewski \cite{AW01}, Druel \cite{Druel04}, and Araujo-Druel-Kov\'acs \cite{ADK08}. These results progressively pass from positivity of the tangent bundle to its exterior and tensor powers, and therefore naturally lead to characterizations that also include smooth hyperquadrics $Q_n$. A particularly powerful characterization is due to Druel and Paris \cite{DruelParis13}, who proved that a nonzero morphism $\det\cF\rightarrow T_X^{\otimes t}$, with $\cF$ an ample vector bundle and $\rk\mathcal F\geq t\geq 1$, is enough to force $X$ to be $\PP^n$ or $Q_n$. Using this, more recently the first author proved the following generalization (along with a symmetric power analogue, see \cite[Theorem~1.2]{Ghosh26}) of several of the above characterizations.

\begin{backgroundtheorem}[{\cite[Theorem~1.1]{Ghosh26}}]
Let $X$ be a smooth complex projective variety of dimension $n$. If there exists an ample vector bundle $\mathcal{E}$ of rank $r$ on $X$ and a positive integer $p\leq \min\{r, n\}$ such that $\bigwedge^p\mathcal{E}\subseteq\bigwedge^pT_X$, then either $X\cong \mathbb P^n$, or $p=n$ and $(X, \mathcal{E})\cong (Q_p, \mathcal{O}_{Q_p}(1)^{\oplus r})$ where $Q_p\subseteq \mathbb{P}^{p+1}$ is a smooth hyperquadric in $\mathbb{P}^{p+1}$. Furthermore, if $X\cong \mathbb{P}^n$, then $\det\mathcal{E}\cong\mathcal{O}_{\mathbb{P}^n}(l)$ with $l=r$ or $r+1$ and:
\begin{enumerate}
\item If $r<n$, then $\mathcal{E}\cong \mathcal{O}_{\mathbb{P}^n}(1)^{\oplus r}$.
\item If $r=n$ and $l=n$, then $\mathcal{E}\cong \mathcal{O}_{\mathbb{P}^n}(1)^{\oplus n}$.
\item If $r=n$ and $l=n+1$:
\begin{enumerate}
\item If $1\leq p\leq n-1$, then $\mathcal{E}\cong T_{\mathbb{P}^n}$.
\item If $p=n$, then either
$\mathcal{E}\cong T_{\mathbb{P}^n}$ or $\mathcal{E}\cong\mathcal{O}_{\mathbb{P}^n}(2)\oplus \mathcal{O}_{\mathbb{P}^n}(1)^{\oplus(n-1)}$.
\end{enumerate}
\end{enumerate}
\end{backgroundtheorem}

There are several other characterizations of $\PP^n$ from the perspective of deformations of rational curves, differential geometry and topology, for which we refer the readers to \cite{LOY11} and \cite{Deb15}. The next simplest class of varieties that one might hope to characterize are products of projective spaces and smooth hyperquadrics. Contrary to the case of $\PP^n$ or $Q_n$, only very few characterizations of products of projective spaces, and even fewer for products of smooth hyperquadrics, are known. 
Ochetta \cite{Occhetta06} obtained a characterization of products of projective spaces using unsplit families of rational curves. Based on this, another characterization using the geometry of minimal rational curves and their tangents was obtained by Araujo \cite[Theorem~1.3]{Araujo06}. This was subsequently utilized by Fujita \cite[Theorem~2.11]{Fujita14} to obtain a criterion for products of projective spaces using lengths of extremal rays.

A general guiding question in the study of characterizations of products of projective spaces is a well-known conjecture of Mukai \cite{Mukai88}. A similar Mukai-type conjecture was recently proposed in \cite{Gongyo26} using an invariant called nef complexity $c_X$ for smooth Fano varieties $X$. This has been settled in the recent work \cite{EnwrightEtAl25}, with weaker assumptions on singularities of the Fano variety, leading to a new characterization of products of projective spaces among Fano varieties, in terms of nef complexity $c_X$. In a recent work \cite{Watanabe26}, Watanabe clarified the full spectrum of possible nef complexities. We refer the reader to \cite{EnwrightEtAl25} for an account of Mukai's conjecture, and several references for related work. 

The purpose of this paper is threefold. Our first main result extends the results of \cite{Ghosh26} to arbitrary Schur functors and allows the tangent bundle to be replaced by an extension by a unipotent vector bundle, thereby yielding the following characterization of projective space and smooth hyperquadrics. Furthermore, Theorem~\ref{thm:Main1} also plays an important role in the inductive steps for the proofs of both Theorem~\ref{thm:Main2} and ~\ref{thm:Main3} below.

\begin{theorem}[=Theorem~\ref{thm:generalizedKovacs}]\label{thm:Main1}
Let $X$ be a smooth projective complex variety of dimension $n>0$, let
$\cU$ be a unipotent vector bundle on $X$ (that is, $\cU$ admits a finite filtration by vector subbundles with trivial successive quotients), and let $\cD$ fit into an exact sequence
\begin{equation}\label{eq:unipotent-extension0}
 0\rightarrow T_X\rightarrow\cD\rightarrow\cU\rightarrow0.
\end{equation}
Further let $\cE$ be an ample vector bundle of rank $r$ on $X$ and $\lambda\vdash m$ a partition of length $\length(\lambda)\leq r$, with associated Schur functor $\bS_{\lambda}$. Assume that there is an injective morphism $\bS_\lambda\cE\hookrightarrow\bS_\lambda\cD$ of coherent $\OO_X$-modules. Then either: 
\begin{enumerate}
    \item[(a)] $(X, \det\cE)\simeq(\PP^n, \OO_{\PP^n}(l))$ with $r\leq l\leq \frac{r(n+1)}{n}$ or,
    \item[(b)] $(X, \cE)\simeq (Q_n, \OO_{Q_n}(1)^{\oplus r})$ with $n\ge2$, $r\leq n$. Furthermore, the morphism $\bS_\lambda\cE\hookrightarrow\bS_\lambda\cD$ factors as $\bS_\lambda\cE\hookrightarrow\bS_\lambda T_{Q_n}\hookrightarrow\bS_\lambda\cD$ and the partition $\lambda\vdash m$ is $(k^r)$, where either $k$ is even or $r=n$.
\end{enumerate}
Furthermore, in case (a), the bundle $\cE$ on $\PP^n$ is determined as:
\begin{enumerate}
 \item[(i)] if $r<n$, then $l=r$ and
       $\cE\simeq\OO_{\PP^n}(1)^{\oplus r}$;
 \item[(ii)] if $r=n$ and $l=n$, then
       $\cE\simeq\OO_{\PP^n}(1)^{\oplus n}$;
 \item[(iii)] if $r=n$ and $l=n+1$, then $\cE\simeq T_{\PP^n}$ or $\cE\simeq \OO_{\PP^n}(2)\oplus\OO_{\PP^n}(1)^{\oplus(n-1)}$. The second alternative can occur only when
       $\lambda=(k^n)$ for some $k\ge1$;
 \item[(iv)] if $r>n$, then $l=r$ and $\cE\simeq\OO_{\PP^n}(1)^{\oplus r}$.
\end{enumerate}
\end{theorem}

Next we characterize products of projective spaces and smooth hyperquadrics among all smooth projective varieties using Schur-functor inclusions of nef vector bundles into almost nef regular foliations, leading to our next main result. 

We briefly recall the terminology used below, see \cite[\S~2.1]{Druel17} and \cite[\S~2.B]{Horing07}. For a smooth complex variety $X$, a vector bundle $\cV$ on $X$ is \emph{almost nef} if there exists a countable union $Z\subset X$ of proper closed subvarieties, such that the restriction of $\cV$ to every irreducible curve not contained in $Z$ is nef. A coherent subsheaf $\cF\subseteq T_X$ is \emph{integrable} (or \emph{involutive}) if it is closed under the natural Lie bracket on $T_X$. A \emph{foliation} is an integrable coherent subsheaf $\cF\subseteq T_X$ that is \emph{saturated}, i.e., the quotient $T_X/\cF$ is torsion-free. A foliation is \emph{regular} if $T_X/\cF$ is locally free, or equivalently, $\cF$ is an integrable
vector subbundle of $T_X$. We allow the zero foliation and $T_X$ itself.

\begin{theorem}[=Theorem~\ref{thm:productcharacterization}]\label{thm:Main2}
Let $X$ be a smooth projective complex variety of dimension $n>0$ and let $r>0$ be an integer. For $1\le i\le r$, let $\cE_i\subseteq T_X$ be an almost nef regular foliation of rank $n_i\geq1$, $\cV_i$ a semiample vector bundle of rank $r_i\ge1$, and $\lambda_i$ a nonzero partition of length $\ell(\lambda_i)\le \min\{r_i,n_i\}$ with associated Schur functor $\bS_{\lambda_i}$. Assume that:
\begin{enumerate}
\item the natural morphism $\bigoplus_{i=1}^r\cE_i\rightarrow T_X$ is a subbundle inclusion,
\item for every $1\leq i\leq r$, there is an injective morphism, $\theta_i:\bS_{\lambda_i}\cV_i\hookrightarrow\bS_{\lambda_i}\cE_i$ of coherent $\OO_X$-modules and,
\item the vector bundle $\cA:=\bigotimes_{i=1}^r\cV_i$ is ample on $X$. 
\end{enumerate}
Then, possibly after reordering, there is an isomorphism
\begin{equation}\label{eq:main-product}
X\simeq\prod_{i=1}^r Z_i,\qquad
\dim X=n=\sum_{i=1}^{r} n_i,\qquad
\cE_i\simeq\pr_i^*T_{Z_i}\subseteq T_X,
\end{equation}
where each $Z_i$ is $\PP^{n_i}$ or a smooth hyperquadric $Q_{n_i}\subset\PP^{n_i+1}$ with $n_i\geq 2$ and $\pr_i:\prod_{i=1}^r Z_i\rightarrow Z_i$ is the natural projection. Furthermore, $r_i\leq n_i$ and for all $1\leq i\leq r$, there are ample bundles $\cT_i$ on $Z_i$ such that $\cV_i\simeq\pr_i^*\cT_i$ and each $\cT_i$ admits an injective morphism $\bS_{\lambda_i}\cT_i\rightarrow\bS_{\lambda_i}T_{Z_i}$ of coherent sheaves. These bundles $\cT_i$ on $Z_i$ are classified as follows.
\begin{enumerate}[(a)]
\renewcommand{\labelenumi}{\textup{(\alph{enumi})}}
\item If $Z_i\simeq\PP^{n_i}$ and $r_i<n_i$, then $\cT_i\simeq\OO_{\PP^{n_i}}(1)^{\oplus r_i}$,
\item If $Z_i\simeq\PP^{n_i}$, $r_i=n_i$ and $\det\cT_i\simeq\OO_{\PP^{n_i}}(n_i)$, then $\cT_i\simeq\OO_{\PP^{n_i}}(1)^{\oplus n_i}$,
\item If $Z_i\simeq\PP^{n_i}$, $r_i=n_i$ and $\det\cT_i\simeq\OO_{\PP^{n_i}}(n_i+1)$, then 
\[\cT_i\simeq T_{\PP^{n_i}} \text{ or } \cT_i\simeq\OO_{\PP^{n_i}}(2)\oplus\OO_{\PP^{n_i}}(1)^{\oplus(n_i-1)}.\]
The second possibility can occur if and only if $\lambda_i=(k^{n_i})$ for some integer $k\ge1$.
\item If $Z_i=Q_{n_i}$ with $n_i\geq 2$, then $\cT_i\simeq\OO_{Q_{n_i}}(1)^{\oplus r_i}$
 and the associated partition is $\lambda_i=(k_i^{r_i})$, where either $k_i$ is even or $r_i=n_i$.
\end{enumerate}
In particular, every $\cV_i$ is globally generated.
\end{theorem}

Finally, restricting to smooth weak Fano varieties we obtain the following characterization of products of projective spaces and smooth hyperquadrics with a weaker bundle hypotheses. Recall that a smooth projective variety $X$ is \emph{weak Fano} if $-K_X$ is nef and big. 

\begin{theorem}[=Theorem~\ref{thm:productcharacterization2}]\label{thm:Main3}
Let $X$ be a smooth weak Fano complex variety of dimension $n>0$, and
let $r\geq1$. For $1\leq i\leq r$, let $\cE_i\subseteq T_X$ be a
vector subbundle of rank $n_i\geq1$, let $\cV_i$ be a nef vector
bundle of rank $r_i\geq1$, and let $\lambda_i$ be a nonzero partition
of length $\ell(\lambda_i)\leq\min\{r_i,n_i\}$. Assume that:
\begin{enumerate}
\item the natural morphism
$\cF:=\bigoplus_{i=1}^r\cE_i\rightarrow T_X$ is a subbundle inclusion with integrable image,
\item for every $1\leq i\leq r$, there is an injective morphism
$\theta_i:\bS_{\lambda_i}\cV_i\hookrightarrow
\bS_{\lambda_i}\cE_i$ of coherent $\OO_X$-modules and,
\item the vector bundle $\cA:=\bigotimes_{i=1}^r\cV_i$ is ample on $X$.
\end{enumerate}
Then $\cF=T_X$ and, possibly after reordering, there is an isomorphism
\begin{equation}\label{eq:wf-main-product}
 X\simeq\prod_{i=1}^rZ_i,\qquad
 n=\sum_{i=1}^r n_i,\qquad
 \cE_i\simeq\pr_i^*T_{Z_i}\subseteq T_X,
\end{equation}
where each $Z_i$ is $\PP^{n_i}$ or a smooth hyperquadric
$Q_{n_i}\subset\PP^{n_i+1}$ with $n_i\geq2$. In particular, $X$ is Fano. Moreover, $r_i\leq n_i$ and there are ample vector bundles $\cT_i$ on $Z_i$
such that $\cV_i\simeq\pr_i^*\cT_i$ and $\cT_i$ admits an injective morphism $\bS_{\lambda_i}\cT_i\rightarrow\bS_{\lambda_i}T_{Z_i}$ of coherent sheaves. The bundles $\cT_i$ and the partitions $\lambda_i$ have precisely the
possibilities listed in Theorem~\ref{thm:Main2}\textup{(a)--(d)}.
\end{theorem}

We emphasize that unlike Theorem~\ref{thm:Main2}, we do not assume almost
nefness or individual integrability of the $\cE_i$, even though the $\cV_i$ are assumed to be nef in both results. Here, a subbundle inclusion $\cF:=\bigoplus_{i=1}^r\cE_i\hookrightarrow T_X$ has integrable image precisely when $\cF$ is a regular foliation. This condition
is different from integrability of the individual summands: if every $\cE_i$ is integrable, then $\cF$ is integrable exactly when $[\cE_i,\cE_j]\subseteq\cF$ for all $i\ne j$. Conversely, integrability of $\cF$ need not imply integrability of its summands
in general. 

For smooth toric weak Fano varieties, the foliation methods of
Section~\ref{sec:CharacterizeProduct2} give a characterization in
terms of the Picard number $\rho(X):=\dim_{\RR}N^1(X)_{\RR}$.

\begin{theorem}[=Theorem~\ref{thm:toric-product}]
\label{thm:MainToric}
Let $X$ be a smooth complex projective toric weak Fano variety of dimension $n>0$. Let $\cE_1,\ldots,\cE_r\subseteq T_X$ be
regular foliations of ranks $n_i\geq1$ such that the natural
morphism $\bigoplus_{i=1}^r\cE_i\to T_X$ is a subbundle inclusion.
Then $r\leq\rho(X)$. If $r=\rho(X)$, there is an isomorphism
\[
 X\simeq\prod_{i=1}^r\PP^{n_i},\qquad
 \cE_i=\pr_i^*T_{\PP^{n_i}}\subseteq T_X,
\]
where $\pr_i$ is the $i^{th}$ projection under this isomorphism. Conversely, every product of positive-dimensional projective
spaces admits such a collection of regular foliations of length $r=\rho(X)$.
\end{theorem}

\subsection{Some special instances}\label{Subsec:introspecial}
As applications of the more technical Theorems~\ref{thm:Main2} and \ref{thm:Main3}, we obtain some new characterizations of products of projective spaces and smooth hyperquadrics with easier-to-state/ natural hypotheses. We prove these in Section~\ref{sec:special}.

The first application is a product analogue of Mori's theorem \cite{Mori79}.

\begin{corollary}[=Theorem~\ref{specialcase1proof}]\label{specialcase1}
 Let $X$ be a smooth projective complex variety of dimension $n>0$ and let $r>0$ be an integer. Furthermore let $n_1,\dots,n_r$ be positive integers such that $T_X$ admits a holomorphic splitting $T_X=\bigoplus_{i=1}^r\cE_i$, where each $\cE_i$ is nef vector bundle of  $\rk\cE_i=n_i$ and $\bigotimes_{i=1}^r\cE_i$ is ample. Then, possibly after reordering, there is an isomorphism 
 \[X\simeq \prod_{i=1}^{r}\PP^{n_i} \text{ such that } \cE_i\simeq \pr_i^*T_{\PP^{n_i}},\] 
 where $\pr_i:X\rightarrow \PP^{n_i}$ is the natural projection. 
\end{corollary}

The next result provides a product analogue of the result of Kobayashi-Ochiai \cite{KO73}.

\begin{corollary}[=Theorem~\ref{specialcase2proof}]\label{specialcase2}
Let $X$ be a smooth projective complex variety of dimension $n>0$ and let $r>0$ be an integer. Let $n_1,\dots,n_r$ be positive integers such that $T_X$ admits a holomorphic splitting $T_X=\bigoplus_{i=1}^r\cE_i$ with $\rk\cE_i=n_i$. Furthermore, assume that for each $1\leq i\leq r$, there exists a nef line bundle $\cL_i$ such that $\bigotimes_{i=1}^r\cL_i$ is ample. Then the following statements hold.
\begin{enumerate}[(i)]
\item If $\cM_i:=\det\cE_i\otimes\cL_i^{-n_i}$ is nef for every $1\leq i\leq r$, then, possibly after reordering, there is an isomorphism
\begin{equation}\label{Eq:productspecialcase}
 X\simeq\prod_{i=1}^r Z_i
 \text{ such that }
 \cE_i\simeq\pr_i^*T_{Z_i}
 \text{ and }
 \cL_i\simeq\pr_i^*H_i,
\end{equation}
where $\pr_i:X\rightarrow Z_i$ is the natural projection and for all $1\leq i\leq r$:
\[
 (Z_i,H_i, \cM_i)\simeq(\PP^{n_i},\OO_{\PP^{n_i}}(1), \cL_i)
 \quad\text{or}\quad
 (Z_i,H_i, \cM_i)\simeq(Q_{n_i},\OO_{Q_{n_i}}(1), \OO_X).
\]
\item If $\det\cE_i\otimes\cL_i^{-(n_i+1)}$ is nef for every $1\leq i\leq r$, then \eqref{Eq:productspecialcase} holds such that for all $1\leq i\leq r$:  \[(Z_i,H_i, \det\cE_i)\simeq(\PP^{n_i},\OO_{\PP^{n_i}}(1), \cL_i^{\otimes (n_i+1)}).\]
\item If $\det\cE_i\equiv\cL_i^{n_i}$ for every $1\leq i\leq r$, then \eqref{Eq:productspecialcase} holds such that for all $1\leq i\leq r$:  \[(Z_i,H_i, \det\cE_i)\simeq(Q_{n_i},\OO_{Q_{n_i}}(1), \cL_i^{\otimes n_i}).\] Here $\equiv$ denotes numerical equivalence of line bundles.
\end{enumerate}

\end{corollary}

Finally, we note a characterization of products of quadrics in terms of nondegenerate quadratic tensors on the tangent summands. 

\begin{definition}
A \textit{holomorphic conformal structure} on a holomorphic vector bundle $\cE$ on a compact complex manifold $X$ is an everywhere nondegenerate symmetric bilinear form with values in a line bundle $\cM$, that is, a section in $H^0(X,\Sym^2\cE^*\otimes\cM)$ inducing an isomorphism $\cE\simeq\cE^*\otimes\cM$. Here $\cM$ is called the \textit{conformal line bundle for the form}. For $\cE=T_X$ this is the usual notion of a holomorphic conformal structure on $X$, see \cite{Ye94}.    
\end{definition} 
The following Corollary gives a product analogue of \cite{Ye94} under the assumption that (i) either the conformal line bundles for holomorphic conformal structures on tangent summands admit a square root, or (ii) all the tangent summands have odd rank. 

\begin{corollary}[=Theorem~\ref{specialcase3proof}]\label{specialcase3}
Let $X$ be a smooth projective complex variety of dimension $n>0$ and let $r>0$ be an integer. Let $n_1,\dots,n_r$ be positive integers such that $T_X$ admits a holomorphic splitting $T_X=\bigoplus_{i=1}^r\cE_i$ with $\rk\cE_i=n_i$. Furthermore, assume that for each $1\leq i\leq r$, there exists a $\cL_i^{\otimes 2}$-valued everywhere nondegenerate holomorphic conformal form $q_i\in H^0\left(X,\Sym^2\cE_i\otimes\cL_i^{-2}\right)$ on $\cE_i$ for some nef line bundle $\cL_i$ such that $\bigotimes_{i=1}^r\cL_i$ is ample. Then, possibly after reordering, there is an isomorphism
\[
 X\simeq\prod_{i=1}^rQ_{n_i}
 \text{ such that }
 \cE_i\simeq\pr_i^*T_{Q_{n_i}}
 \text{ and }
 \cL_i\simeq\pr_i^*\OO_{Q_{n_i}}(1),
\]
where $\pr_i:X\rightarrow Q_{n_i}$ is the natural projection. 
If every $n_i$ is odd, the same conclusion holds if each $\cE_i$ admits an everywhere nondegenerate holomorphic conformal form $g_i\in H^0\left(X,\Sym^2\cE_i^*\otimes\cM_i\right)$ for any choice of nef line bundles $\cM_i$ such that $\bigotimes_{i=1}^r\cM_i$ is ample. In this case,
$\cM_i\simeq\pr_i^*\OO_{Q_{n_i}}(2)$.
\end{corollary}

\subsection{Strategy of proof}

To prove Theorem~\ref{thm:Main1}, we first deduce an induced nonzero morphism $\varphi: (\det \cE)^{\otimes q}\rightarrow \cD^{\otimes R}$ from the hypothesis $\bS_\lambda\cE\hookrightarrow\bS_\lambda\cD$. Using filtrations on $\cD^{\otimes R}$ induced by the defining short exact sequence \eqref{eq:unipotent-extension0}, we reduce to the Druel--Paris characterization of projective space and smooth hyperquadrics. We then determine the ample bundle $\cE$ and the admissible partitions $\lambda$ by restricting to lines and identifying the global sections of Schur bundles with orthogonal invariants of Schur modules.

To prove Theorems~\ref{thm:Main2} and \ref{thm:Main3}, we induct on the number of foliations $\cE_1,\dots, \cE_r$. Applying suitable foliation structure theorems to $\cE_1\subseteq T_X$ (up to relabelling), we produce a smooth fibration $f:X\rightarrow Y$ with rationally connected fibers such that the transverse Schur data 
\[\bigl(X;\cE_2, \dots, \cE_r; \cV_2,\dots,\cV_r\bigr) \text{ of Theorem~\ref{thm:Main2} and \ref{thm:Main3}},\] descend to analogous data $\bigl(Y;\cF_2, \dots, \cF_r; \cW_2,\dots,\cW_r\bigr)$ along $f$. We characterize the fibers of $f$ using Theorem~\ref{thm:Main1}. The key obstruction is descending the ampleness conditions Theorem~\ref{thm:Main2}(3) and \ref{thm:Main3}(3) along $f:X\rightarrow Y$. This is because, if $\cU=\bigotimes_{i=2}^{r}\cW_i$, then ampleness of $\cA\simeq\mathcal V_1\otimes f^*\mathcal U$ on $X$ does not automatically imply ampleness of $\cU$ on $Y$. To overcome this, we bound the degrees of quotient line bundles of pullbacks of $\cU$ to curves, from below by a fixed positive multiple of the degree of a fixed ample line bundle on $Y$ (see Proposition~\ref{prop:relative-ample-base}). This uniform bound proves ampleness of $\cU$. Induction then shows that the foliations exhaust $T_X$, and the resulting connection yields the global product decomposition by \cite[Lemma~3.18]{Horing07}.

In the weak Fano case Theorem~\ref{thm:Main3}, applying Druel's foliation theorem \cite[Theorem~1.3]{Druel17} and H\"oring's split-tangent theorem \cite[Theorem~1.1]{Horing26} to the resulting weak Fano fibers, gives integrability of every partial sum of $\bigoplus_{i=1}^{r}\cE_i$. The Bott partial connection then supplies descent of the transverse Schur data of Theorem~\ref{thm:Main3} along $f:X\rightarrow Y$. The same ampleness and splitting arguments therefore let us conclude the proof.

\subsection{Structure of the paper} The paper is organized as follows. Section~\ref{sec:Schur-preliminaries} reviews the Schur-functor preliminaries and computations used in the proofs of the main results. In Section~\ref{sec:CharacterizeP^nQ_n} we prove Theorem~\ref{thm:Main1}. Section~\ref{sec:CharacterizeProduct1} proves Theorem~\ref{thm:Main2} by induction on the number of the regular foliations $\cE_i$ and also proves the descent and positivity results needed for the inductive proof. Finally, Section~\ref{sec:CharacterizeProduct2} proves
Theorem~\ref{thm:Main3} by the same inductive strategy, replacing the
almost-nef descent step by an integrability argument and applying the
common relative ampleness result of
Proposition~\ref{prop:relative-ample-base}. Furthermore, in Section~\ref{sec:CharacterizeProduct2} we also prove Theorem~\ref{thm:MainToric} using the foliation techniques developed there. Finally in Section~\ref{sec:special}, we prove the special characterizations of Section~\ref{Subsec:introspecial} that we obtain as corollaries of the main results.

\subsection*{Acknowledgement}
We thank Michael Zeng for the suggestion of considering Schur functors in the context of an earlier related work, Joaqu\'\i n Moraga for asking the question leading to this work, and also thank him and J\'anos Koll\'ar for helpful comments.

\begingroup
\renewcommand{\addcontentsline}[3]{}
\subsection*{AI Statement} The content of Remark~\ref{Rem:Druel} and the idea of using it to prove Lemma~\ref{lem:calculation} was suggested by AI models. The rest of the article was developed without the use of any AI.
\par
\endgroup

\section{Schur-functor preliminaries}
\label{sec:Schur-preliminaries}

This section reviews the representation-theoretic facts used in the rest of the paper. We work over $\CC$ and the characteristic-zero hypothesis is essential: the group algebra $\CC[\mathfrak S_m]$ of the symmetric group $\mathfrak{S}_m$ is
semisimple, normalized Young symmetrizers are idempotents, and Schur functors are direct
summands of tensor powers.  We refer the reader to \cite{ABW82}, \cite[Lecture~6]{FultonHarris}, \cite[Chapter~1]{Fulton97}, Macdonald \cite[Chapter~I]{Macdonald95}, and Weyman
\cite[Chapter~2]{Weyman03} for further details. 

We set some general conventions, that we will follow throughout the paper. Unless an auxiliary empty partition occurs in a decomposition, partitions used in positivity statements are nonzero. We use the quotient convention for projectivizations. A vector bundle is semiample when a positive tensor power of its tautological line bundle is globally generated. Finally, in numerical positivity arguments, $c_1(\cE)=0$ means equality in $N^1(X)_{\RR}$.

\subsection{Partitions, tableaux, and Schur modules}

A \emph{partition} is a finite nonincreasing sequence $\lambda=(\lambda_1\geq\lambda_2\geq\cdots\geq\lambda_\ell\geq0)$ of non-negative integers. Its degree, $|\lambda|$, length, denoted by $\operatorname{length}(\lambda)$, and effective length, $\length(\lambda)$ are defined by 
\[
 |\lambda|:=\sum_i\lambda_i,
 \qquad
 \operatorname{length}(\lambda):=\ell,
 \qquad
 \length(\lambda):=\# \text{ of positive entries in }\lambda.
\]
When a length $s$ is fixed, we pad $\lambda$ by zeros to
$(\lambda_1,\ldots,\lambda_s)$.  The transpose partition is denoted $\lambda^t$. For a given partition $\lambda\vdash m$, its associated Young diagram is $D_\lambda=\{(i,j)\mid 1\leq j\leq \lambda_i\}$. A Young tableau $T$ of shape $\lambda$ and set of entries $A$ is a bijection $T:D_\lambda\rightarrow A$. Henceforth, unless otherwise mentioned, we will set the labelling set $A=[m]$. A semistandard tableau is a filling of a tableau by entries in an ordered set, with repetitions allowed, whose labels are weakly increasing along rows and strictly increasing down columns. For a box $(i,j)$ of the Young diagram $D_\lambda$, set
\begin{equation}\label{Eq:hooklength}
 c_{i,j}:=j-i,
 \qquad
 h_{i,j}:=\lambda_i-j+\lambda^t_j-i+1,
\end{equation}
where $h_{i,j}$ is the hook length of the box $(i,j)\in D_\lambda$. Let $H_\lambda:=\prod_{(i,j)\in D_\lambda}h_{i,j}$. Choose a tableau $T$ of shape $\lambda$, and let $R_T,C_T\subseteq\mathfrak S_m$ be the row and column stabilizers of $T$, where the symmetric group $\mathfrak{S}_m$ acts on $T$ by acting on $A=[m]$ by permutations. Let
\[
 a_T:=\sum_{\sigma\in R_T}\sigma,
 \quad
 b_T:=\sum_{\tau\in C_T}\operatorname{sgn}(\tau)\tau \ \in \CC[\mathfrak S_m].
\]
By Young's quasi-idempotence theorem, we have $c_T^2=H_\lambda c_T$, where $c_T:=a_Tb_T\in \CC[\mathfrak{S}_m]$. Thus, $e_T:=c_T/H_\lambda$ is an idempotent in $\CC[\mathfrak{S}_m]$. Let $\rho:\mathfrak S_m\rightarrow \GL(V^{\otimes m})$ be the natural permutation left-action on $V^{\otimes m}$ defined by $\rho(\sigma)(v_1\otimes\dots\otimes v_m):=v_{\sigma^{-1}(1)}\otimes\dots\otimes v_{\sigma^{-1}(m)}$. This further induces the $\CC$-algebra homomorphism $\rho:\CC[\mathfrak S_m]\rightarrow \End(V^{\otimes m})$. If $T$ and $T'$ are two distinct Young tableaux for the same shape $\lambda$, then $\im\rho(e_T)$ and $\im\rho(e_{T'})$ are conjugate subspaces of $V^{\otimes m}$. Thus, up to isomorphism, the image of $\rho(e_T):V^{\otimes m}\rightarrow V^{\otimes m}$ is independent of the choice of the Young tableau $T$ and only depends on the shape $\lambda$. Define:
\begin{equation}\label{Eq:Schurmodule}
 \bS_\lambda V:=\im \Big( \rho(e_T):V^{\otimes m}\rightarrow V^{\otimes m} \Big).
\end{equation}
If $\length(\lambda)\leq \dim V$, $\bS_\lambda V$ is the irreducible $\GL(V)$-module of highest weight $\lambda$ by \cite[Theorem~6.3]{FultonHarris}. We will call this the \emph{Schur $\GL(V)$-module of highest weight $\lambda$}. An alternative categorical definition and construction of $\bS_\lambda V$ due to Deligne \cite{Deligne02} is as follows. The isomorphism classes of irreducible representations of $\mathfrak{S}_m$ over $\CC$ are indexed by partitions of $m$. For each partition $\lambda\vdash m$, let $V_\lambda$ be a choice of an irreducible $\mathfrak{S}_m$-representation of the class $\lambda$. Such a concrete choice for each partition $\lambda\vdash m$ is given by Specht modules $V_\lambda=\Sp_\lambda$, see \cite[Theorem~2.4.6]{Sagan01}. Then by \cite[\S1.4]{Deligne02}, the Schur functor $\bS_\lambda:\operatorname{Vec}_\CC\rightarrow \operatorname{Vec}_\CC$ is defined by
\begin{equation}\label{Eq:Schurmodule2}
   \bS_\lambda V:=\Hom_{\mathfrak{S}_m}\left(V_\lambda, V^{\otimes m}\right) 
\end{equation}
Let $n_\lambda:=\dim V_\lambda$ for all partitions $\lambda\vdash m$. By Schur--Weyl duality one has the semi-simple decomposition of $\GL(V)$-modules, with trivial action on $V_\lambda$ (see \cite[Theorem~6.3(2)]{FultonHarris}):
\begin{equation}\label{Eq:SchurWeyl}
 V^{\otimes m}\simeq\bigoplus_{\lambda\vdash m}\bS_\lambda V\otimes V_\lambda\simeq \bigoplus_{\lambda\vdash m}\big(\bS_\lambda V\big)^{\oplus n_\lambda}.
\end{equation}
\begin{example}\label{Ex:1}
Let $V$ be a $\CC$-vector space of dimension $r$ and let $m>0$ be an integer.
\begin{enumerate}
    \item If $\lambda\vdash m$ is the trivial partition $(m)$, then $\bS_{(m)}V\simeq \Sym^mV$.
    \item If $\lambda\vdash m$ is the partition $(1^m):=(1,1,1,\dots, 1)$, then $\bS_{(1^m)}V\simeq \bigwedge^mV$.
    \item If $m=k\cdot r$, and $\lambda\vdash m$ is the partition $(k^r):=(k, k ,\dots, k)$ with all $r$ entries equal to $k$, then $\bS_{(k^r)}V\simeq(\det V)^{\otimes k}$
\end{enumerate}
\end{example}

The construction of Schur modules can be similarly applied to vector bundles, see \cite[\S 1.2.2]{Jiang23}. In what follows, the vector bundles being considered are over smooth projective complex varieties unless otherwise mentioned. The following proposition determines when Schur modules are non-trivial and their functorial properties.

\begin{proposition}\label{prop:Schur-functorial}
Let $\lambda\vdash m$ be a partition and let $\cV$ be a vector bundle of rank $r$ on a scheme $X$ over $\CC$.
\begin{enumerate}[(i)]
\renewcommand{\labelenumi}{\textup{(\roman{enumi})}}
\item  If $\length(\lambda)>r$, then $\bS_\lambda\cV=0$. If $\length(\lambda)\leq r$, then $\bS_\lambda\cV$ is locally free of rank (where $h_{i,j}$ are given by \eqref{Eq:hooklength})
\[ \rk\bS_\lambda\cV= d_\lambda(r):=
 \prod_{(i,j)\in D_\lambda}\frac{r+j-i}{h_{i,j}}.\]
\item  If $\length(\lambda)\leq r$, then $\bS_\lambda\cV$ is a direct summand of $\cV^{\otimes m}$.
\item  For every morphism of schemes $g:Y\rightarrow X$, we have $g^*(\bS_\lambda\cV)\simeq\bS_\lambda(g^*\cV)$.
\item  Every morphism of vector bundles $\phi:\cV\rightarrow\cW$ induces
$\bS_\lambda(\phi):\bS_\lambda\cV\rightarrow\bS_\lambda\cW$.
\end{enumerate}
\end{proposition}

\begin{proof}
    See \cite[Remark~3.12(1)]{Jiang23} for the non-triviality criterion of $\bS_\lambda \cV$ in (i), and \cite[Exercise~6.4]{FultonHarris} for the local-freeness and rank of $\bS_\lambda \cV$. For (ii), refer to \cite[Lemma~4.26]{FultonHarris}. For the functorial properties (iii) and (iv), refer to \cite[Theorem~2.8]{Jiang23}.
\end{proof}

The following lemma demonstrates the permanence of certain positivity conditions on vector bundles, under the Schur functor.
\begin{lemma}\label{lem:schur-positivity}
Let $\lambda\vdash m$ be a partition and $\cE$ be a vector bundle with $\bS_\lambda\cE\neq0$.
\begin{enumerate}[(i)]
\renewcommand{\labelenumi}{\textup{(\roman{enumi})}}
\item If $\cE$ is globally generated, then $\bS_\lambda\cE$ is globally
generated.
\item If $\cE$ is nef, then $\bS_\lambda\cE$ is nef.
\item If $\cE$ is numerically flat, then $\bS_\lambda\cE$ is numerically flat.
\item If $\cE$ is ample, then $\bS_\lambda\cE$ is ample.
\end{enumerate}
\end{lemma}

\begin{proof}
This is a consequence of Proposition~\ref{prop:Schur-functorial}(ii), whereby $\bS_\lambda\cE$ is a quotient of $\cE^{\otimes m}$ and the fact that global generation/ ampleness/ nefness are preserved under tensor powers and quotients and that $\cE$ is numerically flat if and only if $\cE$ and its dual $\cE^\vee$ are both nef. See \cite[Theorem~6.2.12(i), (iv)]{PosII}, \cite[Corollary~6.1.16(ii)]{PosII}, and \cite[Lemma~17.4.3]{stacks-project}.
\end{proof}

\subsection{Littlewood-Richardson decomposition}

For a partition $\lambda\vdash m$, and vector bundles $\cE, \cF$ on $X$, the Schur module $\bS_\lambda(\cE\oplus\cF)$ decomposes into direct sums of tensor products of Schur modules $\bS_\mu\cE\otimes\bS_\nu\cF$, where $\mu, \nu$ are partitions such that $|\mu|+|\nu|=|\lambda|=m$. The multiplicity of $\bS_\mu\cE\otimes\bS_\nu\cF$ in the direct sum decomposition, denoted by $c^\lambda_{\mu, \nu}$, has a combinatorial description due to the Littlewood-Richardson rule \cite[Chapter 1, \S 9]{Macdonald95}. Thus, one has the following non-canonical isomorphism of vector bundles (see \cite[Proposition~1.8]{Deligne02} or \cite[Theorem~2.3.4]{Weyman03}) 
\begin{equation}\label{Eq:LR1}
\bS_\lambda(\cE\oplus\cF)\simeq \bigoplus_{\substack{\mu,\nu\\
             |\mu|+|\nu|=|\lambda|}}\Big(\bS_\mu\cE\otimes\bS_\nu\cF\Big)^{\oplus c^\lambda_{\mu,\nu}}.
             \end{equation}
One can make the isomorphism \eqref{Eq:LR1} canonical as follows. Let $\Sp_\lambda$ denote the irreducible Specht module for $\mathfrak S_m$ corresponding to the partition $\lambda\vdash m$.  If
$\mu\vdash a$, $\nu\vdash b$ are partitions of integers $a,b>0$, and $\lambda\vdash a+b$, then define the
Littlewood--Richardson multiplicity space as
\begin{equation}\label{Eq:LR2}
 \hskip-.25em
 \operatorname{LR}_{\mu,\nu}^{\lambda}
 :=
 \Hom_{\CC[\mathfrak S_{a+b}]}\!\left(
   \Sp_\lambda,\
   \operatorname{Ind}_{\mathfrak S_a\times\mathfrak S_b}^{\mathfrak S_{a+b}}
      (\Sp_\mu\boxtimes\Sp_\nu)
 \right) \simeq
 \Hom_{\CC[\mathfrak S_a\times\mathfrak S_b]}\!\left(
   \operatorname{Res}_{\mathfrak S_a\times\mathfrak S_b}^{\mathfrak S_{a+b}}
      \Sp_\lambda,\
   \Sp_\mu\boxtimes\Sp_\nu
 \right)\!,\hskip-2.25em
\end{equation}
where the isomorphism is due to Frobenius reciprocity. Thus, $\operatorname{LR}^\lambda_{\mu,\nu}$ is a finite dimensional $\CC$-vector space of dimension $\dim_\CC \operatorname{LR}^\lambda_{\mu,\nu}=  c^\lambda_{\mu,\nu}$ by \cite[(7.3) and (7.5)]{Macdonald95}. Consequently, using \eqref{Eq:Schurmodule2} we obtain the following result (see \cite[\S 1.4, 1.5]{Deligne02}), which yields \eqref{Eq:LR1}.
\begin{proposition}\label{prop:LR-direct-sum}
For a partition $\lambda\vdash m$, and vector bundles $\cE, \cF$ on $X$, one has the following canonical isomorphism of vector bundles.
\[
 \bS_\lambda(\cE\oplus\cF)
 \simeq
 \bigoplus_{\substack{\mu,\nu\\
             |\mu|+|\nu|=|\lambda|}}
 \bS_\mu\cE\otimes\bS_\nu\cF\otimes_{\CC}
 \operatorname{LR}_{\mu,\nu}^{\lambda}.
\]
\end{proposition}

\begin{notation}\label{Notation:partitionordering}
There is a natural partial ordering on the set of partitions of a given integer $m$. For any two partitions $\lambda, \mu\vdash m$, we write $\lambda\geq \mu$ if for all $i\geq 1$:
\[\lambda_1+\dots+\lambda_i\geq \mu_1+\dots+\mu_i.\]
\end{notation}

\begin{definition}\label{Def:composition}
    Given an integer $m>0$, a weak composition of $m$ of length $\ell$ is an ordered tuple $\underline{\alpha}:=(\alpha_1,\dots,\alpha_\ell)\in \ZZ^\ell$ such that $\alpha_i\geq 0$ and $\alpha_1+\dots+\alpha_\ell=m$. A weak composition is called a (strong) composition if $\alpha_i>0$ for all $1\leq i\leq \ell$.
\end{definition}

A composition differs from a partition in that a partition is unordered, or equivalently, always ordered as a non-increasing sequence of integers. Thus, we may also write a partition $\lambda\vdash m$ as a tuple $(\lambda_1,\dots,\lambda_\ell)$ where $\lambda_1\geq \lambda_2\geq\dots\geq \lambda_\ell$. For a (weak) composition of $m$ $\underline{\alpha}:=(\alpha_1,\dots,\alpha_\ell)$ and $\lambda\vdash|\underline{\alpha}|:=\alpha_1+\dots+\alpha_\ell=m$ a partition, one can iterate the definition \eqref{Eq:LR2} of Littlewood-Richardson multiplicity space to define 
\begin{equation}
\operatorname{LR}^{\lambda}_{\underline{\alpha}} :=
  \Hom_{\CC[\mathfrak S_m]}\!\left(
    \Sp_{\lambda},
    \operatorname{Ind}_{\mathfrak S_{\underline{\alpha}}}^{\mathfrak S_m} \bigl(\Sp_{(\alpha_1)}\boxtimes\cdots\boxtimes\Sp_{(\alpha_\ell)}
    \bigr)\right),
\end{equation}
where $\mathfrak S_{\underline{\alpha}}:=\mathfrak S_{\alpha_1}\times\cdots\times\mathfrak S_{\alpha_\ell}\subseteq \mathfrak S_m$ is the standard Young subgroup. Since $(\alpha_i)$ is the single row partition of $\alpha_i$, the Specht module $\Sp^{(\alpha_i)}$ is the trivial representation of
$\mathfrak S_{\alpha_i}$ for all $1\leq i\leq \ell$. Thus, we obtain
\begin{equation}\label{Eq:iteratedLR}
\operatorname{LR}^{\lambda}_{\underline{\alpha}}
  =
  \Hom_{\CC[\mathfrak S_m]}\!\left(
    \Sp_{\lambda},
    \operatorname{Ind}_{\mathfrak S_{\underline{\alpha}}}^{\mathfrak S_m}\mathbf 1_{\underline{\alpha}}
  \right),
\end{equation}
where $\mathbf 1_{\underline{\alpha}}$ is the trivial representation of the Young subgroup $\mathfrak S_{\underline{\alpha}}$ using the conventions $\mathfrak S_0=\{1\}$, $(0)=\varnothing$, and $\Sp_{(0)}=\Sp_{\varnothing}:=\CC$. Furthermore, we will define $\operatorname{LR}^{\lambda}_{\underline{\alpha}}=0$ when
$|\underline{\alpha}|:=\sum_i\alpha_i\neq|\lambda|$. Since $\dim_{\CC}
  \operatorname{LR}^{\lambda}_{\underline{\alpha}}$ is the multiplicity of $\Sp_\lambda$ in the $\mathfrak S_m$ representation $\operatorname{Ind}_{\mathfrak S_{\underline{\alpha}}}^{\mathfrak S_m}\mathbf 1_{\underline{\alpha}}$, by Young's rule \cite[Theorem~2.11.2]{Sagan01}, $\dim_{\CC}\operatorname{LR}^{\lambda}_{\underline{\alpha}}=
  K_{\lambda,\underline{\alpha}}$, where the Kostka number $K_{\lambda,\underline{\alpha}}$ is the number of semistandard tableaux of shape $\lambda$ and weight $\underline{\alpha}$. For a weak composition $\underline{\alpha}$ of $m$, let the corresponding partition obtained by arranging the coordinates of $\underline{\alpha}$ in non-increasing order be $\underline{\alpha}^+$. By \cite[Chapter 1, \S 6(6.5) and \S 7 Example 9(b)]{Macdonald95} the Kostka numbers $K_{\lambda,\underline{\alpha}}>0$ if and only if $\lambda\geq\underline{\alpha}^+$ (see Notation~\ref{Notation:partitionordering}).

\begin{corollary}\label{cor:Schur-sums-of-lines}
Let $\lambda\vdash m$ be the partition $\lambda_1\geq \lambda_2\geq\dots\geq\lambda_r\geq 0$, and let $\cE=\cL_1\oplus\cdots\oplus\cL_r$ be an ordered direct sum of line bundles on a scheme $X$ over $\CC$. There is a natural decomposition
\begin{equation}\label{Eq:Schur-sum-of-lines-natural}
\bS_\lambda\cE\simeq
    \bigoplus_{\substack{\alpha_1,\ldots,\alpha_r\geq 0\\
                         \alpha_1+\cdots+\alpha_r=m\\ \lambda\geq\underline{\alpha}^+}}
    \left( \cL_1^{\otimes\alpha_1}\otimes\cdots\otimes\cL_r^{\otimes\alpha_r}
    \right) \otimes_{\CC}\operatorname{LR}^\lambda_{\underline{\alpha}},
\end{equation}
where the sum runs over all weak
compositions $\underline{\alpha}:=(\alpha_1,\ldots,\alpha_r)$ of $m$ of length $r$ such that $\lambda\geq\underline{\alpha}^+$, where $\underline{\alpha}^+\vdash m$ is the partition associated to $\alpha$. 

If $X=\PP^1$, and $\cE\simeq\bigoplus_{i=1}^{r}\mathcal{O}_{\PP^1}(d_i)$ with $d_1\geq\cdots\geq d_r$ then the maximum and minimum degrees of the line-bundle summands in $\bS_\lambda\cE$ are
\begin{align}
\label{eq:Schur-split-dmaxmin}
    d_{\max}=\sum_{i=1}^{r}\lambda_i d_i,\qquad
    d_{\min}=\sum_{i=1}^{r}\lambda_i d_{r+1-i}.
\end{align}
\end{corollary}

\begin{proof}
Inducting Proposition~\ref{prop:LR-direct-sum} on the number of summands and applying to $\cE=\cL_1\oplus\ldots\oplus\cL_r$ yields,
\begin{align*}
    \bS_\lambda\cE\simeq \bS_\lambda\big(\bigoplus_{i=1}^{r}\cL_i\big)\simeq \bigoplus_{\substack{\alpha_1,\ldots,\alpha_r\geq 0\\
                         \alpha_1+\cdots+\alpha_r=m}}
    \left( \bS_{(\alpha_1)}\cL_1\otimes\cdots\otimes\bS_{(\alpha_r)}\cL_r
    \right) \otimes_{\CC}\operatorname{LR}^\lambda_{\underline{\alpha}}
\end{align*}

We therefore obtain \eqref{Eq:Schur-sum-of-lines-natural} since $\bS_{(\alpha_i)}\cL_i\simeq \Sym^{\alpha_i}\cL_i\simeq\cL_i^{\otimes\alpha_i}$, since $\cL_i$ are line bundles.  The direct sum is further restricted to those weak compositions $\underline{\alpha}$ such that $\lambda\geq\underline{\alpha}^+$ since $\dim\operatorname{LR}^\lambda_{\underline{\alpha}}= K_{\lambda,\underline{\alpha}}=0$ otherwise. If $X=\PP^1$ and $\cE\simeq\bigoplus_{i=1}^{r}\mathcal{O}_{\PP^1}(d_i)$ with $d_1\geq\cdots\geq d_r$, then the degree of the summand indexed by the weak composition $\underline{\alpha}$ is the inner product $\underline{d}\cdot\underline{\alpha}$, where $\underline{d}:=(d_1,\dots, d_r)$. For a vector $\alpha=(\alpha_1,\dots, \alpha_r)$, let $\alpha_{rev}:=(\alpha_r,\alpha_{r-1},\dots, \alpha_1)$. Then among the non-zero summands of $\bS_\lambda(\bigoplus_{i=1}^{r}\OO_{\PP^{1}}(d_i))$, we have
\[\underline{d}\cdot\lambda_{rev}\leq \underline{d}\cdot\underline{\alpha}^+_{rev}\leq\underline{d}\cdot\alpha\leq\underline{d}\cdot\underline{\alpha}^+\leq \underline{d}\cdot\lambda\]
The second and third inequalities above are the rearrangement inequalities, while the first and fourth inequalities are due to $\lambda\geq\underline{\alpha}^+$ to ensure $\operatorname{LR}^\lambda_{\underline{\alpha}}\neq 0$. This yields \eqref{eq:Schur-split-dmaxmin}.

\end{proof}

\subsection{Determinants and injective morphisms of Schur bundles}

\begin{lemma}\label{lem:Schur-basic}
Let $\lambda\vdash m$ be a partition, $\cE$ a vector bundle of rank $r$ and assume $\length(\lambda)\leq r$. Let $d:=d_\lambda(r)=\rk(\bS_\lambda\cE)$ and $q_\lambda(r):=d_{\lambda}(r)m/r$. Then $q_\lambda(r)\in\ZZ$, and:
\begin{enumerate}[(i)]
\renewcommand{\labelenumi}{\textup{(\roman{enumi})}}
\item there is a natural isomorphism $\det(\bS_\lambda\cE) \simeq (\det\cE)^{\otimes q_\lambda(r)}$,
\item for every line bundle $\cL$, there is a natural isomorphism $\bS_\lambda(\cE\otimes\cL) \simeq
 \bS_\lambda\cE\otimes\cL^{\otimes m}$,

\end{enumerate}
\end{lemma}

\begin{proof}
Note that $\det(\bS_\lambda(-))$ is a one-dimensional representation of $\GL_r$ and thus is of the form $\det(-)^{\otimes q}$ for some integer $q$. To compute $q$, it suffices to work fiberwise. To this end, let $V=\CC^r$ and $\rho_\lambda:\GL(V)\rightarrow\GL(\bS_\lambda V)$ be the Schur representation, whereby $\det(\bS_\lambda V)=\det_{\bS_\lambda V}\circ\rho_\lambda$.  Since $\bS_\lambda$ is homogeneous polynomial functor of degree $m$, the scalar matrix $tI_r$ of rank $r$ acts on $\bS_\lambda V$ as $t^mI_d$, where $d=d_\lambda(r)=\rk(\bS_\lambda\cE)$. Consequently,
\[ t^{rq}=\det(tI_r)^q=\det_{\bS_\lambda V}\!(\rho_\lambda(tI_r))=\det_{\bS_\lambda V}\!(t^mI_d)=t^{md},
\]
and hence $q=q_\lambda(r)=d_\lambda(r)m/r\in\ZZ$. This proves (i). For (ii) note that $\bS_\lambda(\cE\otimes \cL):=\Hom_{\mathfrak{S}_m}(\Sp_\lambda, \cE^{\otimes m}\otimes\cL^{\otimes m})\simeq \Hom_{\mathfrak{S}_m}(\Sp_\lambda, \cE^{\otimes m})\otimes\cL^{\otimes m}=:\bS_\lambda(\cE)\otimes \cL^{\otimes m}$, since $\mathfrak{S}_m$ acts by identity on $\cL^{\otimes m}$. 
\end{proof}

In this section we will note some properties and consequences of an injective morphism of Schur bundles $\bS_\lambda\cE\rightarrow \bS_\lambda\cV$, for any given partition $\lambda$ and vector bundles $\cE, \cV$ on a smooth projective complex variety $X$. Furthermore, we define a \emph{Schur-subbundle inclusion} to be an injective morphism $\theta:\bS_\lambda\cE\rightarrow\bS_\lambda\cV$ of Schur bundles with locally free cokernel. Note that a morphism of Schur bundles $\bS_\lambda\cE\rightarrow \bS_\lambda\cV$ need not be induced by a morphism $\cE\rightarrow \cV$ of vector bundles.

\begin{lemma}\label{lem:Schur-rank-comparison}
If $\lambda\vdash m$ is a non-empty partition and $\bS_\lambda\cE\rightarrow\bS_\lambda\cV$
is an injective morphism of Schur bundles, with $\rk\cE=r$, $\rk\cV=s$, and $\length(\lambda)\leq\min\{r,s\}$, then $r\leq s$. If $\rk \bS_\lambda\cV=1$ then $r=s$ and $\lambda=(k^s)$ for some integer $k$.
\end{lemma}

\begin{proof}
By generic injectivity, we have  $d_\lambda(r):=\rk\bS_\lambda\cE\leq \rk \bS_\lambda\cV=:d_\lambda(s)$. By the explicit formula in Proposition~\ref{prop:Schur-functorial}(i), $d_\lambda(-)$ is an increasing function, which proves $r\leq s$. If $\rk \bS_\lambda\cV=1$, then $d_\lambda(s)=1$. By \cite[Theorem~6.3(i)]{FultonHarris},
\[d_\lambda(s)=\prod_{1\leq i<j\leq s}\frac{\lambda_i-\lambda_j+j-i}{j-i}= \prod_{1\leq i<j\leq s}\Big(1+\frac{\lambda_i-\lambda_j}{j-i}\Big),\]
from which it is clear that each of the factors in the product is at least one since $\lambda\vdash m$ is a partition, whereby $\lambda_i\geq\lambda_j$ for $i<j$. Thus $d_\lambda(s)=1$ if and only if for all $1\leq i\leq s$, $\lambda_i=k$ for some integer $k$, i.e. $\lambda=(k^s)$. But $d_\lambda(r)\leq d_\lambda(s)=1$, implies $d_\lambda(r)=1$, whereby the same argument yields $\lambda=(p^r)$ for some positive integer $p$. Thus the number of nonzero entries in $\lambda$ is equal to $s$ and also equal to $r$, yielding $r=s$. 
\end{proof}

\begin{proposition}\label{prop:Schur-det-to-tensor}
Let $\lambda\vdash m$ be a partition, and $\cE$ be a
vector bundle of rank $r\geq\ell(\lambda)$. Let $\cV$ be a vector bundle admitting an injective map $\theta:\bS_\lambda\cE\rightarrow\bS_\lambda\cV$. If  $d=d_\lambda(r)$, $q=q_\lambda(r)$ as defined in Lemma~\ref{lem:Schur-basic} and $R=R_\lambda(r):=md_\lambda(r)=rq_\lambda(r)$ then there is a nonzero morphism
\[
\vartheta: (\det\cE)^{\otimes q}\rightarrow\cV^{\otimes R}.
\]
\end{proposition}

\begin{proof}
Since $d=\rk \bS_\lambda\cE$, taking $d^{th}$ exterior power of the injective map $\theta:\bS_\lambda\cE\rightarrow\bS_\lambda\cV$ induces a non-zero map
\[
 \bigwedge^d\theta:\det(\bS_\lambda\cE)\rightarrow\bigwedge^d(\bS_\lambda\cV).
\]
In characteristic $0$, there is a split injection $\iota_1:\bigwedge^d(\bS_\lambda\cV)\hookrightarrow(\bS_\lambda\cV)^{\otimes d}$ and a split injection $\iota_2:\bS_\lambda\cV\hookrightarrow\cV^{\otimes m}$ by Proposition~\ref{prop:Schur-functorial}(ii). Thus, we obtain the following non-zero composition:
\[
 \vartheta:\det(\bS_\lambda\cE)\xrightarrow{\bigwedge^d\theta}\bigwedge^d(\bS_\lambda\cV)\xrightarrow{\iota_1}(\bS_\lambda\cV)^{\otimes d}\xrightarrow{\iota_2^{\otimes d}} \cV^{\otimes md}=\cV^{\otimes R}.
\]
Finally, $\det(\bS_\lambda\cE)\simeq(\det\cE)^{\otimes q}$ by Lemma~\ref{lem:Schur-basic}(i), whereby the composite $\vartheta$ is the desired non-zero map.
\end{proof}

\subsection{Special Orthogonal invariants of Schur modules}

For a given partition $\lambda\vdash m$, the invariants of the natural $\SO_r$-actions defined by restriction on the $\GL_r$-representation $\bS_\lambda\CC^r$ are related to the global sections of the bundle $\bS_\lambda( T_{Q_n}(-1))\simeq\bS_\lambda T_{Q_n}\otimes\OO_{Q_n}(-m)$, where $T_{Q_n}$ is the tangent bundle of the $n$-dimensional hyperquadric $Q_n\subseteq\PP^{n+1}$. We record this relation in this section.
\begin{lemma}\label{lem:SO-invariants}
Let $r\geq2$ and $\lambda\vdash m$ be a partition of length $\ell(\lambda)\leq r$, such that $\lambda:=(\lambda_1,\ldots,\lambda_r)$, after possibly adding zeroes to the right. Let $\SO_r$ and $\operatorname{O}_r$ be the special orthogonal and orthogonal groups on $\CC^r$ respectively. Then
\begin{align*}
& \dim(\bS_\lambda\CC^r)^{\SO_r}
 =
 \begin{cases}
 1,&\lambda_1\equiv\cdots\equiv\lambda_r\pmod 2,\\
 0,&\text{otherwise}.
 \end{cases}\\
 &\dim(\bS_\lambda\CC^r)^{\operatorname{O}_r}
 =
 \begin{cases}
 1,&\lambda_1\equiv\cdots\equiv\lambda_r\equiv0\pmod 2,\\
 0,&\text{otherwise}.
 \end{cases}
\end{align*}
\end{lemma}
\begin{proof}
    This is \cite[Corollary~3.3, 3.6]{DrenskyHristova22}. An alternative reference is \cite[Pg.~409]{Procesi07}.
\end{proof}

Now let $Q_n\subseteq \PP(W)\simeq\PP^{n+1}$ be a smooth $n$-dimensional quadric with $n\geq 2$. Fix a point $x\in Q_n$, and equip the fiber $V_x:=( T_{Q_n}(-1))_x$ of the twisted tangent bundle $T_{Q_n}(-1)$ with the nondegenerate quadratic form induced by the quadratic form $q$ on $W\cong \CC^{n+2}$ defining $Q_n$. Choosing a nonzero vector $w$ representing $x=[w]\in Q_n$, identifies $(T_{Q_n}(-1))_x\cong w^\perp/\langle w\rangle=:V$ and $q$ induces a nondegenerate quadratic form on $V$.  If $G:=\SO(W,q)$ and $P:=\operatorname{Stab}_G(x)$ is the parabolic subgroup corresponding to $x$, then $Q_n\simeq G/P$.

\begin{lemma}\label{lem:quadric-Schur-sections}
With the above notations, we have the isomorphism  $$H^0\bigl(Q_n,\bS_\lambda (T_{Q_n}(-1))\bigr)
  \xrightarrow{\sim}
  \bigl(\bS_\lambda V\bigr)^{
    \SO(V)}.$$
When $ H^0(Q_n,\bS_\lambda( T_{Q_n}(-1)))\neq 0$, it is generated by a nowhere vanishing section, spanning a trivial homogeneous direct summand $\OO_{Q_n}\subseteq\bS_\lambda( T_{Q_n}(-1))$.
\end{lemma}

\begin{proof}
We first treat $n\geq3$. The homogeneous description in
\cite[Lemma~2 and its proof]{DruelParis13} identifies $T_{Q_n}(-1)\simeq G\times^P V$, where the unipotent radical of $P$ acts trivially on $V$ and the Levi action factors through the standard representation of $\SO(V)$. The proof of \cite[Lemma~2]{DruelParis13} gives, by evaluation at $x\in Q_n$, an isomorphism
\begin{equation}\label{eq:DP-tensor-evaluation}
  \operatorname{ev}_x:
  H^0\bigl(Q_n, T_{Q_n}(-1)^{\otimes m}\bigr)
  \xrightarrow{\sim}
  \bigl(V^{\otimes m}\bigr)^{\SO(V)}.
\end{equation}
We use this argument only for $n\geq3$ since the proof in
\cite[Lemma~2]{DruelParis13} uses the irreducibility of the standard
$\SO(V)$-module and the absence of nontrivial characters of $\SO(V)$, which hold only for $n\geq 3$. 

First, the isomorphism \eqref{eq:DP-tensor-evaluation} is $\mathfrak S_m$-equivariant.  Indeed, $\mathfrak S_m$ acts on both
sides by permuting tensor positions, and for every
$\sigma\in\mathfrak S_m$ and every section $s$ one has $\operatorname{ev}_x(\sigma\cdot s)=\sigma\cdot\operatorname{ev}_x(s)$. Let $e_T\in\CC[\mathfrak S_m]$ be the normalized Young idempotent defining $\bS_\lambda$ in \eqref{Eq:Schurmodule}. The diagonal $\SO(V)$-action on $V^{\otimes m}$
commutes with the permutation action of $\mathfrak S_m$.  Explicitly,
for $g\in\SO(V)$ and $\sigma\in\mathfrak S_m$,
\begin{align}
 g^{\otimes m}\sigma
 (v_1\otimes\cdots\otimes v_m)
 = gv_{\sigma^{-1}(1)}\otimes\cdots\otimes
 gv_{\sigma^{-1}(m)}= \sigma g^{\otimes m}
 (v_1\otimes\cdots\otimes v_m).
\end{align}
Hence $e_T$, being a linear combination of permutations, commutes
with $\SO(V)$. Since $e_T$ is idempotent, we obtain the splitting 
\[
   T_{Q_n}(-1)^{\otimes m}
  =
  \operatorname{im}(e_T)\oplus\ker(e_T)=\bS_\lambda(T_{Q_n}(-1))\oplus\ker(e_T).
\]
Applying global sections to this direct-sum decomposition gives
\begin{equation}\label{eq:H0-image-idempotent}
  H^0\bigl(Q_n,\bS_\lambda(T_{Q_n}(-1))\bigr)
  =
  \operatorname{im}
  \left(
    H^0(e_T):
    H^0(Q_n, T_{Q_n}(-1)^{\otimes m})
    \rightarrow
    H^0(Q_n, T_{Q_n}(-1)^{\otimes m})
  \right).
\end{equation}
Equivalently, $H^0\bigl(Q_n,\bS_\lambda (T_{Q_n}(-1))\bigr)\simeq e_T H^0\bigl(Q_n, T_{Q_n}(-1)^{\otimes m}\bigr)\simeq e_T\bigl(V^{\otimes m}\bigr)^{\SO(V)}$, where the second equality uses \eqref{eq:DP-tensor-evaluation}. Since $e_T$ commutes with the $\SO(V)$-action, we obtain the desired isomorphism $H^0\bigl(Q_n,\bS_\lambda (T_{Q_n}(-1))\bigr)\simeq  \bigl(e_TV^{\otimes m}\bigr)^{\SO(V)}\simeq \bigl(\bS_\lambda V\bigr)^{\SO(V)}$. By Lemma~\ref{lem:SO-invariants}, this vector space is either zero or one-dimensional. Assume it is non-zero, whereby $L:=\bigl(\bS_\lambda V\bigr)^{\SO(V)}\simeq\CC$. Since $\SO(V)$ is reductive, there is an
$\SO(V)$-stable complement $M$ such that $\bS_\lambda V=L\oplus M$.
The $P$-action on $\bS_\lambda V$ factors through $\SO(V)$, so this is also a decomposition of $P$-modules as well.  The group $P$ acts trivially on $L$, and therefore
passing to associated homogeneous bundles gives
\[
\begin{aligned}
 \bS_\lambda (T_{Q_n}(-1))
 \simeq
 G\times^P\bS_\lambda V
 \simeq
 \bigl(G\times^P L\bigr)
 \oplus
 \bigl(G\times^P M\bigr)\simeq
 \OO_{Q_n}\oplus\bigl(G\times^P M\bigr),
\end{aligned}
\]
which gives the desired direct sum decomposition of $ \bS_\lambda T_{Q_n}(-1)$ with a trivial summand (the functor $E\mapsto G\times^P E$ sending a rational $P$ module to a homogeneous bundle on $G/P$ is exact, compatible
with direct sums and tensor products and sends a trivial $P$-module $L$
to the trivial bundle $G\times^P L\simeq\OO_{G/P}\otimes_{\CC}L$-- see \cite[I.5.5, I.5.10--I.5.12, and I.5.15]{Jan03}.) Furthermore, $H^0(Q_n, \bS_\lambda (T_{Q_n}(-1)))\simeq H^0(Q_n, \OO_{Q_n})\oplus H^0(Q_n, (G\times^P M))\simeq H^0(Q_n, \OO_{Q_n})$ since $\dim_\CC H^0(Q_n, \bS_\lambda (T_{Q_n}(-1)))=1$.

It remains to treat $n=2$. Any smooth hyperquadric $Q_2\subseteq\PP^3$ is isomorphic to $\PP^1\times\PP^1$, and thus, $\OO_{Q_2}(1)=\OO_{\PP^1\times\PP^1}(1,1)$. Furthermore, $T_{Q_2}\simeq
  \OO(2,0)\oplus\OO(0,2)$,
and hence
\[
   T_{Q_2}(-1)  \simeq
  \OO(1,-1)\oplus\OO(-1,1).
\]
Since $\rk T_{Q_2}(-1)=2$, if $\ell(\lambda)>2$, then $\bS_\lambda (T_{Q_2}(-1))=0$. Thus, let $\lambda=(a,b)$ with $a\geq b\geq0$ and set $r:=a-b$.  For every rank-two vector bundle $\cE$ one has $\bS_{(a,b)}\cE \simeq\Sym^{a-b}\cE\otimes(\det\cE)^{\otimes b}$ as a special case of the determinant shift identity for Schur functors.  Since
$\det T_{Q_2}(-1)\simeq\OO_{Q_2}$, it follows that
\begin{align}\label{Eq:symn2}
 \bS_{(a,b)} T_{Q_2}(-1)\simeq
 \operatorname{Sym}^{r}
 \bigl(\OO(1,-1)\oplus\OO(-1,1)\bigr)\simeq
 \bigoplus_{j=0}^{r}
 \OO(r-2j,-r+2j).
\end{align}
Since $H^0\bigl(\PP^1\times\PP^1,\OO(c,-c)\bigr)\neq0$ if and only if $c=0$, by \eqref{Eq:symn2} we obtain
\begin{equation}\label{Eq:H0n2}
  H^0\bigl(Q_2,\bS_{(a,b)} (T_{Q_2}(-1))\bigr)
  \simeq
  \begin{cases}
    \CC,&a-b\ \text{is even},\\
    0,&a-b\ \text{is odd}.
  \end{cases}
\end{equation}
In the even case, the unique contributing summand to \eqref{Eq:H0n2} is
$\OO_{Q_2}$, corresponding to $j=r/2$ in \eqref{Eq:symn2}. On the representation-theoretic side, $\SO_2(\CC)\simeq\mathbb G_m$ and after choosing an isotropic basis, its standard representation is $\CC^2\simeq\chi\oplus\chi^{-1}$, where $\chi$ is the weight one character of $\SO_2(\CC)\simeq\mathbb{G}_m$. Since the determinant character is trivial on $\SO_2$, we have
\[
  \bS_{(a,b)}\CC^2
  \simeq
  \operatorname{Sym}^{r}
  (\chi\oplus\chi^{-1})\simeq\bigoplus_{j=0}^{r}\chi^{\otimes (r-2j)}
\]
as an $\SO_2$-module, where $r=a-b$. Thus, its weights are precisely $\{r-2j\mid 0\leq j\leq r\}$, each with multiplicity one.  Hence $\bS_{(a,b)}\CC^2$ contains a weight zero summand exactly when $r=a-b$ is even, in which case the invariant line (i.e., the weight zero piece) is unique. This yields the desired isomorphism
\[
  H^0\bigl(Q_2,\bS_{(a,b)} T_{Q_2}(-1)\bigr)
  \simeq
  \bigl(\bS_{(a,b)}\CC^2\bigr)^{
    \SO_2}.
\]
When $ H^0\bigl(Q_2,\bS_{(a,b)} T_{Q_2}(-1)\bigr)\neq 0$, evaluation identifies the constant section of the unique trivial summand $\OO_{Q_2}$ with the unique invariant
vector in the fiber.  This section is nowhere vanishing and proves
the final assertion for $Q_2$.
\end{proof}

\section{Schur functor Characterization of $\PP^n$ and $Q_n$}\label{sec:CharacterizeP^nQ_n}

Throughout this section, $Q_n$ will denote the smooth hyperquadric of dimension $n$ in $\PP^{n+1}$. The main result of this section is the following generalization of \cite[Theorem~1.1, 1.2]{Ghosh26}. This shall be useful in the proof of the inductive step of Theorems~\ref{thm:Main2} and \ref{thm:Main3} in the later sections.

\begin{definition}
    A unipotent vector bundle on a smooth projective complex variety $X$ is a vector bundle $\cU$ admitting a finite filtration by vector subbundles with trivial successive quotients, i.e, \begin{equation}\label{Eq:unipotentfiltration}
 0=\cU_0\subseteq\cU_1\subseteq\cdots\subseteq\cU_e=\cU
 \qquad\text{with}\qquad
 \cU_i/\cU_{i-1}\simeq\OO_X^{\oplus g_i}
\end{equation}
for some integers $g_i\geq1$. By convention, we will also consider the zero sheaf $\mathbf{0}_X$ as a unipotent vector bundle of rank $0$.
\end{definition}

\begin{theorem}\label{thm:generalizedKovacs}
Let $X$ be a smooth projective complex variety of dimension $n>0$, let
$\cU$ be a unipotent vector bundle on $X$ and let $\cD$ fit into an exact sequence
\begin{equation}\label{eq:unipotent-extension}
 0\rightarrow T_X\rightarrow\cD\rightarrow\cU\rightarrow0.
\end{equation}
Let $\cE$ be an ample vector bundle of rank $r$ on $X$ and $\lambda\vdash m$ be a nonzero partition of $\length(\lambda)\leq r$, and suppose that there is an injective morphism $\bS_\lambda\cE\hookrightarrow\bS_\lambda\cD$ of coherent $\OO_X$-modules. Then either: 
\begin{enumerate}
    \item[(a)] $(X, \det\cE)\simeq(\PP^n, \OO_{\PP^n}(l))$ with $r\leq l\leq \frac{r(n+1)}{n}$ or,
    \item[(b)] $(X, \cE)\simeq (Q_n, \OO_{Q_n}(1)^{\oplus r})$ with $n\ge2$, $r\leq n$. Furthermore, the morphism $\bS_\lambda\cE\hookrightarrow\bS_\lambda\cD$ factors as $\bS_\lambda\cE\hookrightarrow\bS_\lambda T_{Q_n}\hookrightarrow\bS_\lambda\cD$ and the partition $\lambda\vdash m$ is $(k^r)$, where either $k$ is even or $r=n$.
\end{enumerate}
Furthermore, in case (a), the bundle $\cE$ on $\PP^n$ is determined as:
\begin{enumerate}
 \item[(i)] if $r<n$, then $l=r$ and
       $\cE\simeq\OO_{\PP^n}(1)^{\oplus r}$;
 \item[(ii)] if $r=n$ and $l=n$, then
       $\cE\simeq\OO_{\PP^n}(1)^{\oplus n}$;
 \item[(iii)] if $r=n$ and $l=n+1$, then $\cE\simeq T_{\PP^n}$ or $\cE\simeq \OO_{\PP^n}(2)\oplus\OO_{\PP^n}(1)^{\oplus(n-1)}$. The second alternative can occur only when
       $\lambda=(k^n)$ for some $k\ge1$;
 \item[(iv)] if $r>n$, then $l=r$ and $\cE\simeq\OO_{\PP^n}(1)^{\oplus r}$.
\end{enumerate}
\end{theorem}

To prove Theorem~\ref{thm:generalizedKovacs}, we first characterize the variety $X$ and fully determine the quadric case, i.e., prove Theorem~\ref{thm:generalizedKovacs}(b). We then characterize $\cE$ when $X\simeq \PP^n$. 

\subsection{Characterizing the variety and the quadric case}

\begin{theorem}\label{thm:stabilized-schur-kovacs}
Let $X$ be a smooth projective complex variety of dimension $n>0$, let $\cU$ be a unipotent vector bundle of rank $g$ on $X$, i.e., a vector bundle admitting a filtration \eqref{Eq:unipotentfiltration} by vector subbundles, and let $\cD$ be a vector bundle on $X$ which fits into an exact sequence
\begin{equation}\label{eq:stabilized-extension}
 0\rightarrow T_X\rightarrow\cD
 \rightarrow\cU\rightarrow 0.
\end{equation}
Let $\cE$ be an ample vector bundle on $X$ of rank $r$ and $\lambda$ be a nonzero partition of $m:=|\lambda|$ with length $s:=\length(\lambda)\leq r$. Suppose that there is an injective map of coherent sheaves
\begin{equation}\label{Eq:geninjectivecharacterization}
 \bS_\lambda\cE\rightarrow\bS_\lambda\cD.
\end{equation}
Then $(X, \det\cE)\simeq(\PP^n, \OO_{\PP^n}(l))$ with $r\leq l\leq \frac{r(n+1)}{n}$, or $n\geq2$ and $(X, \cE)\simeq (Q_n, \OO_{Q_n}(1)^{\oplus r})$ with $r\leq n$. Furthermore, if $X\simeq Q_n$ then the map \eqref{Eq:geninjectivecharacterization} factors through $\bS_\lambda T_{Q_n}$ by an injective map. Moreover, in this case $\lambda=(k^r)$ is a rectangular partition with either $k$ even or $r=n$.
\end{theorem}

\begin{proof}
Let $m=|\lambda|$, $d=d_\lambda(r)=\rk(\bS_\lambda\cE)$, and $q=q_\lambda(r):=d_\lambda(r)m/r$ as introduced in Lemma~\ref{lem:Schur-basic}. Then by Lemma~\ref{lem:Schur-basic} we have $md=rq$, and we denote the common value by $R:=R_\lambda(r)$.
By Proposition~\ref{prop:Schur-det-to-tensor}, the injective Schur morphism $\bS_\lambda\cE\rightarrow\bS_\lambda\cD$ induces a nonzero morphism $\varphi: (\det \cE)^{\otimes q}\rightarrow \cD^{\otimes R}$. Equip $\cD$ with the two-step descending filtration $F^0\cD\supseteq F^1\cD\supseteq F^2\cD$ induced with $F^0\cD=\cD$, $F^1\cD=T_X$, $F^2\cD=0$. By the tensor-product filtration for finite locally free sheaves \cite[\S~4.1.1]{LaPorta26}, the tensor power $\cD^{\otimes R}$ carries an induced finite descending filtration $F^\bullet$ with
\begin{equation}\label{Eq:gradeddecomposition}
 \operatorname{gr}_F^t(\cD^{\otimes R})
 \simeq
 \bigoplus_{\substack{\varepsilon_1,\ldots,\varepsilon_R\in\{0,1\}\\
                       \varepsilon_1+\cdots+\varepsilon_R=t}}
 \bigotimes_{i=1}^R\operatorname{gr}_F^{\varepsilon_i}(\cD)
 \simeq
 \bigoplus_{\substack{I\subseteq\{1,\ldots,R\}\\ |I|=t}}
 T_X^{\otimes t}\otimes \cU^{\otimes(R-t)}.
\end{equation}
The unipotent filtration
\[
 0=\cU_0\subseteq\cU_1\subseteq\cdots\subseteq\cU_e=\cU,
 \qquad
 \cU_i/\cU_{i-1}\simeq\OO_X^{\oplus g_i},
\]
and the same tensor-product-filtration construction give, for every $N\geq0$, a finite filtration of $\cU^{\otimes N}$ whose associated graded pieces are trivial vector bundles. Consequently, for every $t$, the vector bundle $\operatorname{gr}_F^t(\cD^{\otimes R})$ in \eqref{Eq:gradeddecomposition} carries a finite descending filtration, say $\mathscr H_t^\bullet$, whose associated graded pieces are finite direct sums of copies of $T_X^{\otimes t}$. Choose $t$ maximal such that $\varphi: (\det \cE)^{\otimes q}\rightarrow \cD^{\otimes R}$ factors through $F^t\cD^{\otimes R}$. The induced morphism $\psi_t:(\det\cE)^{\otimes q}\rightarrow
 \operatorname{gr}_F^t(\cD^{\otimes R})$ is nonzero, since otherwise $\varphi$ would factor through
$F^{t+1}\cD^{\otimes R}$, contradicting maximality of chosen $t$. Choose $j$ maximal such that $\psi_t$ factors through $\mathscr H_t^j$. Analogously the induced morphism 
\[
 (\det\cE)^{\otimes q}\rightarrow \operatorname{gr}_{\mathscr H_t}^j
 \bigl(\operatorname{gr}_F^t(\cD^{\otimes R})\bigr)
\]
is nonzero; otherwise $\psi_t$ would factor through $\mathscr H_t^{j+1}$, contradicting maximality of chosen $j$. Since $\operatorname{gr}_{\mathscr H_t}^j
 \bigl(\operatorname{gr}_F^t(\cD^{\otimes R})\bigr)$ is a finite direct sum of copies of $T_X^{\otimes t}$, projecting to a suitable direct summand therefore gives a nonzero morphism $(\det\cE)^{\otimes q}\rightarrow T_X^{\otimes t}$ or equivalently,
\begin{equation}\label{eq:Schur-tangent-tensor}
 H^0\!\left(X,T_X^{\otimes t}\otimes(\det\cE)^{-\otimes q}\right)\neq0.
\end{equation}
If $t=0$, the relevant $\mathscr H_0$-graded piece is a trivial vector bundle, so the preceding nonzero morphism would give a nonzero section of $(\det\cE)^{-\otimes q}$, which is impossible because $\det\cE$ is ample and
$\dim X>0$.  Hence $1\leq t\leq R$, since $F^{R+1}\cD^{\otimes R}=0$. 

Now let $\cH:=\cE^{\oplus q}$, which has rank $\rk\cH=rq=R$ and has determinant $\det\cH=(\det\cE)^{\otimes q}$. Thus, \eqref{eq:Schur-tangent-tensor} is equivalent to $ H^0\!\left(X,T_X^{\otimes t}\otimes(\det\cH)^{-\otimes 1}\right)\neq0$. Then we have two cases as follows.
\begin{enumerate}
    \item If $t<R=\rk\cH$, then write $R=t+k$ for some $k\geq 1$. By \cite[Theorem A]{DruelParis13}, we have $(X, \det\cH)\simeq(\PP^n, \OO_{\PP^n}(l))$, with $R\leq l\leq t(\frac{n+1}{n})$. But since $\det \cH=(\det\cE)^{\otimes q}$ and $R=rq$, we have $(X, \det\cE)\simeq(\PP^n, \OO_{\PP^n}(l'))$, with $r\leq l'\leq t(\frac{n+1}{nq})\leq r(\frac{n+1}{n})$.
    \item  If $t=R=\rk\cH$, then \cite[Theorem B]{DruelParis13} yields either $(X, \det\cH)\simeq(\PP^n, \OO_{\PP^n}(l))$, with $R\leq l\leq R(\frac{n+1}{n})$ or $(X, \cH)\simeq(Q_n, \OO_{Q_n}(1)^{\oplus R})$. Since $\cH:=\cE^{\oplus q}$, we have $(X, \det\cE)\simeq(\PP^n, \OO_{\PP^n}(l'))$, with $r\leq l'\leq r(\frac{n+1}{n})$ (since $R=rq$) or $(X, \cE)\simeq(Q_n, \OO_{Q_n}(1)^{\oplus r})$ by Krull-Schmidt property of the category of vector bundles on $X$, see \cite[Theorem~3]{Atiyah56}.
\end{enumerate}

It remains to analyze the original Schur morphism in the quadric
case, i.e., to show that if $X\simeq Q_n$ then the map \eqref{Eq:geninjectivecharacterization} generically factors through $\bS_\lambda T_{Q_n}$. We first observe that the unipotent cokernel $\cU$ in \eqref{eq:stabilized-extension} is necessarily the trivial bundle $\OO_{Q_n}^{\oplus g}$ when $X\simeq Q_n$. Indeed since $n\geq2$, applying the cohomology long exact sequence to the short exact sequence 
\[0\rightarrow \OO_{\PP^{n+1}}(-2)\rightarrow \OO_{\PP^{n+1}}\rightarrow \OO_{Q_n}\rightarrow 0\] 
yields $H^1(Q_n,\OO_{Q_n})=0$. Clearly, $\cU_1$ in the filtration \eqref{Eq:unipotentfiltration} is trivial of rank $g_1$ and inductively suppose that $\cU_{i-1}$  is trivial of rank $g_1+\dots+g_{i-1}$. The extension
\begin{equation}\label{Eq:unipotentfiltration2}
 0\rightarrow\cU_{i-1}\rightarrow\cU_i
 \rightarrow\OO_{Q_n}^{\oplus g_i}\rightarrow0
\end{equation}
is classified by $\operatorname{Ext}^1_{Q_n}
 \bigl(\OO_{Q_n}^{\oplus g_i},\cU_{i-1}\bigr)
 \simeq H^1(Q_n,\cU_{i-1})^{\oplus g_i}=0$, where the final vanishing follows from triviality of $\cU_{i-1}$ and $H^1(Q_n,\OO_{Q_n})=0$. Thus, the short exact sequence \eqref{Eq:unipotentfiltration2} splits implying $\cU_i\simeq \cU_{i-1}\oplus \OO_{Q_n}^{\oplus g_i}\simeq\OO_{Q_n}^{\oplus g_1+\dots+g_i}$. Induction therefore yields $\cU\simeq\OO_{Q_n}^{\oplus g}$ and thus, on the quadric, \eqref{eq:stabilized-extension} has exactly the same form as in the trivial-cokernel case.

Put $Q:=\cU\simeq\OO_{Q_n}^{\oplus g}$. Since $\cE\simeq\OO_{Q_n}(1)^{\oplus r}$ and
$\bS_\lambda$ is homogeneous of degree $m=|\lambda|$, we have $\bS_\lambda\cE\simeq\bS_\lambda\CC^{\oplus r}\otimes\OO_{Q_n}(m)\simeq\OO_{Q_n}(m)^{\oplus d_\lambda(r)}$, where $d_\lambda(r):=\rk \bS_\lambda\CC^{\oplus r}$. Twisting the given morphism $\bS_\lambda\cE\rightarrow\bS_\lambda\cD$ by $\OO_{Q_n}(-m)$ yields an injective morphism
\begin{equation}\label{eq:twisted-quadric-Schur-map}
 \widetilde\theta:
 \OO_{Q_n}^{\oplus d_\lambda(r)}
 \rightarrow
 \bS_\lambda\cD\otimes\OO_{Q_n}(-m).
\end{equation}

Apply the Schur filtration on $\bS_\lambda \cD$ associated with the short exact sequence \eqref{eq:stabilized-extension}, see \cite[Proposition~2.3.1]{Weyman03} or more specifically, \cite[Proposition~2.7]{Fonarev15}.  The associated graded pieces of the filtration are 
\begin{equation}\label{Eq:filtrationgraded1}
\bigoplus_{\mu\subseteq\lambda}\bS_\mu T_{Q_n}\otimes\bS_{\lambda/\mu}Q,
\end{equation} 
where the sum runs over skew partitions $\lambda/\mu$, i.e., pair of partitions $(\mu, \lambda)$ satisfying the Young diagram inclusion $D_\mu\subseteq D_\lambda$. Here $\bS_{\lambda/\mu}(-)$ represents the skew-Schur functors associated to such skew partitions. By the Littlewood--Richardson decomposition of a skew Schur functor \cite[Proposition~2.6]{Fonarev15} or \cite[Proposition~2.3.6]{Weyman03} we have:
\[
 \bS_{\lambda/\mu}Q \simeq\bigoplus_{\nu}\bS_\nu Q\otimes_{\CC}\operatorname{LR}^{\lambda}_{\mu,\nu},
\]
where $\operatorname{LR}^{\lambda}_{\mu,\nu}$ is defined by \eqref{Eq:LR2}. Grouping the resulting summands according to the degree
$u:=|\mu|$ of the sub-partition $\mu\subseteq\lambda$ in \eqref{Eq:filtrationgraded1} yields a filtration $\mathscr{F}^\bullet$ on $\bS_\lambda\cD$ with associated graded pieces
\begin{equation}\label{Eq:filtrationgraded2}
 \operatorname{gr}_{\mathscr F}^{u}(\bS_\lambda\cD)
 \simeq
 \bigoplus_{\substack{\mu\vdash u,\ \nu\vdash m-u\\\mu\subseteq\lambda}}
 \left(
   \bS_\mu T_{Q_n}\otimes\bS_\nu Q
 \right)
 \otimes_{\CC}\operatorname{LR}^{\lambda}_{\mu,\nu}.
\end{equation}

Since $Q=\OO_{Q_n}^{\oplus g}$ is trivial, twisting the above constructions by $\OO_{Q_n}(-m)$ we obtain a filtration $\mathscr{G}^\bullet$ on $\bS_\lambda\cD\otimes \OO_{Q_n}(-m)$ with associated graded pieces
\begin{equation}\label{Eq:filtrationgraded3}
 \operatorname{gr}_{\mathscr G}^{u}
 \bigl(\bS_\lambda\cD\otimes\OO_{Q_n}(-m)\bigr)
 \simeq
 \bigoplus_{\substack{\mu\vdash u,\ \nu\vdash m-u\\ \\\mu\subseteq\lambda}}
 \left(
   \bS_\mu(T_{Q_n}(-1))\otimes\OO_{Q_n}(u-m)
 \right)
 \otimes_{\CC}
 \left(
   \bS_\nu\CC^g\otimes
   \operatorname{LR}^{\lambda}_{\mu,\nu}
 \right).
\end{equation}
For $0<u<m$, every bundle $\bS_\mu(T_{Q_n}(-1))\otimes\OO_{Q_n}(u-m)$ occurring in \eqref{Eq:filtrationgraded3} is a direct summand of $(T_{Q_n}(-1))^{\otimes u}\otimes\OO_{Q_n}(u-m)\simeq T_{Q_n}^{\otimes u}\otimes\OO_{Q_n}(-m)$ and hence by \cite[Lemma~2(1)]{DruelParis13}
\[
 H^0\!\left(Q_n,\bS_\mu(T_{Q_n}(-1))\otimes\OO_{Q_n}(u-m)
 \right)=0\qquad (0<u<m)
\]
For $Q_2\simeq\PP^1\times\PP^1$, the same vanishing follows directly: the summands of $T_{Q_2}^{\otimes u}\otimes\OO_{Q_2}(-m)$ have bidegrees $(2j-m,2(u-j)-m)$, $0\leq j\leq u$, and both entries cannot be nonnegative when $u<m$. The case $u=0$ is the elementary vanishing $H^0(Q_n,\OO_{Q_n}(-m))=0$.
Consequently,
\begin{equation}\label{Eq:H^0gradedvanishing}
 H^0\!\left(
   Q_n,
   \operatorname{gr}_{\mathscr G}^{u}
   \bigl(\bS_\lambda\cD\otimes\OO_{Q_n}(-m)\bigr) 
 \right)=0 \quad (u<m).
\end{equation}

Writing the filtration as $\mathscr{G}^0:=\bS_\lambda\cD\otimes \OO_{Q_n}(-m)\supseteq \mathscr{G}^1\supseteq\dots\supseteq\mathscr{G}^m\supseteq\mathscr{G}^{m+1}=0$, we have the following short exact sequence for every $u<m$:
\[0\rightarrow\mathscr G^{u+1}\rightarrow\mathscr G^u
 \rightarrow\operatorname{gr}_{\mathscr G}^{u}
 \bigl(\bS_\lambda\cD\otimes\OO_{Q_n}(-m)\bigr)
 \rightarrow0
\]
and \eqref{Eq:H^0gradedvanishing} implies that $H^0(Q_n,\mathscr G^{u+1})\xrightarrow{\sim} H^0(Q_n,\mathscr G^u)$ for $u<m$. Iterating, we see that the inclusion $\mathscr G^m\hookrightarrow\mathscr G^0 =\bS_\lambda\cD\otimes\OO_{Q_n}(-m)$ induces an isomorphism on global sections.  Since $\mathscr G^{m+1}=0$, the top filtration piece $\mathscr{G}^m$ equals $\operatorname{gr}_{\mathscr G}^{m}
 \bigl(\bS_\lambda\cD\otimes\OO_{Q_n}(-m)\bigr)\simeq\bS_\lambda (T_{Q_n}(-1))$ by \eqref{Eq:filtrationgraded3}. Thus, the natural inclusion $\mathscr{G}^m\hookrightarrow \mathscr{G}^0$ induces the isomorphism 
 \begin{equation} \label{Eq:H0Iso1}
 H^0(Q_n, \bS_\lambda (T_{Q_n}(-1)))=H^0(Q_n, \mathscr{G}^m)\xrightarrow{\sim}H^0(Q_n, \mathscr{G}^0)=H^0(Q_n, \bS_\lambda\cD\otimes \OO_{Q_n}(-m))
 \end{equation}
Since $\cE\simeq \OO_{Q_n}(1)^{\oplus r}$, the given morphism $\bS_\lambda\cE\rightarrow \bS_\lambda\cD$ is equivalent to a morphism $\bS_\lambda\OO_{Q_n}^{\oplus r}\rightarrow \bS_\lambda\cD\otimes \OO_{Q_n}(-m)$ by Lemma~\ref{lem:Schur-basic}(ii). Since $\bS_\lambda\OO_{Q_n}^{\oplus r}\simeq \OO_{Q_n}^{\oplus d_\lambda(r)}$, the morphism $\bS_\lambda\cE\rightarrow \bS_\lambda\cD$ is an element of $\Hom\!\left(\OO_{Q_n}^{\oplus d_\lambda(r)}, \bS_\lambda\cD\otimes \OO_{Q_n}(-m)\right)$. Furthermore, by \eqref{Eq:H0Iso1} we obtain
\begin{align*}\label{Eq:H0Iso2}
 \Hom\!\left(\OO_{Q_n}^{\oplus d_\lambda(r)}, \bS_\lambda\cD\otimes \OO_{Q_n}(-m)\right)&\simeq
 H^0(Q_n,\bS_\lambda\cD\otimes \OO_{Q_n}(-m))^{\oplus d_\lambda(r)}\\&\simeq H^0(Q_n, \bS_\lambda (T_{Q_n}(-1)))^{\oplus d_\lambda(r)}\\&\simeq \Hom\!\left(\OO_{Q_n}^{\oplus d_\lambda(r)}, \bS_\lambda (T_{Q_n}(-1))\right)\\
 &\simeq\Hom\!\left(\bS_\lambda\cE, \bS_\lambda T_{Q_n}\right),
\end{align*}
where the last isomorphism is by twisting by $\OO_{Q_n}(m)$. Since the above isomorphism of $\Hom$ groups is induced by the natural inclusion $\bS_\lambda T_{Q_n}\hookrightarrow \bS_\lambda\cD$, we see that the morphism $\bS_\lambda\cE\rightarrow \bS_\lambda\cD$ uniquely factors as
\[
 \bS_\lambda\cE\rightarrow
 \bS_\lambda T_{Q_n} \rightarrow\bS_\lambda\cD.
\]
Since the second arrow is injective and the composite is generically
injective, it follows that $\bS_\lambda\cE\rightarrow
 \bS_\lambda T_{Q_n}$ is injective. By Lemma~\ref{lem:Schur-rank-comparison}, we obtain $r\leq n$, since $\rk T_{Q_n}=n$. Padding with $0$'s on the right if necessary, let $\lambda=(\lambda_1,\dots,\lambda_n)$ be the partition of $m$, with original length of $\lambda$ being $s:=\length(\lambda)\leq r$. Twisting by $\OO_{Q_n}(-m)$, we obtain the injective morphism of vector bundles $\bS_\lambda\OO_{Q_n}^{\oplus r}\rightarrow\bS_\lambda( T_{Q_n}(-1))$ on a smooth projective variety, the induced map $H^0(Q_n, \bS_\lambda\OO_{Q_n}^{\oplus r})\rightarrow H^0(Q_n,\bS_\lambda( T_{Q_n}(-1)))$ is injective. By Lemmas~\ref{lem:SO-invariants} and \ref{lem:quadric-Schur-sections}, we have $\dim_\CC H^0(Q_n,\bS_\lambda( T_{Q_n}(-1)))=1$ and $\lambda_1\equiv\lambda_2\equiv\dots\equiv\lambda_n\pmod 2$. Consequently, we also have $d_\lambda(r)=\dim_\CC H^0(Q_n, \bS_\lambda\OO_{Q_n}^{\oplus r})=1$, whereby the proof of Lemma~\ref{lem:Schur-rank-comparison} yields that $\lambda_1=\dots=\lambda_s=k$ for some positive integer $k$. In particular, $s=r$. If $s<n$, then $\lambda_{s+1}=\dots=\lambda_n=0$, whereby $k$ must be even. Else $s=n$ (which further equals $r$, since $s\leq r\leq n$) and $\lambda_1=\dots=\lambda_n=k\in \ZZ^+$. This completes the proof.
\end{proof}

\subsection{Characterizing $\cE$ when $X\simeq\PP^n$}
In this section we characterize ample vector bundles $\cE$ of rank $r$ on $\PP^n$ admitting an injective morphism $\bS_\lambda\cE\hookrightarrow\bS_\lambda\cD$ of coherent sheaves, where $\cD$ fits into the short exact sequence \eqref{eq:unipotent-extension}. We first classify the unipotent vector bundle $\cU$ and the bundle $\cD$ in \eqref{eq:unipotent-extension} when $X\simeq \PP^n$.

\begin{lemma}\label{Lem:trivialunipotent}
    If $X\simeq \PP^n$, then $\cU\simeq \OO_{\PP^n}^{\oplus g}$ in \eqref{eq:unipotent-extension}. Thus, $\cD\simeq T_{\PP^n}\oplus \OO_{\PP^n}^{\oplus g}$.
\end{lemma}

\begin{proof}
Let the unipotent filtration of $\cU$ be $0=\cU_0\subset\cU_1\subset\cdots\subset\cU_e=\cU$, where $\cU_i/\cU_{i-1}\simeq\OO_{\PP^n}^{\oplus g_i}$, see \eqref{Eq:unipotentfiltration}. Since $H^1(\PP^n,\OO_{\PP^n})=0$, for each $1\leq i\leq e$ the resulting short exact sequence $0\rightarrow \cU_{i-1}\rightarrow \cU_i\rightarrow \OO_{\PP^n}^{\oplus g_i}\rightarrow 0$ splits inductively since $\operatorname{Ext}^1(\OO_{\PP^n}^{g_i}, \cU_{i-1})\simeq H^1(\PP^n, \cU_{i-1})^{\oplus g_i}=0$ (by inductive hypothesis $\cU_{i-1}\simeq \OO_{\PP^{n}}^{g_1+\dots+g_{i-1}}$). Hence $\cU\simeq\OO_{\PP^n}^{\oplus g}$. Furthermore since $\operatorname{Ext}^1(\OO_{\PP^n}, T_{\PP^n})\simeq H^1(\PP^n, T_{\PP^n})=0$, it follows that the extension \eqref{eq:unipotent-extension} is split when $X\simeq \PP^n$, yielding $\cD\simeq T_{\PP^n}\oplus \OO_{\PP^n}^{\oplus g}$.
\end{proof}

We now note a theorem of Elencwajg--Hirschowitz--Schneider on uniform vector bundles on $\PP^n$, which will be useful.

\begin{theorem}[\cite{EHS80}]\label{thm:EHS}
Let $\cF$ be a uniform vector bundle of rank $r\leq n$ on $\PP^n_{\mathbb{C}}$.
Then either $\cF$ is a direct sum of line bundles, or $r=n$ and $\cF$ is (up to tensoring by a line bundle)
the tangent bundle $T_{\PP^n}$ or the cotangent bundle $\Omega^1_{\PP^n}$.
\end{theorem}

By Theorem~\ref{thm:stabilized-schur-kovacs}, $\det\cE\simeq\OO_{\PP^n}(l)$ with $r\leq l\leq r(\frac{n+1}{n})$. We now break the characterization of $\cE$ into the following two cases.

\subsubsection{Case I: If $r\leq n$} In this case the bound $r\leq l\leq r(\frac{n+1}{n})=r+r/n$ implies that $l=r$ if $r<n$ and $l=r$ or $r+1$ if $r=n$. Let $L\simeq \PP^1\subseteq\PP^n$ be a line and let $\cE|_L \cong \bigoplus_{i=1}^r \mathcal{O}_L(a_i)$ with $a_1\geq a_2\geq \cdots \geq a_r$. Since $\cE$ is ample, so is the restriction $\cE|_L$ to any line $L\subset\PP^n$ and thus, we have $a_r\geq 1$. Then $\OO_{\PP^1}(l)\simeq\det(\cE)|_L\simeq \det(\cE|_L)\simeq\OO_{\PP^1}(a_1+\dots+a_r)$, whereby $a_1+\dots+a_r=l$. Consequently, 
\begin{itemize}
    \item If $r<n$, then $l=r$ and $a_1+\dots+a_r=l=r$ with $a_i\geq 1$ implies $a_i=1$ for all $1\leq i\leq r$. Thus, $\cE|_L\simeq\OO_{\PP^1}(1)^{\oplus r}$ for all lines $L\subseteq \PP^n$. By Theorem~\ref{thm:EHS}, we conclude $\cE\simeq\OO_{\PP^n}(1)^{\oplus r}$.
    \item If $r=n$ and $l=r$, then $a_i=1$ for all $1\leq i\leq r=n$ and so $\cE|_L\simeq\OO_{\PP^1}(1)^{\oplus n}$ for all lines $L\subseteq \PP^n$. By Theorem~\ref{thm:EHS}, we conclude $\cE\simeq\OO_{\PP^n}(1)^{\oplus n}$.
    \item If $r=n$ and $l=r+1$, then $a_1=2$ and $a_i=1$ for all $2\leq i\leq r=n$ and so $\cE|_L\simeq\OO_{\PP^1}(2)\oplus\OO_{\PP^1}(1)^{\oplus (n-1)}$ for all lines $L\subseteq \PP^n$. By Theorem~\ref{thm:EHS}, we conclude $\cE\simeq T_{\PP^n}$ or $\OO_{\PP^n}(2)\oplus\OO_{\PP^n}(1)^{\oplus (n-1)}$. 
    
    We now prove that $\cE\simeq \cE_0:= \OO_{\PP^n}(2)\oplus\OO_{\PP^n}(1)^{\oplus(n-1)}$ is possible only when the partition $\lambda\vdash m$ is $\lambda=(k^n)$ for some $k\geq1$. Adding zeroes to $\lambda$ on the right if necessary, let $\lambda=(\lambda_1,\ldots,\lambda_n)$ with $m=|\lambda|$. By Corollary~\ref{cor:Schur-sums-of-lines}, the highest degree summand in the line-bundle decomposition of $\bS_\lambda\cE_0$ is 
\[
 \cL_\lambda
 :=\OO_{\PP^n}(2)^{\otimes\lambda_1}\otimes\OO_{\PP^n}(1)^{\otimes(\lambda_2+\cdots+\lambda_n)}\simeq\OO_{\PP^n}(m+\lambda_1),
\]
since $\lambda_1\geq\lambda_2\geq\dots$. Since the given morphism $\bS_\lambda\cE_0\hookrightarrow\bS_\lambda\cD$ is injective, its restriction to $\cL_\lambda$ is nonzero. By Lemma~\ref{Lem:trivialunipotent}, $\cD\simeq T_{\PP^n}\oplus \OO_{\PP^n}^g$ and by Proposition~\ref{prop:LR-direct-sum} we therefore obtain
\begin{equation}\label{Eq:filtrationgraded4}
 \operatorname{gr}_{\mathscr F}^{u}(\bS_\lambda\cD)\simeq
 \bigoplus_{\substack{\mu\vdash u,\ \nu\vdash m-u}} \left(\bS_\mu T_{\PP^n}\otimes\bS_\nu\OO_{\PP^n}^{\oplus g}\right)\otimes_{\CC}
 \operatorname{LR}^{\lambda}_{\mu,\nu}.
\end{equation}
Since $\bS_\nu\OO_{\PP^n}^{\oplus g}$ is a trivial vector bundle of rank $d_\nu(g)$, it follows that $\operatorname{gr}_{\mathscr F}^{u}(\bS_\lambda\cD)$ is a direct sum of copies of $\bS_\mu T_{\PP^n}$, where $\mu\vdash u$ for all $u\leq m$. Thus, the non-zero map $\cL_\lambda\to\bS_\lambda\cD$ induces a non-zero map $\cL_\lambda\to \operatorname{gr}_{\mathscr F}^{u}(\bS_\lambda\cD)$ for some $u\leq m$. Composing this with a suitable projection to a direct summand of the chosen piece $\operatorname{gr}_{\mathscr F}^{u}(\bS_\lambda\cD)$, therefore gives a nonzero morphism $\cL_\lambda\simeq\OO_{\PP^n}(m+\lambda_1)\rightarrow\bS_\mu T_{\PP^n}$ for some subpartition $\mu\subseteq \lambda$ with $u:=|\mu|\leq m$. The case $u=0$ is
impossible because $m+\lambda_1>0$. Since $\bS_\mu T_{\PP^n}$ is a direct summand of $T_{\PP^n}^{\otimes u}$ when $u>0$, we obtain
\[
 H^0\!\left(
   \PP^n,
   T_{\PP^n}^{\otimes u}(-m-\lambda_1)
 \right)\neq0.
\]
By \cite[Lemma~1]{DruelParis13}, $m+\lambda_1\leq\frac{u(n+1)}{n}\leq\frac{m(n+1)}{n}$, and hence $\lambda_1\leq m/n$.  On the other hand, $m=\lambda_1+\cdots+\lambda_n\leq n\lambda_1$, so $\lambda_1\geq m/n$. Thus, equality holds throughout, and $\lambda_1=\cdots=\lambda_n=k$, or equivalently $\lambda=(k^n)$ for some $k\geq1$.
\end{itemize}
Conversely, when $\lambda=(k^n)$, the split bundle $\cE_0$ does occur: $\bS_{(k^n)}\cE_0\simeq(\det\cE_0)^{\otimes k}\simeq(\det T_{\PP^n})^{\otimes k}\simeq\bS_{(k^n)}T_{\PP^n}$, which injects into $\bS_{(k^n)}\cD$. This completes the proof of Theorem~\ref{thm:generalizedKovacs}(i), (ii), and (iii).

\subsubsection{Case II: If $r> n$} 
 Restricting the injection $\bS_\lambda\cE\rightarrow \bS_\lambda\cD$ to a generic line $L\simeq\PP^1\subseteq \PP^n$, we obtain the injection $\bS_\lambda(\cE|_L)\rightarrow \bS_\lambda(\cD|_L)$ by Proposition~\ref{prop:Schur-functorial}(iii). Here the general line meets the locus where the original map has full rank, so its restricted kernel is torsion and hence zero. By Lemma~\ref{Lem:trivialunipotent} and $T_{\PP^n}|_L\simeq\OO_{\PP^1}(2)\oplus\OO_{\PP^1}(1)^{\oplus (n-1)}$, we obtain
\begin{equation}\label{eq:D-line-splitting}
 \cD|_L\simeq
 \OO_{\PP^1}(2)\oplus\OO_{\PP^1}(1)^{\oplus(n-1)}\oplus\OO_{\PP^1}^{\oplus g}.
\end{equation}
Let $\cE|_L\simeq\bigoplus_{i=1}^r\OO_{\PP^1}(a_i)$ with $a_1\geq\cdots\geq a_r\geq1$ be the splitting type of $\cE$.

We will repeatedly use the following elementary observation. If
$\cV=\bigoplus_j\OO_{\PP^1}(b_j)$ is a vector bundle on $\PP^1$ and $c\in\mathbb{Z}$, denote $\cV_{\geq c}:=\bigoplus_{b_j\geq c}\OO_{\PP^1}(b_j)$ for the degree $\geq c$ summands of $\cV$. Since $\Hom(\OO_{\PP^1}(a),\OO_{\PP^1}(b))=0$ for $a>b$, every (injective) morphism $\cV\rightarrow \cV'$ of bundles on $\PP^1$ induces the (injective) morphism $\cV_{\geq c}\rightarrow \cV'_{\geq c}$ and hence
\begin{equation}\label{eq:HN-rank-inequality}
 \rk\cV_{\geq c}\leq\rk\cV'_{\geq c}
 \qquad\text{for every }c\in\ZZ.
\end{equation}

Set $\lambda_i=0$ for $i>\ell(\lambda)$ and let $m=|\lambda|$. By Corollary~\ref{cor:Schur-sums-of-lines}, the highest degree in the line-bundle decomposition of $\bS_\lambda(\cE|_L)$ is $M_{\cE}=\sum_{i=1}^r\lambda_i a_i$, whereas by \eqref{eq:D-line-splitting} the corresponding highest degree summand in $\bS_\lambda(\cD|_L)$ is
\begin{equation}\label{Eq:Ddeg}
M_{\cD}=2\lambda_1+\sum_{i=2}^n\lambda_i
 =m+\lambda_1-\sum_{i>n}\lambda_i
 \leq m+\lambda_1.
 \end{equation}
The injection $\bS_\lambda(\cE|_L)\rightarrow \bS_\lambda(\cD|_L)$ implies $M_{\cE}\leq M_{\cD}$. Suppose, for
a contradiction, that $a_1\geq2$. Then
\[
 m+\lambda_1\geq M_\cD\geq M_{\cE}
 =m+\sum_{i=1}^r\lambda_i(a_i-1)
 \geq m+\lambda_1,
\]
where the first inequality is \eqref{Eq:Ddeg}. Consequently equality holds throughout above, yielding
\begin{equation}\label{eq:equality-consequences}
 a_1=2,\qquad
 \lambda_i(a_i-1)=0\quad(i\geq2),\qquad \lambda_i=0\quad(i>n).
\end{equation}
We now derive a contradiction from \eqref{eq:equality-consequences}.

Assume first that the Young diagram of $\lambda$ is not a row, or equivalently, $\lambda_2> 0$ and put
$\overline\lambda:=(\lambda_2,\lambda_3,\dots, \lambda_n)$. Then \eqref{eq:equality-consequences} gives $\cE|_L\simeq\OO_{\PP^1}(2)\oplus\OO_{\PP^1}(1)^{\oplus(r-1)}$. By Corollary~\ref{cor:Schur-sums-of-lines}, the maximal-degree part of $\bS_\lambda(\cE|_L)$ is
\[
 \bigl(\bS_\lambda(\cE|_L)\bigr)_{\geq m+\lambda_1}
 \simeq
 \OO_{\PP^1}(m+\lambda_1)\otimes_{\mathbb{C}}
 \bS_{\overline\lambda}\mathbb{C}^{r-1}.
\]
Likewise, the maximal-degree part of
$\bS_\lambda(\cD|_L)$ is
\[
 \bigl(\bS_\lambda(\cD|_L)\bigr)_{\geq m+\lambda_1}
 \simeq
 \OO_{\PP^1}(m+\lambda_1)\otimes_{\mathbb{C}}
 \bS_{\overline\lambda}\mathbb{C}^{n-1}
\]
Applying \eqref{eq:HN-rank-inequality} with $c=m+\lambda_1$ gives
\[
 \dim_{\mathbb{C}}\bS_{\overline\lambda}\mathbb{C}^{r-1}
 \leq
 \dim_{\mathbb{C}}\bS_{\overline\lambda}\mathbb{C}^{n-1}.
\]
But this is impossible if $\overline\lambda=(\lambda_2,\dots,\lambda_n)$ is a non-empty partition since
$r-1>n-1\geq\ell(\overline\lambda)$, and the dimension of a nonzero
Schur module is strictly increasing with the dimension of the
underlying vector space (see Proposition~\ref{prop:Schur-functorial}(i)). Thus, it remains to consider the case when $\lambda=(m)$ is a row partition, in which case $\bS_\lambda=\Sym^m$ by Example~\ref{Ex:1}(1). Thus, the restricted injection $\bS_\lambda(\cE|_L)\rightarrow \bS_\lambda(\cD|_L)$ is now an injection $\Sym^m(\cE|_L)\rightarrow\Sym^m(\cD|_L)$ of symmetric powers. By \eqref{eq:equality-consequences} (where we now have $\lambda_2=\dots=\lambda_n=0$) and the ordering of the $a_i$, there
is an integer $j\geq1$ such that
\[
 \cE|_L\simeq
 \OO_{\PP^1}(2)^{\oplus j}\oplus\OO_{\PP^1}(1)^{\oplus(r-j)}.
\]
Thus, the degree-$2m$ (which is the highest degree) part of $\operatorname{Sym}^m(\cE|_L)$ has rank
$\dim\operatorname{Sym}^m\mathbb{C}^j$, while the degree-$2m$ part of
$\operatorname{Sym}^m(\cD|_L)$ has rank one from \eqref{eq:D-line-splitting}. Applying \eqref{eq:HN-rank-inequality} with $c=2m$ therefore forces $j=1$, whereby $\cE|_L\simeq\OO_{\PP^1}(2)\oplus\OO_{\PP^1}(1)^{\oplus(r-1)}$.
Now the degree-$\geq2m-1$ parts of
$\operatorname{Sym}^m(\cE|_L)$ and $\operatorname{Sym}^m(\cD|_L)$ have ranks, respectively,
\[
 1+(r-1)=r
 \qquad\text{and}\qquad
 1+(n-1)=n.
\]
A final application of \eqref{eq:HN-rank-inequality}, with
$c=2m-1$, gives $r\leq n$, contradicting the standing assumption
$r>n$. Thus we have contradicted our assumption $a_1\geq 2$, whereby $a_1=1$, and $\cE|_L\simeq\OO_{\PP^1}(1)^{\oplus r}$. The splitting on a general line gives $\det\cE\simeq\OO_{\PP^n}(r)$. On every line, ampleness makes all splitting degrees at least one; their sum is $r$, so $\cE(-1)$ is trivial on every line $L\subseteq \PP^n$.  Hence $\cE(-1)\simeq\OO_{\PP^n}^{\oplus r}$ by \cite[Chapter~I, Theorem~3.2.1]{OSS88}. Thus $\cE\simeq\OO_{\PP^n}(1)^{\oplus r}$ and $l=r$. This completes the proof of Theorem~\ref{thm:generalizedKovacs}(iv), and thus of Theorem~\ref{thm:generalizedKovacs}.

\section{Characterizing products of projective spaces and hyperquadrics among smooth projective varieties}\label{sec:CharacterizeProduct1}

The main theorem of this section is the following.

\begin{theorem}\label{thm:productcharacterization}
Let $X$ be a smooth projective complex variety of dimension $n>0$ and let $r>0$ be an integer. For $1\le i\le r$, let $\cE_i\subseteq T_X$ be an almost nef regular foliation of rank $n_i\geq1$, $\cV_i$ be a nef vector bundle of rank $r_i\ge1$, and let $\lambda_i$ be a nonzero partition of length $s_i:=\ell(\lambda_i)\le \min\{r_i,n_i\}$. Assume that
\begin{enumerate}
\item the natural morphism $\bigoplus_{i=1}^r\cE_i\rightarrow T_X$ is a subbundle inclusion,
\item for every $1\leq i\leq r$, there is an injective morphism, $\theta_i:\bS_{\lambda_i}\cV_i\hookrightarrow\bS_{\lambda_i}\cE_i$ of coherent $\OO_X$-modules and,
\item the vector bundle $\cA:=\bigotimes_{i=1}^r\cV_i$ is ample on $X$. 
\end{enumerate}
Then, possibly after reordering, there is an isomorphism
\begin{equation}\label{eq:main-product2}
X\simeq\prod_{i=1}^r Z_i,\qquad
\dim X=n=\sum_{i=1}^{r} n_i,\qquad
\cE_i\simeq\pr_i^*T_{Z_i}\subseteq T_X,
\end{equation}
where each $Z_i$ is $\PP^{n_i}$ or a smooth hyperquadric $Q_{n_i}\subset\PP^{n_i+1}$ with $n_i\geq2$ and $\pr_i:\prod_{i=1}^r Z_i\rightarrow Z_i$ is the natural projection. Furthermore $r_i\leq n_i$ and for all $1\leq i\leq r$, there are ample bundles $\cT_i$ on $Z_i$ such that $\cV_i\simeq\pr_i^*\cT_i$ and $\cT_i$ admits an injective morphism $\bS_{\lambda_i}\cT_i\rightarrow\bS_{\lambda_i}T_{Z_i}$ of coherent sheaves. These bundles $\cT_i$ on $Z_i$ are classified as follows.
\begin{enumerate}[(a)]
\renewcommand{\labelenumi}{\textup{(\alph{enumi})}}
\item If $Z_i\simeq\PP^{n_i}$ and $r_i<n_i$, then $\cT_i\simeq\OO_{\PP^{n_i}}(1)^{\oplus r_i}$,
\item If $Z_i\simeq\PP^{n_i}$, $r_i=n_i$ and $\det\cT_i\simeq\OO_{\PP^{n_i}}(n_i)$, then $\cT_i\simeq\OO_{\PP^{n_i}}(1)^{\oplus n_i}$,
\item If $Z_i\simeq\PP^{n_i}$, $r_i=n_i$ and $\det\cT_i\simeq\OO_{\PP^{n_i}}(n_i+1)$, then 
\[\cT_i\simeq T_{\PP^{n_i}} \text{ or } \cT_i\simeq\OO_{\PP^{n_i}}(2)\oplus\OO_{\PP^{n_i}}(1)^{\oplus(n_i-1)}.\]
The second possibility can occur if and only if $\lambda_i=(k^{n_i})$ for some integer $k\ge1$.
\item If $Z_i=Q_{n_i}$ with $n_i\geq 2$, then $\cT_i\simeq\OO_{Q_{n_i}}(1)^{\oplus r_i}$
 and the associated partition is $\lambda_i=(k_i^{r_i})$, where either $k_i$ is even or $r_i=n_i$.
\end{enumerate}
In particular, every $\cV_i$ is globally generated.
\end{theorem}

\subsection{A few useful lemmas}

The goal of this section is to record a few lemmas which will be useful at various stages of the proof of Theorem~\ref{thm:productcharacterization}. We first note the following structure theorem of almost nef regular foliations due to \cite{Iwai22}, which will allow us to prove Theorem~\ref{thm:productcharacterization} by inducting on the number $r$ of almost nef regular foliations $\cE_i\subseteq T_X$.

\begin{theorem}[{\cite[Theorem~4.1]{Iwai22}}]\label{thm:Iwai}
Let $X$ be a smooth $n$-dimensional projective variety and $\cF\subseteq T_X$ be an almost-nef regular foliation on $X$ with $c_1(\cF)\neq0$. Then there exist a smooth morphism $f:X\rightarrow Y$ between smooth projective varieties and a numerically flat regular foliation $\cG$ on $Y$, such that all the fibers of $f$ are rationally connected, $\cF=f^{-1}\cG$ and there is an exact sequence 
\begin{equation}\label{eq:Iwai}
0\rightarrow T_{X/Y}\rightarrow\cF\rightarrow f^*\cG\rightarrow0.
\end{equation}
\end{theorem}

\begin{remark}
    Note that the fibers of $f:X\rightarrow Y$ in Theorem~\ref{thm:Iwai} must be positive dimensional, as otherwise $f$ would be an isomorphism, whereby $\cF=f^{-1}\cG$ would be numerically flat, contradicting $c_1(\cF)\neq 0$.
\end{remark}

\begin{lemma}\label{lem:nefzero}
Let $X$ be a smooth projective variety of dimension $n>0$ with an ample line bundle $H$ and a nef line bundle $L$. If $L.H^{n-1}=0$, then $L\equiv0$ in $N^1(X)_\RR$.
\end{lemma}
\begin{proof}
For $n=1$, the statement follows from definition. For $n>1$, the assumption, $L.H^{n-1}=0$ along with generalized Hodge Index Theorem implies that $L^2.H^{n-2}<0$ if $L\not\equiv0$ in $N^1(X)_\RR$. But $L^2.H^{n-2}\geq 0$ since $H$ is ample and $L$ is nef, whereby we must have $L\equiv 0$ in $N^1(X)_\RR$.
\end{proof}

\begin{lemma}\label{lem:NF-RC}
Let $X$ be any smooth projective complex variety and $F$ be a smooth projective rationally connected complex variety.
\begin{enumerate}
    \item[(i)] An almost nef vector bundle $E$ on $X$ with $c_1(E)=0$ is numerically flat.
    \item[(ii)]  Every numerically flat bundle on $F$ is trivial. 
    \item[(iii)] $\Pic(F)$ is torsion free and $H^q(F,\OO_F)=0$ for all $q>0$.
\end{enumerate}
 
\end{lemma}

\begin{proof}
Statement (i) is \cite[Corollary~5.4]{HP19}. To prove (ii), let $\gamma:\PP^1\rightarrow F$ be a map to the smooth projective rationally connected complex variety $F$. Thus, if $\gamma^\ast E\simeq \bigoplus_{i=1}^{r}\OO_{\PP^1}(a_i)$ and $E$ is numerically flat, then $a_i=0$ for all $i=1,\dots, r$ since $\gamma^\ast E$ and $\gamma^\ast E^\vee$ are nef. This implies $\gamma^\ast E$ is trivial for all $\gamma:\PP^1\rightarrow F$ and thus, $E$ is trivial on $F$ by \cite[Proposition~2.1]{BiswasPedro09}. 

We now prove (iii). If $L$ is any torsion line bundle on $F$, then $c_1(L)=0$ and $L$ is nef, whereby (i) implies that $L$ is numerically flat, and (ii) implies that $L$ is trivial. This proves that $\Pic(F)$ is torsion free. The vanishing $H^q(F, \OO_F)=0$ for all $q>0$ is \cite[Corollary~4.18(a)]{Deb01}.
\end{proof}

\begin{lemma}\label{lem:positive-block}
\begin{enumerate}
 \item[(i)]  Let $\cV$ be a rank $r$ nef vector bundle. For a non-zero partition $\lambda$ of length $\length(\lambda)\leq r$ with $|\lambda|=m$, let $\iota:\bS_\lambda\cV\hookrightarrow\cP$ be an injective morphism of coherent sheaves into a numerically flat vector bundle $\cP$. Then $\cV$ is numerically flat.
\item[(ii)] In particular let $\cP=\bS_\lambda\cE$ for a rank $s$ almost nef vector bundle $\cE$ on $X$ and a non-zero partition $\lambda$ of length $\length(\lambda)\leq \min\{r,s\}$ and $\theta:\bS_\lambda\cV\hookrightarrow\bS_\lambda\cE$ be an injective morphism of coherent sheaves. If $c_1(\cV)\neq0$, then $c_1(\cE)\neq0$.
\end{enumerate}

\end{lemma}

\begin{proof}
 For (i), the case $\dim X=0$ is immediate. Set $n=\dim X>0$ and choose an ample line bundle $H$ on $X$. Let $U\subseteq X$ be the open locus where $\iota$ has rank equal to $\rk\bS_\lambda\cV$. Since $\iota$ is injective map of coherent sheaves, $U\subseteq X$ is nonempty. For $n\geq 2$, choose an integer $k\gg 0$ and sufficiently general divisors $H_1,\ldots,H_{n-1}\in|kH|$ such that by Bertini's theorem, $C:=H_1\cap\cdots\cap H_{n-1}$ is a smooth integral curve intersecting $U$. For $n=1$, one can take $C=X$ and $k=1$. In either case, we have
    \begin{equation}\label{Eq:positivedegree}
  \deg(\cV|_C)=c_1(\cV)\cdot C=k^{n-1}c_1(\cV)\cdot H^{n-1}\geq0.
\end{equation}
Since $C\cap U\neq \emptyset$ by construction, $\iota|_C: \bS_\lambda(\cV|_C)\hookrightarrow\cP|_C$ is generically injective on $C$ and thus, $\ker\iota|_C$ is a torsion subsheaf of the vector bundle $\bS_\lambda(\cV|_C)$, implying $\ker\iota|_C=0$. Thus, $\iota_C:\bS_\lambda(\cV|_C)\rightarrow\cP|_C$ is injective as a morphism of coherent sheaves. Let $\cI:=\im(\iota_C)\subseteq\cP|_C$, whereby $\cI\simeq\mathcal \bS_\lambda(\cV|_C)$. Let $\cT:=((\cP|_C)/\mathcal I)_{\mathrm{tors}}$ and $\cQ:=((\cP|_C)/\cI)/\cT$ and consider the saturation of $\cI$ in $\cP|_C$ defined by $\cS:=\ker(\cP|_C\rightarrow\mathcal Q)$. By definition, $\cQ$ is torsion-free while $\cS$, being a subsheaf of $\cP|_C$, is also torsion-free. Both are coherent and hence locally free on the smooth curve $C$. We have the following short exact sequences.
The definitions yield exact sequences
\begin{equation}\label{Eq:saturation-sequences}
\begin{gathered}
  0\rightarrow\bS_\lambda(\cV|_C)\rightarrow\cS\rightarrow\cT\rightarrow 0,\\
  0\rightarrow\cS\rightarrow\cP|_C\rightarrow\cQ\rightarrow 0.
\end{gathered}
\end{equation}
In the first sequence, $\cT$ has finite support on $C$ since it is a torsion sheaf. Consequently, $\bS_\lambda(\cV|_C)$ and $\cS$ have the same rank equal to $d_\lambda(r)$. Furthermore, we have
\[\det \cS=\det\bS_\lambda(\cV|_C)\otimes\OO_C(D),\]
where $D=\sum_{p\in C}\operatorname{length}_{\OO_C,p}(\cT_p)p$ is the effective divisor on $C$ corresponding to the torsion sheaf $\cT$. Thus, we obtain
\begin{equation}\label{Eq:degreeinequality}
\deg\cS=\deg \bS_\lambda(\cV|_C)+\deg D\geq \deg \bS_\lambda(\cV|_C).
\end{equation}
By the second short exact sequence of \eqref{Eq:saturation-sequences}, $\cQ$ is a quotient of the nef bundle $\cP|_C$ and therefore, nef itself by \cite[Proposition~1.15(i)]{DPS94} whereby $\deg\cQ\geq 0$. Since $\cV$ is nef, so is the restriction $\cV|_C$ and hence $\bS_\lambda(\cV|_C)$ is nef by Lemma~\ref{lem:schur-positivity}(ii). Thus, taking degrees in the second short exact sequence of \eqref{Eq:saturation-sequences}, we have
\[0\leq\deg\bS_\lambda(\cV|_C)\leq \deg\cS= \deg\cP|_C-\deg\cQ=-\deg\cQ\leq0,\]
since $\cP|_C$ is numerically flat. Thus, $\deg\bS_\lambda(\cV|_C)=0$. By Lemma~\ref{lem:Schur-basic}(i), $\deg\bS_\lambda(\cV|_C)=q_\lambda(r)\deg(\cV|_C)$ with $q_\lambda(r)>0$, whereby we obtain $\deg(\cV|_C)=0$ or equivalently, $c_1(\cV)H^{n-1}=0$. Since $\det\cV$ is nef, Lemma~\ref{lem:nefzero} yields $c_1(\cV)=0$ and then Lemma~\ref{lem:NF-RC}(i) implies that $\cV$ is numerically flat.  

For (ii), assume for contradiction that $c_1(\cE)=0$. By Lemma~\ref{lem:NF-RC}(i), $\cE$ is numerically flat and thus, $\bS_\lambda(\cE)$ is numerically flat by Lemma~\ref{lem:schur-positivity}(iii). Thus letting $\cP=\bS_\lambda(\cE)$, Lemma~\ref{lem:positive-block}(i) yields $c_1(\cV)=0$, which is the desired contradiction. 
\end{proof}

\begin{lemma}\label{lem:fiber-trivial}
Let $f:X\rightarrow Y$ be a smooth projective surjective morphism between smooth projective complex varieties, with rationally connected fibers. Let $\cV$ be a nef vector bundle of rank $r$ on $X$. Let $\cE$ be a vector bundle on $Y$ such that for a nonempty partition $\lambda\vdash m$ of length $\length(\lambda)\leq r$, there is an injection $\iota:\bS_\lambda\cV\hookrightarrow f^*\cE$ of coherent sheaves. Then $\cV|_F$ is trivial for every fiber $F$ of $f:X\rightarrow Y$. Consequently, $\cW:=f_*\cV$ is locally free of rank $r$, and the adjunction evaluation is an isomorphism $f^*\cW\xrightarrow{\sim}\cV$. The descended bundle $\cW$ is unique up to isomorphism.
\end{lemma}

\begin{proof}
By generic flatness \cite[Proposition~29.28.1]{stacks-project}, there exists a dense open subset $Y^\circ\subseteq Y$ over which the cokernel $\operatorname{coker}\iota$ of the injection $\iota:\bS_\lambda\cV\hookrightarrow f^*\cE$ is flat. Therefore for $y\in Y^\circ$, the restriction to the fiber $\iota_y:\bS_\lambda(\cV|_{F_y})\hookrightarrow (f^*\cE)|_{F_y}$ is injective. But $(f^*\cE)|_{F_y}$ is trivial bundle, whereby Lemma~\ref{lem:positive-block}(i) makes
$\cV|_{F_y}$ numerically flat and Lemma~\ref{lem:NF-RC}(ii) makes it
trivial. Thus, $\cL:=\det\cV$ satisfies $\cL_{F_y}\simeq \OO_{F_y}$ for $y\in Y^\circ$. By separatedness of the relative Picard scheme for projective flat morphism with geometrically integral fibers (see \cite[Théorème~3.1]{Grothendieck62}) we obtain $\cL_{F_y}\simeq \OO_{F_y}$ for all $y\in Y$. But $\cV$ is nef, whereby triviality of $\cL_{F_y}$ makes $\cV|_{F_y}$ numerically flat for all $y\in Y$ by Lemma~\ref{lem:NF-RC}(i), and furthermore $\cV|_{F_y}\simeq\OO_{F_y}^{\oplus r}$ by Lemma~\ref{lem:NF-RC}(ii). By Grauert's theorem \cite[Chapter III, Corollary~12.9]{Hartshorne77}, it follows that  $\cW:=f_*\cV$ is locally free of rank $r$, and the adjunction evaluation is an isomorphism $f^*\cW\xrightarrow{\sim}\cV$. The uniqueness of $\cW$ is a consequence of $f_\ast\OO_X\simeq\OO_Y$ and projection formula.
\end{proof}

\begin{setup}\label{Setup1}
 Now consider the situation of Theorem~\ref{thm:productcharacterization}. Since $\cA:=\cV_1\otimes \cV_2\otimes\dots \otimes \cV_r$ is ample on $X$, there exists $1\leq i\leq r$ such that $\cV_i$ is numerically non-trivial: this is because $\cV_i$ are nef and $c_1(\cA)=R\sum_{i=1}^{r}c_1(\cV_i)/r_i$, where $R:=\prod_{i=1}^{r}r_i$. Then applying Lemma~\ref{lem:positive-block} to the almost nef vector bundle $\cE_i$ in the hypothesis of Theorem~\ref{thm:productcharacterization}, admitting an injective morphism $\bS_{\lambda_i}\cV_i\rightarrow \bS_{\lambda_i}\cE_i$, we obtain $c_1(\cE_i)\neq 0$. By Theorem~\ref{thm:Iwai}, there exists a smooth projective variety $Y_i$ admitting a numerically flat regular foliation $\cG_i$ and a smooth morphism $f_i:X\rightarrow Y_i$ with rationally connected fibers such that $\cE_i=f_i^{-1}\cG_i$. Furthermore, there is a short exact sequence
 \begin{equation}\label{Equation:ses1}
 0\rightarrow T_{X/Y_i}\rightarrow \cE_i\rightarrow f_i^\ast\cG_i\rightarrow 0.\end{equation}  
 Furthermore, the fiber dimension of $f_i$ is strictly positive, as else, $f_i:X\rightarrow Y_i$ would be an isomorphism implying $\cE_i=f_i^{-1}\cG_i$ is numerically flat, contradicting $c_1(\cE_i)\neq 0$.
\end{setup} 

 \begin{lemma}\label{lem:fib-restriction}
     Let $F$ be any fiber of the morphism $f_i:X\rightarrow Y_i$ in Setup~\ref{Setup1} and $df_i:T_X\rightarrow f_i^\ast T_{Y_i}$ be the induced differential. Then
     \begin{enumerate}
         \item for all $j\neq i$, the restriction of $df_i$ to $\cE_j$ induces a subbundle inclusion $\cE_j\hookrightarrow f_i^\ast T_{Y_i}$
         \item $\cV_j|_F\cong \mathcal{O}_F^{\oplus\rk\cV_j}$ for all $j\neq i$,
         \item $\cV_i|_F$ is an ample vector bundle.
     \end{enumerate}
 \end{lemma}

\begin{proof}
    Since $f_i$ is smooth, we have the relative tangent sequence
\[
0\rightarrow T_{X/Y_i}\rightarrow T_X\xrightarrow{df_i} f_i^\ast T_{Y_i}\rightarrow 0.
\]
By \eqref{Equation:ses1}, $T_{X/Y_i}$ is a subbundle of $\cE_i$. Together with the given subbundle inclusion $\bigoplus_{k=1}^r\cE_k\hookrightarrow T_X$, this implies that $T_{X/Y_i}\oplus\cE_j$ is a subbundle of $T_X$ for $j\neq i$. Quotienting by $T_{X/Y_i}$ therefore gives the subbundle inclusion $\cE_j\hookrightarrow T_X/T_{X/Y_i}\simeq f_i^\ast T_{Y_i}$
induced by $df_i$, proving (1).

    Now for all $j\neq i$, the composed map $\bS_{\lambda_j}\cV_j\rightarrow \bS_{\lambda_j}\cE_j\hookrightarrow \bS_{\lambda_j}(f_i^\ast T_{Y_i})\simeq f_i^\ast (\bS_{\lambda_j}T_{Y_i})$ (using Proposition~\ref{prop:Schur-functorial}(iii)) is injective. Thus (2) follows from Lemma~\ref{lem:fiber-trivial}.

  By assumption $\cA:=\cV_1\otimes \cV_2\otimes\dots\otimes \cV_r$ is an ample vector bundle, whereby restricting to any fiber $F$, we have $(\cV_i|_F)^{\oplus N}\cong (\cV_1\otimes \cV_2\otimes \dots\otimes \cV_r)|_F$ is ample for $N=\prod_{j\neq i}\rk\cV_j|_F$ since by (2), $\cV_j|_F$ is a trivial vector bundle on $F$ for all $j\neq i$. Consequently, $\cV_i|_F$ is ample, proving (3). 
\end{proof}

The following lemma characterizes the fibers of $f_i:X\rightarrow Y_i$.

\begin{lemma}\label{lem:fib-characterization}
Let $f_i:X\rightarrow Y_i$ be the smooth morphism as in Setup~\ref{Setup1} with fiber dimension $d\geq1$. Then either:
 \begin{enumerate}
       \item[(i)] all fibers of $f_i:X\rightarrow Y_i$ are isomorphic to $\PP^d$, or
       \item[(ii)] all fibers of $f_i:X\rightarrow Y_i$ are isomorphic to $Q_d$, where $Q_d\subseteq \PP^{d+1}$ is a smooth hyperquadric.
  \end{enumerate}
\end{lemma}
\begin{proof}
By Setup~\ref{Setup1}, $Y_i$ is a smooth projective variety admitting a numerically flat regular foliation $\cG$ and the smooth morphism $f_i:X\rightarrow Y_i$ has rationally connected fibers of dimension $d>0$ such that $\cE_i\cong f_i^{-1}\cG$. Restricting the associated short exact sequence \eqref{Equation:ses1} to any fiber $F$ of $f_i:X\rightarrow Y_i$ yields
   \[0\rightarrow T_F\rightarrow \cE_i|_F\rightarrow \OO_F^{\oplus\rk\cG}\rightarrow 0,\]
   since $f_i^\ast\cG|_F\simeq\OO_F^{\oplus\rk\cG}$. By generic flatness \cite[Proposition~29.28.1]{stacks-project}, there exists a dense open subset $Y_i^\circ\subseteq Y_i$ over which the cokernel $\operatorname{coker} \theta_i$ of the injective map $\theta_i:\bS_{\lambda_i}\cV_i\rightarrow\bS_{\lambda_i}\cE_i$ is flat. Thus, for all $y\in Y_i^\circ$, $\theta_i$ restricts to an injective map $\theta_i|_{F_y}:\bS_{\lambda_i}(\cV_i|_{F_y})\rightarrow\bS_{\lambda_i}(\cE_i|_{F_y})$ of coherent sheaves on the fiber $F_y$ by Proposition~\ref{prop:Schur-functorial}(iii). By Lemma~\ref{lem:fib-restriction}(3), the vector bundle $\cV_i|_{F_y}$ is ample on $F_y$. By Theorem~\ref{thm:generalizedKovacs}, we therefore have $F_y\simeq \PP^d$ or $(F_y, \cV_i|_{F_y})\simeq(Q_d, \OO_{Q_d}(1)^{\oplus r_i})$ for all $y\in Y^\circ$. Furthermore by non-deformability of $\PP^d$ for $d\geq 1$ and $Q_d$ for $d\geq 3$ (see \cite[Theorem~1']{HwangMok98} for general non-deformability of irreducible Hermitian symmetric spaces of compact type or alternatively, \cite{Siu89, Siu92_Errata} and \cite{Hwang95} for projective space and quadrics), we obtain the following:
   \begin{enumerate}
       \item[(i)] all fibers of $f_i:X\rightarrow Y_i$ are isomorphic to $\PP^d$, or
       \item[(ii)] all fibers of $f_i:X\rightarrow Y_i$ are isomorphic to $Q_d$ when $d\geq 3$.
  \end{enumerate}    
  Since $\PP^1\simeq Q_1$, the only remaining case is if $F_y\simeq Q_2\simeq \PP^1\times\PP^1$ for $y\in Y_i^\circ$. Note that by \cite[Remark~22]{Araujo09}, $Q_2$ can degenerate to any Hirzebruch surface $\mathbb{F}_a$ for $a$ even. However, we will show that this is not possible in our setup since by Theorem~\ref{thm:generalizedKovacs}, we have $\cV_i|_{F_y}\simeq\OO_{Q_2}(1)^{\oplus r_i}$ or equivalently $\OO_{\PP^1\times\PP^1}(1,1)^{\oplus r_i}$ when $F_y\simeq Q_2\simeq \PP^1\times\PP^1$ for all $y\in Y_i^\circ$. To this end, let $\cL_i:=\det\cV_i$ and let
\[
  \cM_i:=\cL_i^{\otimes2}\otimes\omega_{X/Y_i}^{\otimes r_i}.
\]
Since $\cL_i^{\otimes 2}|_{F_y}\simeq \OO_{\PP^1\times\PP^1}(2r_i, 2r_i)$ and $\omega_{X/Y_i}|_{F_y}\simeq \omega_{F_y}\simeq\OO_{\PP^1\times\PP^1}(-2, -2)$ for $y\in Y_i^\circ$, we obtain $\cM_i|_{F_y}\simeq\mathcal O_{F_y}$ for all $y\in Y_i^\circ$. By separatedness of the relative Picard scheme for projective flat morphism with geometrically integral fibers (see \cite[Théorème~3.1]{Grothendieck62}) we obtain $\cM_i|_{F_y}\simeq\mathcal O_{F_y}$, or equivalently, $\omega_{F_y}^{\otimes -r_i}\simeq \cL_i^{\otimes 2}|_{F_y}$ for all $y\in Y_i$. Since $\cV_i|_{F_y}$ is ample for all $y\in Y_i$ by Lemma~\ref{lem:fib-restriction}(3), it follows that $\omega_{F_y}^{-1}$ is ample for all $y\in Y_i$, and thus every fiber of $f_i:X\rightarrow Y_i$ is a smooth del Pezzo surface. The degree $K^2_{F_y}$ of the surfaces $F_y$ is constant in the family $f_i:X\rightarrow Y_i$ by constancy of Euler characteristic on connected varieties and the surface Riemann-Roch theorem \cite[Chapter V, Theorem~1.6]{Hartshorne77}, whereby every fiber of $f_i:X\rightarrow Y_i$ is a degree $8$ smooth del Pezzo surface. By \cite[\S 0, pg 102-03]{MM83} smooth complex del Pezzo surfaces of degree eight form two distinct deformation classes, represented by $\PP^1\times\PP^1$ and the Hirzebruch surface $\mathbb F_1$. Consequently, in the connected smooth flat family $f_i:X\rightarrow Y_i$ of degree $8$ del Pezzo surfaces, the presence of one fiber isomorphic to $\mathbb P^1\times\mathbb P^1$ implies that every fiber is isomorphic to $\mathbb P^1\times\mathbb P^1$. This proves the remaining alternative that all the fibers of $f_i:X\rightarrow Y_i$ are isomorphic to $Q_2$.

\end{proof}

The following proposition gives a condition for the descent of an almost nef subbundle of the pulled back tangent bundle.

\begin{proposition}\label{Prop:almost-nef-descent}
Let $f:X\rightarrow Y$ be a smooth surjective morphism with rationally connected fibers between smooth irreducible projective
complex varieties, and let $\cE\hookrightarrow f^*T_Y$
be a rank-$r$ subbundle such that $\cE$ is almost nef. Then there
is a unique rank-$r$ subbundle $\cF\hookrightarrow T_Y$
such that $\cE\cong f^*\cF$ as subbundles of $f^*T_Y$.
\end{proposition}

\begin{proof}
If $f$ has zero-dimensional fibers, it is an isomorphism and there is nothing to prove. Assume the fiber dimension is positive. Let $\operatorname{Gr}_Y(r,T_Y)$ be the relative Grassmannian of rank-$r$ subspaces of $T_Y$, and let $\mathcal{S}\hookrightarrow \pi^*T_Y$ be the universal rank-$r$ subbundle, where $\pi:\operatorname{Gr}_Y(r,T_Y)\rightarrow Y$ is the structure
morphism. The given subbundle $\cE\hookrightarrow f^*T_Y$ defines a classifying morphism $ \phi: X\rightarrow \operatorname{Gr}_Y(r,T_Y)$ over $Y$ such that $\phi^*\mathcal{S}\cong \cE$. Let $Z=\bigcup_{m\geq 1} Z_m\subsetneq X$ be a countable union of proper closed subvarieties such that $\cE|_C$ is nef for every irreducible curve $C\subset X$ not contained in
$Z$. For every $m$, let $W_m:=\{y\in Y\mid X_y\subseteq Z_m\}$.
Since $f$ is smooth, it is universally open and hence $Y\setminus W_m=f(X\setminus Z_m)$ is a nonempty open subset. Thus, $W_m\subseteq Y$ is a possibly empty proper closed subset for all $m\geq 1$. 

Let $F$ be a very general fiber of $f:X\rightarrow Y$ such that $F\not\subseteq Z_m$ for every $m$ or equivalently, each $Z_m\cap F$ is a proper closed subset of $F$.  Hence, by uncountability of the base field $\CC$, one may choose a very general complete intersection curve $C\subset F$ such that $C\not\subseteq \bigcup_{m\geq 1} (Z_m\cap F)$. By almost-nefness of $\cE$ on $X$, the bundle $\cE|_C$ is nef on $C$. In particular, $\deg\det(\cE|_C)\geq 0$.

Restricting the subbundle inclusion $\cE\hookrightarrow f^\ast T_Y$ to the fiber $F$, we have the subbundle inclusion $\cE|_F\hookrightarrow f^\ast T_Y|_F\cong T_{Y,y}\otimes\mathcal{O}_F$ into a trivial vector bundle, where $T_{Y,y}$ is the fiber of the tangent bundle $T_Y$ at $y=f(F)\in Y$. Dualizing, we see that $(\cE|_F)^\vee$ is a quotient of a trivial bundle, and is thus globally generated. Hence $\det(\cE|_F)^\vee$ is globally generated and nef, whereby $\deg\det (\cE|_C)^\vee\geq 0$, or equivalently, $\deg\det(\cE|_C)\leq 0$ for a very general complete intersection curve $C\subset F$. Combining with the inequality from the previous paragraph, we have $\deg\det(\cE|_C)= 0$ if furthermore $C\not\subseteq \bigcup_{m\geq 1} (Z_m\cap F)$. Equivalently, we have $c_1(\cE|_F).H^{d-1}=0$ where $H$ is an ample line bundle on $F$ and $\dim F=d$, which implies $c_1(\cE|_F)=0$ by Lemma~\ref{lem:nefzero} applied to the nef line bundle $\det(\cE|_F)^\vee$. Since $(\cE|_F)^\vee$ is nef, Lemma~\ref{lem:NF-RC}(i) makes it numerically flat. Thus, $\cE|_F$ is numerically flat, and hence trivial on $F$ by Lemma~\ref{lem:NF-RC}(ii), since $F$ is rationally connected. Consequently, $\det\cE|_F\cong \OO_F$. Let $\mathcal{O}_{\operatorname{Gr}_Y(r,T_Y)}(1)\cong \det(\mathcal{S}^\vee)$ be the Plücker line bundle on $\operatorname{Gr}_Y(r,T_Y)$, where $\mathcal{S}\hookrightarrow \pi^*T_Y$ is the universal rank-$r$ subbundle. Thus,
\[\phi^\ast\mathcal{O}_{\operatorname{Gr}_Y(r,T_Y)}(1)\cong \phi^\ast \det(\mathcal{S}^\vee)\cong \det (\phi^\ast \mathcal{S}^\vee)\cong \det(\cE^\vee) \cong \det(\cE)^\vee\]
Restricting to a very general fiber $F=X_y$ as above, we have $(\phi|_F)^\ast\mathcal{O}_{\operatorname{Gr}_Y(r,T_{Y,y})}(1)\cong  \det(\cE|_F)^\vee\cong \mathcal{O}_F$. Since $\mathcal{O}_{\operatorname{Gr}_Y(r,T_{Y,y})}(1)$ is ample on $\operatorname{Gr}_Y(r,T_{Y,y})$, it follows that $\phi|_F$ is constant.

Let $S$ be the dense set of closed points of $Y\setminus\bigcup_m W_m$. Let $R:=X\times_YX$, and denote the two projections by $p_1,p_2: R\rightrightarrows X$. Consider the $Y$-morphisms $\phi\circ p_1$, $\phi\circ p_2: R\rightarrow\operatorname{Gr}_Y(r,T_Y)$. For every $y\in S$, the constancy of $\phi|_{X_y}$ is equivalent to
\[
  (\phi\circ p_1)|_{X_y\times X_y}=(\phi\circ p_2)|_{X_y\times X_y}.
\]
Thus, $\phi\circ p_1$ and $\phi\circ p_2$ agree on the subset $R_S:=R\times_Y S=h^{-1}(S)\subseteq R$, where $h:=f\circ p_1=f\circ p_2: R\rightarrow Y$ is the natural smooth and surjective morphism. Therefore $h$ is open and $R_S=h^{-1}(S)$ is dense in $R$ since $S$ is dense in $Y$. Furthermore, since the relative Grassmannian $\operatorname{Gr}_Y(r,T_Y)\rightarrow Y$ is separated, the scheme-theoretic equalizer $D\subseteq R$ of $\phi\circ p_1, \phi\circ p_2: R\rightrightarrows \operatorname{Gr}_Y(r,T_Y)$ is a closed subscheme by \cite[Lemma~26.21.5]{stacks-project}. Since the underlying closed set of $D$ contains the dense subset $R_S\subset R$, it follows that $|D|=|R|$, i.e., $D$ and $R$ have the same underlying topological space. But since $R$ is reduced, it follows that $D=R$ scheme-theoretically, and hence
\begin{equation}\label{Equation:descent}
  \phi\circ p_1=\phi\circ p_2
  \qquad\text{on }X\times_YX.
\end{equation}
Since $f:X\rightarrow Y$ is smooth and surjective morphism of varieties, it is an fpqc cover. Thus, \eqref{Equation:descent} is the fpqc descent condition for the morphism $\phi:X\rightarrow \operatorname{Gr}_Y(r,T_Y)$. Since morphism of schemes satisfy fpqc descent (see \cite[Lemma~35.35.11 and Remark 35.35.12]{stacks-project}), we get an induced morphism $\overline{\phi}:Y\rightarrow \operatorname{Gr}_Y(r,T_Y)$ such that $\phi=\overline{\phi}\circ f$. Since $\phi$ is a $Y$-morphism, the induced morphism $\overline{\phi}$ is in fact a section to the structure morphism $\pi:\operatorname{Gr}_Y(r,T_Y)\rightarrow Y$. By the universal property of the relative Grassmannian scheme, define $\cF:=\overline{\phi}^\ast\mathcal{S}\hookrightarrow T_Y$ on $Y$, whereby $f^\ast\cF\cong f^\ast\overline{\phi}^\ast\mathcal{S}\cong\phi^\ast\mathcal{S}\cong \cE$, and the subbundle inclusion $\cF\hookrightarrow T_Y$ is pulled back to the given subbundle inclusion $\cE\cong f^\ast\cF\hookrightarrow f^\ast T_Y$. Uniqueness of the bundle $\cF$ follows from the uniqueness of fpqc descent.

\end{proof}

\subsection{Proof of Theorem~\ref{thm:productcharacterization}}\label{subsec:proofproduct} We prove Theorem~\ref{thm:productcharacterization} by induction on $r$.

\subsubsection{Base case $r=1$}\label{subsubsec:product-base-case} In this case, the hypotheses of Theorem~\ref{thm:productcharacterization} boil down to an injective morphism $\bS_{\lambda_1}\cV_1\rightarrow \bS_{\lambda_1}T_X$ (since $\bS_{\lambda_1}(-)$ preserves subbundle inclusions), where $\cV_1$ is an ample vector bundle on $X$. Thus, the hypotheses of Theorem~\ref{thm:productcharacterization} for $r=1$ imply those of Theorem~\ref{thm:generalizedKovacs} with the unipotent bundle $\cU$ being $\mathbf{0}_X$, whereby all the conclusions of the base case $r=1$ of Theorem~\ref{thm:productcharacterization} except $\cE_1\simeq T_X$ follow from Theorem~\ref{thm:generalizedKovacs}. To prove $\cE_1\simeq T_X$, we consider the following two cases:
\begin{enumerate}
    \item If $X\not\simeq Q_2$: Then $X\simeq \PP^n$ or $X\simeq Q_n$ with $n\geq 3$, whereby $X$ has Picard number one. Since $\cV_1$ is ample, we have $c_1(\cV_1)\not\equiv0$. Consequently, Lemma~\ref{lem:positive-block}(ii) and Theorem~\ref{thm:Iwai} yield a smooth morphism $f:X\rightarrow Y$ with positive-dimensional rationally connected fibers. If $\dim Y>0$, let $H$ be an ample line bundle on $Y$, whereby $f^\ast H$ is a nef line bundle on $X$ with $c_1(f^\ast H)\neq 0$. But since $\rho(X)=1$, we have $f^\ast H=cL$ for some $c>0$, where $L$ is the ample generator of $N^1(X)_\RR$. This contradicts the fact that $(f^\ast H).C=0$ for any curve $C$ in fiber of $f:X\rightarrow Y$. Thus, $Y$ is a point whereby $\cG=0$ in the associated short exact sequence \eqref{Equation:ses1}, implying $\cE_1=T_X$.
    \item If $X\simeq Q_2\simeq \PP^1\times\PP^1$: Assume $\rk\cE_1=1$, whereby $\rk \cV_1=1$ by Lemma~\ref{lem:Schur-rank-comparison}. From Theorem~\ref{thm:generalizedKovacs}, the partition $\lambda_1\vdash m$ is the trivial partition $\lambda_1=(m)$, whereby $\bS_{\lambda_1}\cF\simeq\cF^{\otimes m}$ for $\cF=\cV_1, \cE_1$. Using the isomorphism $Q_2\simeq \PP^1\times \PP^1$, let $\cV_1=\OO(a,b)$ and $\cE_1=\OO(c,d)$ for $a,b,c,d\in\ZZ$. Ampleness of $\cV_1$ implies $a,b>0$. Furthermore, the injective morphism $\cV_1^{\otimes m}\to\cE_1^{\otimes m}$ is a nonzero section of $\OO(m(c-a),m(d-b))$, so $c\geq a>0$ and $d\geq b>0$. But
\[
\begin{aligned}
 \Hom(\cE_1,T_{Q_2})=H^0\bigl(Q_2,\OO(2-c,-d)\bigr)\oplus H^0\bigl(Q_2,\OO(-c,2-d)\bigr)=0,
\end{aligned}
\]
since $T_{\PP^1\times\PP^1}=\OO(2,0)\oplus\OO(0,2)$ and $-d,-c<0$. This contradicts the given subbundle inclusion $\cE_1\hookrightarrow T_{Q_2}$. Hence we must have $\rk\cE_1=2$, and since the cokernel of $\cE_1\hookrightarrow T_{Q_2}$ is locally free, it must have rank zero. This proves the isomorphism $\cE_1\simeq T_{Q_2}$ and completes the proof of the base case of the inductive proof of Theorem~\ref{thm:productcharacterization}.
\end{enumerate}

\subsubsection{Inductive step}

The following proposition provides the inductive reduction step. Following \cite{Araujo09}, a \emph{quadric bundle} is defined to be a flat projective morphism between quasi-projective varieties all of whose fibers are isomorphic to irreducible and reduced hyperquadrics in projective space.

\begin{proposition}\label{Proposition:4}
Let $X$ be a smooth projective complex variety of dimension $n>0$ and let $r>0$ be an integer. For $1\le i\le r$, let $\cE_i\subseteq T_X$ be an almost nef regular foliation of rank $n_i\geq1$, $\cV_i$ be a nef vector bundle of rank $r_i\ge1$, and let $\lambda_i$ be a nonzero partition of length $s_i:=\ell(\lambda_i)\le \min\{r_i,n_i\}$. Assume that
\begin{enumerate}
\item the natural morphism $\bigoplus_{i=1}^r\cE_i\rightarrow T_X$ is a subbundle inclusion,
\item for every $1\leq i\leq r$, there is an injective morphism, $\theta_i:\bS_{\lambda_i}\cV_i\hookrightarrow\bS_{\lambda_i}\cE_i$ of coherent $\OO_X$-modules and,
\item the vector bundle $\cA:=\bigotimes_{i=1}^r\cV_i$ is ample on $X$. 
\end{enumerate}
Then, after reordering the indices, there exists a Brauer-Severi scheme or smooth quadric bundle $f:X\rightarrow Y$, with connected fibers of positive dimension, where $Y$ is a smooth projective variety and $\cE_1=f^{-1}\cG$ for a numerically flat regular foliation $\cG\subseteq T_Y$. Thus the exact sequence \eqref{Equation:ses1} holds for $f$. The following properties are satisfied:
\begin{enumerate}
    \item[(a)] there exist nef vector bundles $\cW_2,\dots,\cW_r$ on $Y$ such that $\cV_j\cong f^*\cW_j$ for all $j\ge 2$,
    \item[(b)] there exist almost nef regular foliation $\cF_j\subseteq T_Y$ such that $\cE_j\cong f^\ast \cF_j$ for all $j\geq 2$. Furthermore the injective morphisms $\theta_j:\bS_{\lambda_j}\cV_j\rightarrow \bS_{\lambda_j}\cE_j$ of coherent sheaves descend to injective morphisms $\bS_{\lambda_j}\cW_j\rightarrow \bS_{\lambda_j}\cF_j$ for all $j\geq 2$ and,
    \item[(c)] there is a subbundle inclusion $\cG\oplus \cF_2\oplus\cdots\oplus \cF_r\hookrightarrow T_Y$.
\end{enumerate}
\end{proposition}

\begin{proof}
   By assumption, one of the vector bundles $\cV_1,\dots, \cV_r$ is numerically non-trivial, and thus by Lemma~\ref{lem:positive-block}, one of the $\cE_1,\dots, \cE_r$ has non-trivial first Chern class. Without loss of generality, up to relabelling, we may assume $\cV_1$ is numerically non-trivial, whereby $c_1(\cE_1)\neq 0$. Thus, we are in Setup~\ref{Setup1} with $i=1$, i.e., we have a smooth projective variety $Y$ admitting a numerically flat regular foliation $\cG$ and a smooth morphism $f:X\rightarrow Y$ with rationally connected fibers such that $\cE_1\cong f^{-1}\cG$. By Lemma~\ref{lem:fib-characterization}, either all the fibers of $f:X\rightarrow Y$ are isomorphic to $\PP^d$ for $d\geq 1$, in which case $f:X\rightarrow Y$ is a Brauer-Severi scheme by \cite[Corollary~ 8.3]{BSI}; or all the fibers of $f:X\rightarrow Y$ are isomorphic to the smooth quadric $Q_d$ for $d\geq 1$, in which case $f:X\rightarrow Y$ is a smooth quadric bundle, as defined in \cite{Araujo09}.

   By Lemma~\ref{lem:fib-restriction}(2), $\cV_j$ restricted to any fiber of $f:X\rightarrow Y$ is the trivial vector bundle of rank $\rk\cV_j$ for all $j=2,\dots, r$. Thus, for each $j\geq 2$, $\cV_j\cong f^\ast \cW_j$ for a unique vector bundle $\cW_j$ on $Y$ by Lemma~\ref{lem:fiber-trivial}. By \cite[Proposition~1.10]{DPS94}, $\cW_j$ is nef for all $j\geq 2$ since $\cV_j$ is. This proves (a).

    By Lemma~\ref{lem:fib-restriction}(1), there is a subbundle inclusion $\cE_j\hookrightarrow f^\ast T_Y$ for all $j\geq 2$. By Proposition~\ref{Prop:almost-nef-descent}, there are subbundles $\cF_j\hookrightarrow T_Y$ such that $\cE_j\cong f^\ast \cF_j$ as subbundles of $f^\ast T_Y$ for all $j\geq 2$. By \cite[Lemma~2.6(2)]{IMZ25}, $\cF_j$ is almost nef. Since $\cF_j$ is already a subbundle of $T_Y$, to prove that is a regular foliation, it suffices to prove that it is closed under the Lie bracket. Let $u,v$ be local sections of $\cF_j$. Under $df|_{\cE_j}:\cE_j\xrightarrow{\sim} f^*\cF_j$, their pullbacks have unique lifts $\widetilde u,\widetilde v \in \cE_j$. These are projectable vector fields by definition. Since $\cE_j$ is a foliation, $[\widetilde u,\widetilde v]\in\cE_j$. The compatibility of brackets of projectable vector fields with $df$ yields $f^*[u,v]=df[\widetilde u,\widetilde v]\in f^\ast\cF_j$, whereby $[u,v]\in\cF_j$ since $f$ is faithfully flat. This proves that $\cF_j$ is a regular foliation. Furthermore, since $f_*\OO_X=\OO_Y$, adjunction and the projection formula identify $\Hom_X(f^*\cA,f^*\cB)$ with $\Hom_Y(\cA,\cB)$ for vector bundles $\cA,\cB$ on $Y$. Thus, the injection $\bS_{\lambda_j}\cV_j\rightarrow \bS_{\lambda_j}\cE_j$, which is equivalent to the injection $f^\ast(\bS_{\lambda_j}\cW_j)\rightarrow f^\ast(\bS_{\lambda_j}\cF_j)$ by Proposition~\ref{prop:Schur-functorial}(iii), descends to an injection $\bS_{\lambda_j}\cW_j\rightarrow \bS_{\lambda_j}\cF_j$ for all $j\geq 2$. The descended morphisms are injective by faithful flatness. This proves (b).

    Quotienting the given subbundle inclusion $\bigoplus_{i=1}^{r}\cE_i\hookrightarrow T_X$ by $T_{X/Y}\subseteq \cE_1$ (where the inclusion is from \eqref{Equation:ses1}), we obtain the subbundle inclusion $f^\ast\cG\oplus\bigoplus_{j=2}^{r}\cE_j\hookrightarrow f^\ast T_Y$, where $\cG$ is a numerically flat regular foliation on $T_Y$ such that $\cE_1/T_{X/Y}\cong f^\ast \cG$ by \eqref{Equation:ses1}. By (b), we therefore have $f^\ast(\cG\oplus\bigoplus_{j=2}^{r}\cF_j)\hookrightarrow f^\ast T_Y$, which implies $\cG\oplus\bigoplus_{j=2}^{r}\cF_j\hookrightarrow T_Y$ since $f$ is faithfully flat. This proves (c).
\end{proof}

\begin{proposition}\label{Proposition:5}
    In the setup of Proposition~\ref{Proposition:4},  $\cU:=\cW_2\otimes\cdots\otimes \cW_r$ is an ample vector bundle on $Y$.
\end{proposition}

Before we prove Proposition~\ref{Proposition:5}, we require the following remark and lemma.

\begin{remark}\label{Rem:Druel}
     The following consequence of the proofs of \cite[Lemma~5 and Proposition~7]{DruelParis13} will be useful. Let $S$ be a smooth complex projective variety, $\cE$ an ample vector bundle of rank $s\geq q>0$. Let $h:S\rightarrow C$ be a surjective smooth morphism onto a smooth connected curve $C$ with integral fibers. Let $N$ be a nef vector bundle of positive rank on $C$. Assume that $h$ is a smooth projective bundle or a smooth quadric bundle. Let $A$ be an ample line bundle on $S$, whose restriction to any fiber $F$ of $h:S\rightarrow C$ is $A|_F\simeq\OO_F(s)$ with $s\geq q\geq1$. Then:
    \begin{equation}\label{Eq:DPvanishing}
 H^0\bigl(S,T_{S/C}^{\otimes q}\otimes A^{-1}\otimes h^*N^\vee\bigr)=0.
\end{equation}
Although the statements of \cite[Lemma~5 and Proposition~7]{DruelParis13} use $A=\det\cE$ for an ample bundle $\mathcal E$ of rank $s\geq q$, the only property of $\cE$ used is the ampleness of $\det\mathcal E$ and its degree on the fibers $F$. Thus, \eqref{Eq:DPvanishing} follows by a similar proof.
 \end{remark}

\begin{lemma}\label{lem:calculation}
Let $f:X\rightarrow Y$ be a smooth Brauer-Severi scheme or smooth quadric bundle with fiber dimension $d\geq 1$. Let $\cL$ be a nef line bundle such that $\cL|_F\simeq\OO_F(\ell)$ on every fiber $F$ of $f:X\rightarrow Y$, for a fixed $\ell>0$. For integers $a,q\geq1$ satisfying $s:=a\ell\geq q$, let
\[
 \cK_q=f_*\bigl(T_{X/Y}^{\otimes q}\otimes \cL^{-\otimes a}\bigr).
\]
Then:
\begin{enumerate}
    \item[(i)] $\cK_q$ is locally free, whose formation commutes with base change and $\cK_q^\vee$ is nef.
    \item[(ii)]  For a finite morphism $\nu:C\rightarrow Y$ from a smooth projective curve, if the pullback $\cL|_S$ on $S:=X\times_Y C$ is ample and $\cK_q\neq0$, then $\nu^*\cK_q^\vee$ is ample.
\end{enumerate} 
\end{lemma}

\begin{proof}
  For any fiber $F$ of $f:X\rightarrow Y$, we have $(T_{X/Y}^{\otimes q}\otimes\cL^{-\otimes a})|_F\simeq T_F^{\otimes q}\otimes\OO_F(-a\ell)$, whereby $\dim H^0(F,T_F^{\otimes q}\otimes\OO_F(-a\ell))$ is constant for all fibers $F$ of $f:X\rightarrow Y$. By Grauert's theorem \cite[Chapter III, Corollary~12.9]{Hartshorne77} it follows that $\cK_q$ is a vector bundle on $Y$. Furthermore, the formation of $\cK_q$ commutes with base change. 
 
 We now prove that $\cK_q^\vee$ is nef by contradiction. Thus, by \cite[Proposition~6.1.18(i)]{PosII}, there is a finite morphism $\gamma:C\to Y$ from a smooth projective integral curve such that there is a surjection $\gamma^*\cK_q^\vee\twoheadrightarrow M^{-1}$ to a line bundle $M$ on $C$ with $\deg M>0$. Dualizing yields an injection $M\hookrightarrow\gamma^*\cK_q$. Let \[
 h:S=X\times_Y C\rightarrow C,\qquad g:S\rightarrow X, \qquad \cL|_S=g^*\cL.
\]
Using the base-change isomorphism $\gamma^*\cK_q\simeq h_*(T_{S/C}^{\otimes q}\otimes \cL|_S^{-\otimes a})$ along with projection formula turns the injection $M\hookrightarrow\gamma^*\cK_q$ into a nonzero section
\begin{equation}\label{eq:414-positive-line-section}
 0\ne\tau\in H^0\!\left(S,T_{S/C}^{\otimes q}\otimes \cL|_S^{-\otimes a} \otimes h^*M^{-1}\right).
\end{equation}    
We claim that $\cL|_S^{\otimes a}\otimes h^\ast M$ is ample on $S$. Note that $\cL|_S$ is $h$-ample by assumption $\cL|_F\simeq\OO_F(\ell)$ for $\ell>0$ on every fiber $F$ of $f:X\rightarrow Y$. Since $\deg M>0$, $M$ is ample on $C$, whereby $\cL|_S^{\otimes a}\otimes h^\ast M^{\otimes N}$ is ample on $S$ for $N\gg0$ by \cite[Proposition~1.7.10]{PosI}. Thus, 
\[(\cL|_S^{\otimes a}\otimes h^\ast M)^{\otimes N}\simeq (\cL|_S^{\otimes a}\otimes h^\ast M^{\otimes N})\otimes \cL|_S^{\otimes a(N-1)}\]
is the tensor product of an ample and nef line bundle (since $\cL$ is nef by assumption), whereby $\cL|_S^{\otimes a}\otimes h^\ast M$ is ample by \cite[Corollary~1.4.10]{PosI}. Furthermore, the degree of $\cL|_S^{\otimes a}\otimes h^\ast M$ restricted to any fiber of $h:S\rightarrow C$ is $s:=a\ell$, which is at least $q$ by hypothesis. Thus, applying \eqref{Eq:DPvanishing} using $A=\cL|_S^{\otimes a}\otimes h^\ast M$ and $N=\OO_C$, we have $H^0\!\left(S,T_{S/C}^{\otimes q}\otimes \cL|_S^{-\otimes a} \otimes h^*M^{-1}\right)=0$ since the general fiber of $h$ is projective space or hyperquadric. But this contradicts \eqref{eq:414-positive-line-section}, which completes the proof that $\cK^\vee_q$ is nef and completes proof of (i).

We now prove (ii). Using the same notation as above, suppose that $\cL|_S$ is ample and $\cK_q\ne0$. Let
\[
 \cE_q:=\nu^*\cK_q^\vee,\qquad
 \cF_q:=T_{S/C}^{\otimes q}\otimes\cL|_S^{-\otimes a}.
\]
The bundle $\cE_q$ is nef by the preceding argument, and
base change yields $\cE_q^\vee\simeq h_*\cF_q$. Assume $\cE_q$ is not ample, whereby Hartshorne's criterion \cite[Theorem~2.4]{Har71} yields a positive-rank vector-bundle quotient $\cE_q\twoheadrightarrow\cQ$ with $\deg\cQ\leq0$. Since $\cQ$ is nef, we must have $\deg\cQ=0$ and therefore, $\cQ$ is numerically flat. Dualizing the quotient and applying adjunction yields a non-zero element of 
\begin{equation}\label{Eq:nonvanishing}
\Hom_C(\cQ^\vee,h_*\cF_q)\simeq
 H^0\!\left(S,T_{S/C}^{\otimes q}\otimes\cL|_S^{-a}\otimes h^*\cQ\right)\neq 0.
\end{equation}
But $\cL|_S^{\otimes a}$ is ample by assumption and restricts to
$\OO_F(a\ell)$ on each fiber of $h:S\rightarrow C$, where $a\ell\geq q$. Since $\cQ^\vee$ is nef, applying \eqref{Eq:DPvanishing} using $A=\cL|_S^{\otimes a}$ and $N=\cQ^\vee$, yields $H^0\!\left(S,T_{S/C}^{\otimes q}\otimes \cL|_S^{-\otimes a} \otimes h^*\cQ\right)=0$ since the general fiber of $h$ is projective space or hyperquadric. However this contradicts \eqref{Eq:nonvanishing}, whereby $\nu^*\cK_q^\vee$ is ample.
\end{proof}

The following proposition will be useful in the proof of both Theorems~\ref{thm:Main2} and \ref{thm:Main3}.

\begin{proposition}\label{prop:relative-ample-base}
Let $f:X\rightarrow Y$ be a smooth Brauer-Severi scheme or smooth quadric bundle between smooth projective complex varieties, with connected fibers of dimension $d>0$. Let $\cV$ be a nef vector bundle of rank $s>0$ on $X$, and put $\cL:=\det\cV$. Suppose that $\cL|_F\simeq\OO_F(\ell)$ on every fiber, for a fixed integer $\ell>0$. Assume that there are integers $a,q\geq1$ with
$a\ell\geq q$, a numerically flat vector bundle $\cP$ on $Y$, and
a nonzero morphism
\begin{equation}\label{eq:relative-twisted-tensor}
 \psi:\cL^{\otimes a}\rightarrow
 T_{X/Y}^{\otimes q}\otimes f^*\cP.
\end{equation}
Let $\cU$ be a nef vector bundle of positive rank on $Y$. If
$\cB:=\cV\otimes f^*\cU$ is ample, then $\cU$ is ample.
\end{proposition}

\begin{proof}
If $Y$ is a point, then $\cU$ is ample since every vector bundle on a point is ample. So assume $\dim Y>0$ and let $\cK_q=f_*(T_{X/Y}^{\otimes q}\otimes\cL^{-\otimes a})$ as in Lemma~\ref{lem:calculation}. By the projection formula, \eqref{eq:relative-twisted-tensor} yields a nonzero section $\sigma\in H^0(Y,\cK_q\otimes\cP)$. In particular, $\cK_q\neq0$ and is therefore by Lemma~\ref{lem:calculation}(i), $\cK_q$ is locally free with $\cK_q^\vee$ being nef. Since $\cP$ is numerically flat, it follows that
\[
 \cH:=\cK_q^\vee\otimes\cP^\vee
\]
is nef as well. A nonzero section of the dual of a nef vector bundle
is nowhere vanishing by \cite[Proposition~1.16]{DPS94}. Thus $\sigma$
gives a quotient
\begin{equation}\label{eq:relative-trivial-quotient}
 \cH\twoheadrightarrow\OO_Y.
\end{equation}

Now fix an ample line bundle $H$ on $Y$ and choose an integer $M\geq1$ such that
\begin{equation}\label{eq:relative-uniform-ample}
 \Sym^M\cB\otimes f^*H^{-1}\quad\text{is ample}.
\end{equation}
Indeed, for $M\gg0$, the bundle $\Sym^{M-1}\cB\otimes f^*H^{-1}$ is globally generated by \cite[Section~2, Definition]{Hartshorne66}. Since $\cB$ is ample, the tensor product $\cB\otimes\Sym^{M-1}\cB\otimes f^*H^{-1}$ is ample and surjects onto the bundle in \eqref{eq:relative-uniform-ample}, thereby making it ample.

Let $\nu:C\rightarrow Y$ be a finite nonconstant morphism from a
smooth projective integral curve, and let
$\nu^*\cU\twoheadrightarrow N$ be a quotient line bundle. Let
\[
 S=X\times_Y C,\qquad h:S\rightarrow C,\qquad g:S\rightarrow X,
\]
whereby the morphism $g$ is finite. Since $f\circ g=\nu\circ h$, the quotient
$\nu^*\cU\twoheadrightarrow N$ induces
\[
 g^*\cB\simeq g^*\cV\otimes h^*\nu^*\cU
 \twoheadrightarrow g^*\cV\otimes h^*N.
\]
Symmetric powers preserve surjections and commute with pullback.
Consequently, taking the $M$-th symmetric power and tensoring with
$h^*\nu^*H^{-1}$ gives
\begin{align*}
 g^*\!\left(\Sym^M\cB\otimes f^*H^{-1}\right)
 &\twoheadrightarrow
 \Sym^M(g^*\cV\otimes h^*N)\otimes h^*\nu^*H^{-1}\\
 &\simeq
 \Sym^M(g^*\cV)\otimes
 h^*\!\left(N^{\otimes M}\otimes\nu^*H^{-1}\right)
 =:\cB_C,
\end{align*}
where the last isomorphism follows from
Lemma~\ref{lem:Schur-basic}(ii), since $N$ is a line bundle.
The source is ample by \eqref{eq:relative-uniform-ample} and
finiteness of $g$, so its locally free quotient $\cB_C$ is ample. Let $r_M=\rk(\Sym^M\cV)$ and
$e=Mr_M/s\in\ZZ_{>0}$. By Lemma~\ref{lem:Schur-basic}, we obtain
\begin{equation}\label{eq:relative-curve-determinant}
 \det\cB_C=(g^*\cL)^{\otimes e}\otimes
 h^*\bigl(N^{\otimes Mr_M}\otimes(\nu^*H)^{-\otimes r_M}\bigr).
\end{equation}
If $M\deg N<\deg\nu^*H$, then $(N^{\otimes Mr_M}\otimes(\nu^*H)^{-\otimes r_M})^\vee$ has positive degree on $C$, and its pullback along $h$ is nef. Note that $\det\cB_C$ is ample since $\cB_C$ is, whereby \eqref{eq:relative-curve-determinant} implies that $(g^*\cL)^{\otimes e}$, and hence $g^*\cL$, is ample. By Lemma~\ref{lem:calculation}(ii) we then obtain that
$\nu^*\cK_q^\vee$ is ample. The bundle $\nu^*\cP^\vee$ is nef, so
$\nu^*\cH=\nu^*\cK_q^\vee\otimes\nu^*\cP^\vee$ would also be
ample by \cite[Theorem~6.2.12]{PosII}. This contradicts the quotient
$\nu^*\cH\twoheadrightarrow\OO_C$ obtained from
\eqref{eq:relative-trivial-quotient}. Therefore we must have
\begin{equation}\label{eq:relative-uniform-bound}
 M\deg N\geq\deg\nu^*H
\end{equation}
for every such curve and quotient, with the same $M$ and $H$. Let $\pi:\PP_Y(\cU)\rightarrow Y$ be the projectivization of $\cU$, whereby using the quotient convention, we have the universal quotient $\pi^\ast\cU\twoheadrightarrow\xi$, where $\xi=\OO_{\PP_Y(\cU)}(1)$. By definition, $\cU$ is ample if $\xi$ is ample. If $j:C\rightarrow \PP_Y(\cU)$ is the normalization of any irreducible curve $\Gamma\subset\PP_Y(\cU)$ not contained in a fiber of $\pi$, then $\nu:=\pi\circ j:C\rightarrow Y$ is finite non-constant morphism. Pulling back the universal quotient along $j$ yields the quotient $\nu^\ast\cU\twoheadrightarrow N:=j^\ast\xi$. Let $B:=\xi^{\otimes M}\otimes\pi^\ast H^{-1}$, where $M, H$ are as above. Then
\[B\cdot\Gamma:=\deg_C(j^\ast B)=\deg_C(N^{\otimes M}\otimes(\nu^\ast H)^{-1})=M\deg N-\deg\nu^\ast H\geq 0,\] by \eqref{eq:relative-uniform-bound}. On curves $\Gamma\subseteq \PP_Y(\cU)$ contained in a fiber $F$ of $\pi$, $\xi$ is ample and $H$ is trivial, whereby $B$ is ample and therefore $B.\Gamma>0$. Thus, by the curve criterion for nefness
\[
 B:=\xi^{\otimes M}\otimes\pi^*H^{-1}\quad\text{is nef}.
\]
Since $\xi$ is $\pi$-ample and $H$ is ample on the base $Y$, by \cite[Proposition~1.7.10]{PosI}, there exist $k\gg0$ such that $\xi\otimes\pi^*H^k$ is ample. Then
\[
 \xi^{\otimes(Mk+1)}=
 (\xi\otimes\pi^*H^k)\otimes
 (\xi^{\otimes M}\otimes\pi^*H^{-1})^{\otimes k}
\]
is ample, since it is the tensor product of an ample and a nef line bundle.
Thus $\xi$, and consequently $\cU$, is ample.

\end{proof}

\noindent
We can now prove Proposition~\ref{Proposition:5}.

\begin{proof}[Proof of Proposition~\ref{Proposition:5}]
The statement is tautological if $Y$ is a point, whereby we assume $\dim Y\geq 1$. By Proposition~\ref{Proposition:4}(a), each $\cW_j$ is nef. Therefore $\cU:=\cW_2\otimes\cdots\otimes\cW_r$ is nef by \cite[Theorem~6.2.12]{PosII}. We verify the hypotheses of Proposition~\ref{prop:relative-ample-base}.
\begin{itemize}
    \item By hypothesis $\cV_1$ is nef and  $\cA=\cV_1\otimes\dots\cV_r=\cV_1\otimes f^*\cU$ is ample.
    \item By Lemma~\ref{lem:fib-restriction}(3), $\cV_1|_F$ is ample on every fiber $F$ of $f:X\rightarrow Y$ and by Lemma~\ref{lem:fib-characterization}, either all fibers of $f:X\rightarrow Y$ are $\PP^d$ or all fibers are $Q_d$ for some $d\geq 1$.
    \item Let $\cL:=\det\cV_1$, whereby proof of Lemma~\ref{lem:fib-characterization} via Theorem~\ref{thm:generalizedKovacs} yields $\cL|_F\simeq \OO_{F}(\ell)$ for general fiber $F$ of $f:X\rightarrow Y$, where $\ell\geq r_1$ is a fixed integer given by the single-factor classification. We claim this can be extended to all fibers $F$: let $\iota=d+1$ for $F\simeq\PP^d$ and $\iota=d$ for $F\simeq Q_d$, whereby $\cL^{\otimes\iota}\otimes\omega_{X/Y}^{\otimes \ell}$ restricts to $\OO_F$ for the general fiber $F$. By separatedness of the relative Picard scheme for projective flat morphism with geometrically integral fibers (see \cite[Théorème~3.1]{Grothendieck62}), we obtain $(\cL^{\otimes\iota}\otimes\omega_{X/Y}^{\otimes \ell})|_F\simeq\OO_F$ or equivalently, $\cL|_F\simeq \OO_F(\ell)$ (using torsion freeness of $\Pic(F)$ by Lemma~\ref{lem:NF-RC}(iii)) for all fibers $F$ of $f:X\rightarrow Y$.
\end{itemize}
Let $a=q_{\lambda_1}(r_1)>0$, where $r_1=\rk\cV_1$ and $R= r_1a$.
Applying Proposition~\ref{prop:Schur-det-to-tensor} to $\theta_1:\bS_{\lambda_1}\cV_1\rightarrow \bS_{\lambda_1}\cE_1$, yields
the non-zero morphism $\varphi:\cL^{\otimes a}\rightarrow\cE_1^{\otimes R}$. Analogous to \eqref{Eq:gradeddecomposition}, the tensor filtration $F^\bullet$ of $\cE_1^{\otimes R}$ induced by the associated short exact sequence \eqref{Equation:ses1}, results in the graded pieces
\[
 \operatorname{gr}_F^t(\cE_1^{\otimes R})
 \simeq
 \bigl(T_{X/Y}^{\otimes t}\otimes f^*(\cG^{\otimes(R-t)})\bigr)
       ^{\oplus\binom Rt},\qquad 0\leq t\leq R.
\]
Analogous to the argument leading to \eqref{eq:Schur-tangent-tensor}, choose the largest $q$ such that $\varphi: \cL^{\otimes a}\rightarrow\cE_1^{\otimes R}$ factors through $F^q$. Thus, the induced map $\cL^{\otimes a}\rightarrow\operatorname{gr}_F^q(\cE_1^{\otimes R})$ is nonzero, since otherwise $\varphi$ would factor through $F^{q+1}$. Projecting to a suitable summand of $\operatorname{gr}_F^q(\cE_1^{\otimes R})$ yields a non-zero map, 
\begin{equation}\label{Eq:map2}
 \psi:\cL^{\otimes a}\rightarrow T_{X/Y}^{\otimes q}\otimes f^*\cP_q,
\end{equation}
where $\cP_q:=\cG^{\otimes(R-q)}$ is a numerically flat vector bundle on $Y$. Since $\cL|_F$ is ample on the fibers $F$ of $f:X\rightarrow Y$, we have $f_*\cL^{-\otimes a}=0$. Thus, we cannot have $q=0$ in \eqref{Eq:map2} as the projection formula yields
\[
 \Hom_X(\cL^{\otimes a},f^*(\cG^{\otimes R}))
 \simeq H^0(Y,f_*\cL^{-\otimes a}\otimes\cG^{\otimes R})=0.
\]
Hence $1\leq q\leq R$. With the convention $\cG^{\otimes0}=\OO_Y$, the bundle
$\cP_q=\cG^{\otimes(R-q)}$ is numerically flat, and
\[
 a\ell\geq ar_1=R\geq q.
\]
Thus \eqref{Eq:map2} supplies the nonzero morphism required in Proposition~\ref{prop:relative-ample-base}, with
$\cV=\cV_1$ and $\cP=\cP_q$. Since $\cA=\cV_1\otimes f^*\cU$ is ample, that proposition proves that $\cU$ is ample.
\end{proof}

 We can now complete the inductive step for the proof of Theorem~\ref{thm:productcharacterization}.

 \begin{proof}[Proof of Theorem~\ref{thm:productcharacterization}]
    Assume $r\ge 2$ and that Theorem~\ref{thm:productcharacterization} has been proved for $r-1$ almost nef regular foliations of rank $n_2, n_3,\dots, n_r\geq 1$. Applying Propositions~\ref{Proposition:4} and~\ref{Proposition:5}, after possibly relabeling, there exists either a smooth Brauer-Severi scheme or a smooth quadric bundle $f:X\rightarrow Y$. Furthermore for all $j\geq 2$, we have:
\begin{enumerate}
    \item nef vector bundles $\cW_2,\dots,\cW_r$ on $Y$ such that $\cV_j\cong f^*\cW_j$ and $\cW_2\otimes \dots\otimes \cW_r$ ample on $Y$,
    \item almost nef regular foliations $\cF_j\subseteq T_Y$ such that $\cE_j\cong f^\ast\cF_j$ with subbundle inclusion $\cF_2\oplus\cdots\oplus \cF_r\hookrightarrow T_Y$ and,
    \item injective morphisms $\bS_{\lambda_j}\cW_j\rightarrow \bS_{\lambda_j}\cF_j$ of coherent sheaves for all $j\geq 2$.
\end{enumerate}  
Hence the induction hypothesis applies to the tuple $\bigl(Y;\cF_2, \dots, \cF_r; \cW_2,\dots,\cW_r\bigr)$, whereby we conclude that $Y\cong Z_2\times\dots Z_r$ such that for all $2\leq j\leq r$, $\cF_j\simeq p_j^\ast T_{Z_j}$ where $p_j: Y\cong Z_2\times\dots\times Z_r\rightarrow Z_j$ is the natural projection map. Furthermore, for all $2\leq j\leq r$, we have $Z_j\simeq \PP^{n_j}$ or a smooth hyperquadric $Q_{n_j}\subseteq \PP^{n_j+1}$ of dimension at least $2$ such that there are ample bundles $\cT_j$ on $Z_j$ satisfying $\cW_j\simeq p_j^\ast \cT_j$ and admitting injective morphism $\bS_{\lambda_j}\cT_j\rightarrow\bS_{\lambda_j}T_{Z_j}$ of coherent sheaves, where $\cT_j$ are classified as in Theorem~\ref{thm:productcharacterization}. 

However, by Proposition~\ref{Proposition:4}(c), there is a numerically flat regular foliation $\mathcal{G}$ on $Y$ and a subbundle inclusion $\mathcal{G}\oplus \cF_2\oplus\dots \cF_r\hookrightarrow T_Y$, which is only possible if $\mathcal{G}=0$. Then the corresponding short exact sequence \eqref{Equation:ses1} applied to $f:X\rightarrow Y$ yields $\cE_1\cong T_{X/Y}$, whereby the fibers of $f:X\rightarrow Y$ have dimension equal to $\rk\cE_1=n_1$. Thus, $\dim X=n_1+\dim Y=n_1+n_2+\dots+n_r$, which is also the rank of $T_X$. Thus, by comparing ranks it follows that the subbundle inclusion $\bigoplus_{i=1}^{r}\cE_i\hookrightarrow T_X$ in hypothesis of Theorem~\ref{thm:productcharacterization} is an isomorphism. In other words, $f:X\rightarrow Y\cong Z_2\times \dots\times Z_r$ admits a connection, that is $T_X\cong T_{X/Y}\oplus V$, where $V=\bigoplus_{j=2}^r\cE_j\simeq f^*T_Y$. Since $Y$ is rationally connected, by \cite[Lemma 3.18]{Horing07}, we have $X\cong Y\times F$, where $F$ is a general fiber of $f:X\rightarrow Y$, and $f$ is the projection onto $Y$. The cited lemma also gives integrability of $V$, and the product isomorphism respects $T_{X/Y}$ and $V$; the individual $\cE_j$ are the horizontal lifts of $\cF_j$. By Lemma~\ref{lem:fib-characterization}, it follows that $F\cong \PP^{n_1}$ or $Q_{n_1}$, whereby we obtain $X\cong Z_1\times\dots\times Z_r$ where  $Z_i\simeq \PP^{n_i}$ or a smooth hyperquadric $Q_{n_i}\subseteq \PP^{n_i+1}$ of dimension at least $2$. Let $\pi_i:X\cong Z_1\times\dots\times Z_r\rightarrow Z_i$ be the natural projection onto the $i^{th}$ factor. Then $\pi_j=p_j\circ f$ for all $j\geq 2$ and thus, $\cE_j\cong f^\ast\cF_j\cong \pi_j^\ast T_{Z_j}$ for $j\geq 2$. For $j=1$, $\cE_1\cong T_{X/Y}\cong \pi_1^\ast T_{Z_1}$. Similarly for $j\geq 2$, $\cV_j\cong f^\ast \cW_j\cong \pi_j^\ast\cT_j$, where $\cT_j$ are ample vector bundles on $Z_j$, classified as in Theorem~\ref{thm:productcharacterization} by induction hypothesis. Therefore, what remains is to prove that $\cV_1\simeq \pi_1^\ast\cT_1$ for some ample vector bundle $\cT_1$ on $Z_1$ and to classify $\cT_1$.

By hypothesis, we have the following injection of coherent sheaves, where we use $\cE_1\simeq\pi_1^\ast T_{Z_1}$ and Proposition~\ref{prop:Schur-functorial}(iii):
\begin{equation}\label{eq:V1-injection}
 \theta_1:\bS_{\lambda_1}\cV_1\rightarrow\bS_{\lambda_1}(\pi_1^*T_{Z_1})
     \simeq\pi_1^*(\bS_{\lambda_1}T_{Z_1}).
\end{equation}
Since $\cV_1$ is nef by hypothesis and the fibers of $\pi_1:X\rightarrow Z_1$ are rationally connected, one can apply Lemma~\ref{lem:fiber-trivial} to obtain the unique vector bundle $\cT_1:=\pi_{1,*}\cV_1$ on $Z_1$ such that $\pi_1^\ast\cT_1\simeq\cV_1$ on $X$. For any $y\in Y$ let $\iota_y:Z_1\to Z_1\times Y$ be the section of $\pi_1$ defined by
$z\mapsto(z,y)$. Since $\pi_1\circ\iota_y=\id_{Z_1}$, we obtain $\cT_1\simeq\iota_y^*\pi_1^\ast\cT_1\simeq\iota_y^\ast\cV_1$.
Substituting $\pi_1^\ast\cT_1\simeq\cV_1$ in
\eqref{eq:V1-injection} we obtain the injection $\theta_1:\pi_1^\ast(\bS_{\lambda_1}\cT_1)\rightarrow\pi_1^*(\bS_{\lambda_1}T_{Z_1})$, which by $\pi_{1,*}\OO_X=\OO_{Z_1}$, adjunction, and the projection formula descends uniquely. By faithful flatness of $\pi_1$, the descended morphism is injective. Thus we obtain the unique injection of coherent sheaves on $Z_1$:
\begin{equation}\label{eq:T1-injection}
\overline\theta_1:\bS_{\lambda_1}\cT_1\rightarrow\bS_{\lambda_1}T_{Z_1},
 \qquad \pi_1^*\overline\theta_1=\theta_1.
\end{equation}
We now show that $\cT_1$ is an ample bundle on $Z_1$. Restricting the ample bundle $\cA=\cV_1\otimes\dots\otimes\cV_r$ on $X\simeq Z_1\times Y$ to the closed subvariety $Z_1\times\{y\}$ for any $y\in Y$, we see that:
\[\iota_y^*\cA
  \simeq\cT_1\otimes_\CC\left(\bigotimes_{j=2}^r\cT_{j,p_j(y)}\right)
  \simeq\cT_1^{\oplus M},\]
is an ample bundle on $Z_1$, where $M:=\prod_{j=2}^{r}r_j$ and $\cT_{j,p_j(y)}$ is the fiber of the bundle $\cT_j$ at $p_j(y)\in Z_j$. This is because $\pi_j\circ\iota_y:Z_1\rightarrow Z_j$ is the constant map to $p_j(y)$ for all $j\geq 2$. This proves that $\cT_1$ is ample on $Z_1$. Thus, applying Theorem~\ref{thm:generalizedKovacs} to \eqref{eq:T1-injection}, we obtain the desired classification of $\cT_1$. Every classified $\cT_i$ is globally generated (using the Euler sequence in the tangent-bundle case), so its pullback $\cV_i\simeq\pi_i^*\cT_i$ is globally generated. This completes the proof of Theorem~\ref{thm:productcharacterization} by induction on the number $r$ of almost nef regular foliations.
 \end{proof}

 \section{Characterizing products of projective spaces and hyperquadrics among smooth weak Fano varieties}\label{sec:CharacterizeProduct2}

 The main result of this section is the following.

 \begin{theorem}\label{thm:productcharacterization2}
Let $X$ be a smooth weak Fano complex variety of dimension $n>0$, and
let $r\geq1$. For $1\leq i\leq r$, let $\cE_i\subseteq T_X$ be a
vector subbundle of rank $n_i\geq1$, let $\cV_i$ be a nef vector
bundle of rank $r_i\geq1$, and let $\lambda_i$ be a nonzero partition
of length $\ell(\lambda_i)\leq\min\{r_i,n_i\}$. Assume that:
\begin{enumerate}
\item the natural morphism
$\cF:=\bigoplus_{i=1}^r\cE_i\rightarrow T_X$ is a subbundle inclusion with integrable image,
\item for every $1\leq i\leq r$, there is an injective morphism
$\theta_i:\bS_{\lambda_i}\cV_i\hookrightarrow
\bS_{\lambda_i}\cE_i$ of coherent $\OO_X$-modules and,
\item the vector bundle $\cA:=\bigotimes_{i=1}^r\cV_i$ is ample on $X$.
\end{enumerate}
Then $\cF=T_X$ and, possibly after reordering, there is an isomorphism
\begin{equation}\label{eq:wf-main-product2}
 X\simeq\prod_{i=1}^rZ_i,\qquad
 n=\sum_{i=1}^r n_i,\qquad
 \cE_i\simeq\pr_i^*T_{Z_i}\subseteq T_X,
\end{equation}
where each $Z_i$ is $\PP^{n_i}$ or a smooth hyperquadric
$Q_{n_i}\subset\PP^{n_i+1}$ with $n_i\geq2$. In particular, $X$ is Fano. Moreover, $r_i\leq n_i$ and there are ample vector bundles $\cT_i$ on $Z_i$
such that $\cV_i\simeq\pr_i^*\cT_i$ and $\cT_i$ admits an injective morphism $\bS_{\lambda_i}\cT_i\rightarrow\bS_{\lambda_i}T_{Z_i}$ of coherent sheaves. The bundles $\cT_i$ and the partitions $\lambda_i$ have precisely the
possibilities listed in Theorem~\ref{thm:Main2}\textup{(a)--(d)}.
\end{theorem}

In Section~\ref{subsec:toric-foliations}, we also give a simpler characterization of products of projective spaces when $X$ is further assumed to be toric.

\subsection{A few useful lemmas}

\begin{lemma}\label{lem:wf-targetfiber}
Let $f:X\rightarrow Y$ be a smooth projective surjective morphism with
connected fibers between smooth projective complex varieties. If $X$
is weak Fano, then $Y$ and every positive-dimensional fiber of $f$ are
weak Fano as well. In particular, these varieties are all rationally connected.
\end{lemma}

\begin{proof}
By \cite[Theorem~1.1]{FG12}, $Y$ is weak Fano. To prove that the fibers of $f$ are weak Fano, let $L:=\OO_X(-K_X)$ and $d>0$ be the dimension of the fibers. By adjunction, $L|_{X_y}\simeq\OO_{X_y}(-K_{X_y})$, which is nef for
every $y\in Y$. Since $-K_X$ is big, by Kodaira's lemma \cite[Proposition~2.2.6 and Corollary~2.2.7(iii)]{PosI}, there exist an integer $m_0>0$, an ample Cartier divisor $A$, and an effective Cartier divisor $E$ on $X$ such that
\begin{equation}\label{Eq:bigkodaira}
 -m_0K_X\sim A+E.
\end{equation}
Since $f$ is smooth, it is an open morphism and hence $Y^\circ:=f(X\setminus\operatorname{Supp}E)$ is a dense open subset of $Y$. For $y\in Y^\circ$, $E|_{X_y}$ is effective and $A|_{X_y}$
is ample. Thus, by pulling back \eqref{Eq:bigkodaira} to the fiber $X_y$, \cite[Corollary~2.2.7(iii)]{PosI} implies that $L|_{X_y}$ is big.

For every integer $m\geq0$, the invertible sheaf $L^{\otimes m}$
is flat over $Y$ and thus, by \cite[Lemma~36.32.2]{stacks-project} $\chi(X_y,L^{\otimes m}|_{X_y})$ is independent of $y$. By asymptotic Riemann--Roch \cite[Theorem~1.1.24]{PosI}, we obtain
\[
 \chi(X_y,L^{\otimes m}|_{X_y})
 =
 \frac{(L|_{X_y})^d}{d!}\,m^d+O(m^{d-1}),
\]
whereby the top intersection number $(L|_{X_y})^d$ is independent of $y$. It is positive on a general fiber, since $L|_{X_y}$ is nef and big for $y\in Y^\circ$. Consequently $(L|_{X_y})^d$ is positive for every fiber. By the numerical criterion for bigness of nef divisors \cite[Theorem~2.2.16]{PosI}, we therefore obtain bigness of $L|_{X_y}$ for every $y\in Y$. This proves that $\OO_{X_y}(-K_{X_y})$ is big and nef for all $y\in Y$, thereby establishing that every fiber of $f$ is weak Fano. By \cite[Theorem~1]{Zhang06}, every weak Fano variety is rationally connected, which establishes the last claim of the lemma.
\end{proof}

\begin{definition}[{\cite[Definition~1.1]{AraujoMassarenti16}, \cite[Definition~2.14]{Horing26}}]
Let $X$ be a normal projective $\QQ$-factorial variety. We say that $X$ is log Fano or Fano type if there exists an effective $\QQ$-divisor $\Delta$ such that $-(K_X+\Delta)$ is ample, and the pair $(X, \Delta)$ is klt. \end{definition}

\begin{lemma}\label{lem:wf-partial-integrability}
Let $X$ be a smooth weak Fano variety and let $\cF=\bigoplus_{i=1}^r\cE_i\subseteq T_X$ be a regular foliation. Then every partial sum $\cE_I:=\bigoplus_{i\in I}\cE_i$ defines a regular foliation, and thus, a subfoliation $\cE_I\subseteq \cF\subseteq T_X$, for
any $I\subseteq\{1,\ldots,r\}$.
\end{lemma}

\begin{proof}
If $\cF=0$, then every $\cE_i=0$ and the assertion is immediate.
We may therefore assume that $\rk\cF>0$. By Druel's theorem \cite[Theorem~1.3]{Druel17} and its proof, there is a smooth projective morphism $p:X\rightarrow B$ with connected fibers, onto a smooth projective variety $B$ such that $\cF=T_{X/B}$. In particular, the fibers of $p$ are positive dimensional.

By \cite[Lemma~2.5]{AraujoMassarenti16}, every smooth weak Fano variety is log Fano. Thus, by Lemma~\ref{lem:wf-targetfiber}, every fiber $W$ of $p:X\rightarrow B$ is weak Fano, and hence log Fano. For a nonempty proper subset $I\subseteq\{1,\dots, r\}$, restricting $\cF=T_{X/B}$ to $W$ yields the decomposition
\[
 T_{X/B}|_W=T_W=\cF|_W=\cE_I|_W\oplus\cE_{I^c}|_W,
\]
where $I^c:=\{1,\dots, r\}\setminus I$. Thus, $\cE_I|_W$ defines a regular foliation on $W$ by \cite[Theorem~1.1]{Horing26}. Now recall that the Lie bracket induces an $\OO_X$-linear morphism $\bigwedge^2\cE_I\to T_X/\cE_I$ called the Frobenius tensor, whose vanishing is equivalent to integrability of $\cE_I$ (see \cite[\S~2.B]{Horing07}). Since $\cF$ is integrable, this morphism factors as
\[
 \beta_I:\bigwedge\nolimits^2\cE_I
       \rightarrow\cF/\cE_I,
 \qquad u\wedge v\longmapsto[u,v]\bmod\cE_I.
\]
For a fiber $W$ of $p$, we have $(\cF/\cE_I)|_W=(T_{X/B}/\cE_I)|_W=T_W/(\cE_I|_W)$. Since $\cE_I\subseteq T_{X/B}$, its local sections are $\OO_B$-derivations, i.e., derivations annihilating functions pulled back from $B$, see \cite[Lemma~29.33.2]{stacks-project}.
Let $\mathcal{J}_W\subseteq \OO_X$ be the ideal sheaf of a fiber $W$ of $p$. This ideal is locally generated by pullbacks of functions on $B$ vanishing at $p(W)$, so the Leibniz rule implies $u(\mathcal{J}_W)\subseteq\mathcal{J}_W$ for every local section $u$ of $\cE_I$. Thus $u$ induces a derivation $\bar u$ on $\OO_W$.
For two such sections $u,v$, the definition of the Lie bracket as a commutator yields $\overline{[u,v]}=[\bar u,\bar v]$. Consequently, the restriction $\beta_I|_W$ is the Frobenius tensor of $\cE_I|_W$, which vanishes by the integrability of $\cE_I|_W$ established above. Thus $\beta_I$ vanishes at every closed point of $X$, and hence, $\beta_I=0$ since $X$ is reduced and $\beta_I$ is a morphism of vector bundles. This proves that $\cE_I$ is integrable, or equivalently, a foliation on $X$. Since $\cE_I\hookrightarrow \cF\hookrightarrow T_X$ is an inclusion of subbundles, it follows that $\cE_I$ is a regular foliation on $X$ and a subfoliation of $\cF$. The cases $I=\varnothing$ and $I=\{1,\ldots,r\}$ are immediate.
\end{proof}

The following lemma gives the descent
step corresponding to Proposition~\ref{Prop:almost-nef-descent}.

\begin{lemma}\label{lem:wf-transverse-descent}
Let $f:X\rightarrow Y$ be a smooth projective surjective morphism with
connected fibers between smooth projective complex varieties. Let
$\cE\subseteq T_X$ be a subbundle such that
$T_{X/Y}\oplus\cE\rightarrow T_X$ is a subbundle inclusion. Suppose
that
\begin{equation}\label{eq:wf-mixed-bracket}
 [T_{X/Y},\cE]\subseteq T_{X/Y}\oplus\cE.
\end{equation}
Then there is a unique subbundle $\cG\subseteq T_Y$ such that
$df|_{\cE}$ identifies $\cE$ with $f^*\cG$ as a subbundle of
$f^*T_Y$. If $\cE$ is integrable, then so is $\cG$.
\end{lemma}

\begin{proof}
Let $r:=\rk\cE$ and $\cH:=df(\cE)\subseteq f^*T_Y$. As in the
proof of Proposition~\ref{Prop:almost-nef-descent}, let $\pi:\Gr_Y(r,T_Y)\rightarrow Y$ be the relative Grassmannian with universal subbundle $\cS\hookrightarrow\pi^*T_Y$. The subbundle
$\cH\hookrightarrow f^*T_Y$ defines a classifying morphism
\[
 \varphi:X\rightarrow\Gr_Y(r,T_Y),
 \qquad \varphi^*\cS\simeq\cH.
\]

The Bott partial connection (see \cite[Section~2.2]{Druel17}) on
$T_X/T_{X/Y}\simeq f^*T_Y$ along $T_{X/Y}$ is given by
\[
 \nabla_u(df(v))=df([u,v]),
 \qquad u\in T_{X/Y},\quad v\in T_X.
\]
For $v\in \cE\subseteq T_X$, $\nabla_u(df(v))=df([u,v])\in df(\cE)$ by
\eqref{eq:wf-mixed-bracket}, whereby the connection $\nabla_u$ preserves $\cH$. Thus $\nabla$ restricts to a relative holomorphic connection
\[
 \nabla^{\cH}:\cH\rightarrow\Omega^1_{X/Y}\otimes\cH.
\]
Moreover, for local vertical vector fields $u,w\in T_{X/Y}$ and a local section $v$ of $T_X$, the Jacobi identity gives
\[
 \begin{aligned}
 \bigl(\nabla_u\nabla_w-\nabla_w\nabla_u-\nabla_{[u,w]}\bigr)(df(v)) =df\bigl([u,[w,v]]-[w,[u,v]]-[[u,w],v]\bigr)=0,
 \end{aligned}
\]
which implies that the curvature of the restricted relative connection $\nabla^{\cH}$ is zero and hence $\nabla^{\cH}$ is flat. Now fix a fiber $F=X_y$, with ideal sheaf $\cI_F$. This ideal is locally
generated by pullbacks of functions on $Y$ vanishing at $y$, which
are annihilated by vertical vector fields, i.e, sections of $T_{X/Y}$. The Leibniz rule therefore yields
\[
 \nabla^{\cH}(\cI_F\cH)\subseteq\cI_F\bigl(\Omega^1_{X/Y}\otimes\cH\bigr).
\]
Consequently, $\nabla^{\cH}$ induces a connection on $\cH|_F$.
Since $f$ is smooth and $\Omega^1_{X/Y}|_F\simeq\Omega^1_F$, we obtain a flat holomorphic connection
\[
 \nabla^F:\cH|_F\rightarrow\Omega^1_F\otimes\cH|_F,
\]
since its curvature vanishes by the preceding calculation. By \cite[Lemma~2.6 and its proof]{Druel17}, we obtain $c_1(\cH|_F)=0$ in $H^2(F,\CC)$. In particular,
\[
 (\varphi|_F)^*\OO_{\Gr(r,T_{Y,y})}(1)
 \simeq\det(\cH|_F)^\vee
 \equiv_{\mathrm{num}}0.
\]
Since the Pl\"ucker line bundle is ample, $\varphi|_F$ is constant,
as in the proof of Proposition~\ref{Prop:almost-nef-descent}. The faithfully flat descent argument in that proof now yields a
unique section
$\overline\varphi:Y\rightarrow\Gr_Y(r,T_Y)$ satisfying
$\varphi=\overline\varphi\circ f$. Thus
$\cG:=\overline\varphi^*\cS\subseteq T_Y$ is the required unique
subbundle, with
$f^*\cG=\cH\simeq\cE$ under $df|_{\cE}$.
If $\cE$ is integrable, the projectable-vector-field argument in
the proof of Proposition~\ref{Proposition:4}(b) shows that $\cG$
is integrable as well.
\end{proof}

We record the fiber-classification argument that will be used in the present section. The proof is an amalgamation of statements already proven in Sections~\ref{sec:CharacterizeP^nQ_n} and~\ref{sec:CharacterizeProduct1}.

\begin{lemma}\label{lem:wf-fibers}
Let $f:X\rightarrow Y$ be a smooth, surjective morphism
between smooth projective complex varieties, with rationally
connected fibers of dimension $d>0$. Let $\cV$ be a nef vector bundle
of rank $s$, ample on every fiber of $f$, and assume the existence of an injection
\[
 \bS_\lambda\cV\hookrightarrow\bS_\lambda T_{X/Y}
\]
of coherent sheaves, for a nonzero partition $\lambda$ of length at most $\min\{s,d\}$. Then all
fibers are isomorphic to $\PP^d$, or all are isomorphic to $Q_d$ with
$d\geq2$. Furthermore letting $\cL:=\det\cV$, there is a fixed integer $\ell\geq s$
such that $\cL|_F\simeq\OO_F(\ell)$ on every fiber $F$.
\end{lemma}

\begin{proof}
 The proof of all fibers of $f:X\rightarrow Y$ being isomorphic to $\PP^d$, or all fibers being isomorphic to $Q_d$ with $d\geq2$ is the same as that of Lemma~\ref{lem:fib-characterization} with the numerically flat quotient $\cG$ equal to zero. By applying Theorem~\ref{thm:generalizedKovacs} to a general fiber $F$ of $f:X\rightarrow Y$, we obtain
$\cL|_F\simeq\OO_F(\ell)$, with $\ell\geq s$. The relative Picard argument in the proof of Proposition~\ref{Proposition:5}, applied to
$\cL^{\otimes\iota}\otimes\omega_{X/Y}^{\otimes\ell}$, where
$\iota=d+1$ for projective spaces and $\iota=d$ for quadrics, extends
this identity to every fiber. Torsion freeness of $\Pic(F)$, as in
Lemma~\ref{lem:NF-RC}(iii), then yields the assertion about $\cL$.
\end{proof}

The morphism in Lemma~\ref{lem:wf-fibers} is a Brauer-Severi scheme
or a smooth quadric bundle, as in the proof of
Proposition~\ref{Proposition:4}. Thus
Proposition~\ref{prop:relative-ample-base} applies once its tensor
morphism and numerical inequality are verified. We make this
verification in the inductive proof below.

\subsection{Proof of Theorem~\ref{thm:productcharacterization2}}\label{subsec:wf-proof}
We prove Theorem~\ref{thm:Main3} by induction on $r$, similar to Section~\ref{subsec:proofproduct}.

\subsubsection{Base case $r=1$}
In this case, $\cV_1$ is ample and $\cE_1$ is a regular foliation.
Composing $\theta_1:\bS_{\lambda_1}\cV_1\rightarrow\bS_{\lambda_1}\cE_1$ with the morphism induced by the subbundle inclusion $\cE_1\hookrightarrow T_X$ gives an injection $\bS_{\lambda_1}\cV_1\hookrightarrow\bS_{\lambda_1}T_X$. Thus
Theorem~\ref{thm:generalizedKovacs}, with unipotent quotient $\cU=0$,
gives the desired classification of $X$ and $\cV_1$ in Theorem~\ref{thm:productcharacterization2}. It remains to prove
$\cE_1=T_X$. As in Section~\ref{subsubsec:product-base-case}, we
consider the following two cases.
\begin{enumerate}
\item If $X\not\simeq Q_2$, then $X\simeq\PP^n$ or
$X\simeq Q_n$ with $n\geq3$, so $\rho(X)=1$. By
\cite[Theorem~1.3]{Druel17}, there is a smooth projective morphism
$f:X\rightarrow Y$ with connected fibers such that
$\cE_1=T_{X/Y}$. Its fibers have positive dimension since
$n_1\geq1$. The Picard-number-one argument in
Section~\ref{subsubsec:product-base-case}(1) therefore shows that
$Y$ is a point. Hence $\cE_1=T_X$.
\item If $X\simeq Q_2$, the rank-one case is excluded by the same
calculation as in Section~\ref{subsubsec:product-base-case}(2): since that
calculation only uses ampleness of $\cV_1$, the Schur injection, and
the subbundle inclusion $\cE_1\hookrightarrow T_{Q_2}$. Thus
$\rk\cE_1=2$, and the locally free cokernel of this inclusion has
rank zero. Hence $\cE_1=T_{Q_2}$.
\end{enumerate}
Consequently, $n=n_1$, and
Lemma~\ref{lem:Schur-rank-comparison} gives $r_1\leq n_1$.
Taking $Z_1=X$, $\cT_1=\cV_1$, and $\overline\theta_1=\theta_1$
proves the base case.

\subsubsection{Inductive step}
The following proposition provides the inductive reduction for Theorem~\ref{thm:productcharacterization2}, analogous to Proposition~\ref{Proposition:4}.

\begin{proposition}\label{prop:wf-inductive-reduction}
Let $X$ be a smooth weak Fano complex variety of dimension $n>0$, and
let $r\geq1$. For $1\leq i\leq r$, let $\cE_i\subseteq T_X$ be a
vector subbundle of rank $n_i\geq1$, let $\cV_i$ be a nef vector
bundle of rank $r_i\geq1$, and let $\lambda_i$ be a nonzero partition
of length $\ell(\lambda_i)\leq\min\{r_i,n_i\}$. Assume that:
\begin{enumerate}
\item the natural morphism
$\cF:=\bigoplus_{i=1}^r\cE_i\rightarrow T_X$ is a subbundle inclusion with integrable image,
\item for every $1\leq i\leq r$, there is an injective morphism
$\theta_i:\bS_{\lambda_i}\cV_i\hookrightarrow
\bS_{\lambda_i}\cE_i$ of coherent $\OO_X$-modules and,
\item the vector bundle $\cA:=\bigotimes_{i=1}^r\cV_i$ is ample on $X$.
\end{enumerate}
Then, after reordering the indices, there exists a Brauer-Severi scheme or smooth quadric bundle $f:X\rightarrow Y$ with connected fibers of dimension $n_1$, where $Y$ is a smooth weak Fano variety and $\cE_1=T_{X/Y}$. The following properties are satisfied:
\begin{enumerate}
\item[(a)] there exist nef vector bundles $\cW_2,\dots,\cW_r$ on
$Y$ such that $\cV_j\simeq f^*\cW_j$ for all $j\geq2$,
\item[(b)] there exist regular foliations $\cG_j\subseteq T_Y$
such that $df|_{\cE_j}$ identifies $\cE_j$ with $f^*\cG_j$ for
all $j\geq2$. Furthermore, the  injective morphisms $\theta_j:\bS_{\lambda_j}\cV_j\hookrightarrow
\bS_{\lambda_j}\cE_j$ of coherent sheaves descend to
injective morphisms 
\begin{equation}\label{eq:wf-main-descent}
 \overline\theta_j:
 \bS_{\lambda_j}\cW_j\hookrightarrow\bS_{\lambda_j}\cG_j,
 \qquad f^*\overline\theta_j=\theta_j\quad(j\geq2),
\end{equation}
\item[(c)] the natural morphism
$\bigoplus_{j=2}^r\cG_j\rightarrow T_Y$ is a subbundle inclusion
with integrable image,
\item[(d)] $\cU:=\cW_2\otimes\cdots\otimes\cW_r$ is ample on $Y$.
\end{enumerate}
\end{proposition}

\begin{proof}
By Lemma~\ref{lem:wf-partial-integrability}, all partial sums of
the $\cE_i$'s are regular foliations. Applying
\cite[Theorem~1.3]{Druel17} to $\cE_1$ produces a smooth projective
morphism $f:X\rightarrow Y$ with connected fibers and
$\cE_1=T_{X/Y}$. By Lemma~\ref{lem:wf-targetfiber}, $Y$ and all the fibers of $f$ are weak Fano and rationally connected. As in the proof of Lemma~\ref{lem:fib-restriction}(1), the differential
induces subbundle inclusions
\[
 df|_{\cE_j}:\cE_j\hookrightarrow f^*T_Y\quad(j\geq2),
 \qquad
 \bigoplus_{j=2}^r\cE_j\hookrightarrow f^*T_Y.
\]
Composing with the Schur morphisms and applying
Lemma~\ref{lem:fiber-trivial}, exactly as in the proof of
Lemma~\ref{lem:fib-restriction}(2), yields
$\cV_j\simeq f^*\cW_j$ for vector bundles $\cW_j$ on $Y$.
Nefness of $\cW_j$ follows from \cite[Proposition~1.10]{DPS94}. This proves~(a).

By Lemma~\ref{lem:wf-partial-integrability}, $\cE_1\oplus\cE_j$ is integrable, i.e., $[T_{X/Y},\cE_j]\subseteq T_{X/Y}\oplus\cE_j$. Lemma~\ref{lem:wf-transverse-descent} therefore gives regular foliations $\cG_j\subseteq T_Y$ such that $\cE_j\simeq f^*\cG_j$ under $df$. Since $f_*\OO_X=\OO_Y$, the projection formula descends the injective morphisms $\theta_j:\bS_{\lambda_j}\cV_j\hookrightarrow\bS_{\lambda_j}\cE_j$ to the morphisms  \eqref{eq:wf-main-descent}, which are injective by faithful flatness. This is the same descent step as in the
proof of Proposition~\ref{Proposition:4}(b), and proves~(b).

The subbundle inclusion in~(c) follows by faithful flatness.
Applying Lemma~\ref{lem:wf-transverse-descent} to the integrable
partial sum $\bigoplus_{j=2}^r\cE_j$ also proves integrability
of the descended image of $\bigoplus_{j=2}^r\cG_j$ in $T_Y$. Thus~(c) holds.

We now prove (d). Restricting the ample bundle $\cA:=\bigotimes_{i=1}^r\cV_i$ to any fiber $F$ of $f:X\rightarrow Y$, we obtain $\cA|_F\simeq(\cV_1|_F)^{\oplus M}$, where $M:=\prod_{j=2}^r r_j$, since for all $j\geq 2$, $\cV_j|_F\simeq\OO_F^{\oplus r_j}$ by part (a) and therefore, $\cV_1|_F$ is ample. Since $\cE_1=T_{X/Y}$, by assumption, we have the injective morphism
\[
 \theta_1:\bS_{\lambda_1}\cV_1
 \hookrightarrow\bS_{\lambda_1}\cE_1=\bS_{\lambda_1}T_{X/Y}.
\]
By Lemma~\ref{lem:wf-fibers}, we conclude that $f$ is a Brauer-Severi
scheme or smooth quadric bundle, as in the proof of
Proposition~\ref{Proposition:4}. Furthermore, if $\cL:=\det\cV_1$ then by Lemma~\ref{lem:wf-fibers}, there is a fixed integer $\ell\geq r_1$ such that $\cL|_F\simeq\OO_F(\ell)$ for every fiber $F$. Let $a=q_{\lambda_1}(r_1)=|\lambda_1|d_{\lambda_1}(r_1)/r_1$ and $R=ar_1$. Applying Proposition~\ref{prop:Schur-det-to-tensor} to $\theta_1$ yields a nonzero morphism $\cL^{\otimes a}\rightarrow T_{X/Y}^{\otimes R}$.
Since $a\ell\geq ar_1=R$ and
$\cA=\cV_1\otimes f^*\cU$ is ample,
Proposition~\ref{prop:relative-ample-base} applies with
$\cV=\cV_1$, $q=R$, and $\cP=\OO_Y$. It follows that
$\cU$ is ample, proving~(d).
\end{proof}

We can now complete the inductive step for the proof of
Theorem~\ref{thm:productcharacterization2}.

\begin{proof}[Proof of Theorem~\ref{thm:productcharacterization2}]
Assume $r\geq2$ and that Theorem~\ref{thm:productcharacterization2} has been proved for $r-1$ subbundles of ranks $n_2,n_3,\dots,n_r\geq1$. Applying Proposition~\ref{prop:wf-inductive-reduction}, there exists a smooth Brauer-Severi scheme or a smooth quadric bundle $f:X\rightarrow Y$, where $Y$ is a smooth weak Fano variety and $\cE_1=T_{X/Y}$. Furthermore, for all $2\leq j\leq r$, we have:
\begin{enumerate}
    \item nef vector bundles $\cW_2,\dots,\cW_r$ on $Y$ such that $\cV_j\simeq f^*\cW_j$ and $\cW_2\otimes\cdots\otimes\cW_r$ is ample on $Y$,
    \item regular foliations $\cG_j\subseteq T_Y$ of rank $n_j$ such that $df|_{\cE_j}$ identifies $\cE_j$ with $f^*\cG_j$, with a subbundle inclusion $\cG_2\oplus\cdots\oplus\cG_r\hookrightarrow T_Y$ whose image is integrable, and
    \item injective morphisms $\overline\theta_j:\bS_{\lambda_j}\cW_j
    \hookrightarrow\bS_{\lambda_j}\cG_j$ of coherent sheaves such that
    $f^*\overline\theta_j=\theta_j$.
\end{enumerate}
Hence the induction hypothesis applies to the tuple $\bigl(Y;\cG_2,\dots,\cG_r;\cW_2,\dots,\cW_r;\lambda_2,\dots,\lambda_r\bigr)$, whereby we conclude that $Y\simeq Z_2\times\cdots\times Z_r$ and $T_Y=\bigoplus_{j=2}^r\cG_j$, with $\cG_j\simeq p_j^*T_{Z_j}$ for all $2\leq j\leq r$, where $p_j:Y\simeq Z_2\times\cdots\times Z_r\rightarrow Z_j$ is
the natural projection map. Furthermore, for all $2\leq j\leq r$, we have $Z_j\simeq\PP^{n_j}$ or a smooth hyperquadric $Q_{n_j}\subseteq\PP^{n_j+1}$ of dimension at least $2$, and there are ample bundles $\cT_j$ on $Z_j$ satisfying $\cW_j\simeq p_j^*\cT_j$ and admitting injective morphisms $\bS_{\lambda_j}\cT_j\hookrightarrow\bS_{\lambda_j}T_{Z_j}$ of coherent sheaves. These morphisms pull back to $\overline\theta_j$, and the bundles $\cT_j$ and partitions $\lambda_j$ are classified as in Theorem~\ref{thm:Main2}.

By Proposition~\ref{prop:wf-inductive-reduction}, $\cE_1=T_{X/Y}$, whereby the fibers of $f:X\rightarrow Y$ have dimension equal to $\rk\cE_1=n_1$. Thus, $\dim X=n_1+\dim Y=n_1+n_2+\cdots+n_r$, which is also the rank of $T_X$. By comparing ranks, it follows that the subbundle inclusion
$\bigoplus_{i=1}^r\cE_i\hookrightarrow T_X$ in the hypothesis of Theorem~\ref{thm:Main3} is an isomorphism. In other words, $f:X\rightarrow Y\simeq Z_2\times\cdots\times Z_r$ admits a connection:
\begin{equation}\label{eq:wf-whole-tangent}
 T_X=T_{X/Y}\oplus\bigoplus_{j=2}^r\cE_j =\bigoplus_{i=1}^r\cE_i,
\end{equation}
where $df$ identifies the horizontal summand $\bigoplus_{j=2}^r\cE_j$ with $f^*T_Y$. This summand is integrable by Lemma~\ref{lem:wf-partial-integrability}. Since $Y$ is rationally connected, \cite[Lemma~3.18]{Horing07} gives an isomorphism $X\simeq F\times Y$, where $F$ is a fiber of
$f:X\rightarrow Y$ and $f$ is the projection onto $Y$. The isomorphism respects the relative tangent bundle and the horizontal summand, as in \cite[Theorem~3.17]{Horing07}. By Proposition~\ref{prop:wf-inductive-reduction}, $F\simeq\PP^{n_1}$ or $Q_{n_1}$, whereby, setting $Z_1:=F$,
we obtain $X\simeq Z_1\times\cdots\times Z_r$, where each $Z_i$ is $\PP^{n_i}$ or a smooth hyperquadric $Q_{n_i}\subseteq\PP^{n_i+1}$ of dimension at least $2$. Let $\pi_i:X\simeq Z_1\times\cdots\times Z_r\rightarrow Z_i$ be the natural projection onto the $i^{th}$ factor. Then
$\pi_j=p_j\circ f$ for all $j\geq2$, and thus $\cE_j\simeq f^*\cG_j\simeq\pi_j^*T_{Z_j}$ for $j\geq2$. For $j=1$, we have $\cE_1=T_{X/Y}\simeq\pi_1^*T_{Z_1}$. This proves \eqref{eq:wf-main-product2}.
Similarly, for $j\geq2$, $\cV_j\simeq f^*\cW_j\simeq\pi_j^*\cT_j$, where $\cT_j$ are ample vector bundles on $Z_j$, classified as in Theorem~\ref{thm:Main2} by the induction hypothesis. The corresponding assertions about $\theta_j$ follow from the induction hypothesis and \eqref{eq:wf-main-descent}. Therefore, what remains is to prove that
$\cV_1\simeq\pi_1^*\cT_1$ for an ample vector bundle $\cT_1$ on $Z_1$, to descend $\theta_1$, and to classify $\cT_1$.

By hypothesis, we have the following injection of coherent sheaves, where we use $\cE_1\simeq\pi_1^*T_{Z_1}$ and Proposition~\ref{prop:Schur-functorial}(iii):
\begin{equation}\label{eq:wf-V1-injection}
 \theta_1:\bS_{\lambda_1}\cV_1
 \hookrightarrow\bS_{\lambda_1}(\pi_1^*T_{Z_1})
 \simeq\pi_1^*(\bS_{\lambda_1}T_{Z_1}).
\end{equation}
Since $\cV_1$ is nef and the fibers of $\pi_1:X\rightarrow Z_1$ are rationally connected, we can apply Lemma~\ref{lem:fiber-trivial} to obtain the unique vector bundle $\cT_1:=\pi_{1,*}\cV_1$ on $Z_1$ such that $\pi_1^*\cT_1\simeq\cV_1$ on $X$. For any $y\in Y$, let $\iota_y:Z_1\rightarrow Z_1\times Y$ be the section of $\pi_1$ defined by $z\mapsto(z,y)$. Since $\pi_1\circ\iota_y=\id_{Z_1}$, we obtain $\cT_1\simeq\iota_y^*\pi_1^*\cT_1\simeq\iota_y^*\cV_1$. Substituting $\pi_1^*\cT_1\simeq\cV_1$ in \eqref{eq:wf-V1-injection}, we obtain an injection $\theta_1:\pi_1^*(\bS_{\lambda_1}\cT_1)\hookrightarrow\pi_1^*(\bS_{\lambda_1}T_{Z_1})$. Since $\pi_{1,*}\OO_X\simeq\OO_{Z_1}$, the projection formula gives a unique morphism on $Z_1$ whose pullback is $\theta_1$, whose injectivity follows by faithful flatness of $\pi_1$.
Thus we have an injection of coherent sheaves
\begin{equation}\label{eq:wf-T1-injection}
 \overline\theta_1:\bS_{\lambda_1}\cT_1
 \hookrightarrow\bS_{\lambda_1}T_{Z_1},
 \qquad
 \pi_1^*\overline\theta_1=\theta_1.
\end{equation}

We now show that $\cT_1$ is an ample bundle on $Z_1$. Restricting the ample bundle $\cA=\cV_1\otimes\cdots\otimes\cV_r$ on $X\simeq Z_1\times Y$ to the closed subvariety $Z_1\times\{y\}$ for any $y\in Y$, we see that
\[
 \iota_y^*\cA\simeq\cT_1\otimes_{\CC}\left(\bigotimes_{j=2}^r\cT_{j,p_j(y)}\right)
 \simeq\cT_1^{\oplus M}
\]
is an ample bundle on $Z_1$, where $M:=\prod_{j=2}^r r_j$ and $\cT_{j,p_j(y)}$ is the fiber of the bundle $\cT_j$ at $p_j(y)\in Z_j$.
This is because $\pi_j\circ\iota_y:Z_1\rightarrow Z_j$ is
the constant map to $p_j(y)$ for all $j\geq2$. Consequently, $\cT_1$ is ample, being a quotient of $\cT_1^{\oplus M}$. Applying Theorem~\ref{thm:generalizedKovacs} to \eqref{eq:wf-T1-injection} gives the desired classification of $\cT_1$ and the partition $\lambda_1$, while
Lemma~\ref{lem:Schur-rank-comparison} gives $r_1\leq n_1$. The corresponding inequalities for $j\geq2$ follow from the induction hypothesis. Finally, every factor $Z_i$ is Fano, so $X$ is Fano. This completes the proof of
Theorem~\ref{thm:productcharacterization2} by induction on the number $r$
of subbundles.
\end{proof}

\subsection{Toric weak Fano varieties}
\label{subsec:toric-foliations}
The following result shows that for toric weak Fano varieties, the existence of enough ``independent" regular foliations on $X$ is sufficient to force it to be a product of projective spaces, without any positivity conditions. 

\begin{theorem}\label{thm:toric-product}
Let $X$ be a smooth complex projective toric weak Fano variety of dimension $n>0$, and let $\cE_1,\ldots,\cE_r\subseteq T_X$
be regular foliations of ranks $n_i\geq1$ such that the natural
morphism $\bigoplus_{i=1}^r\cE_i\to T_X$ is a subbundle inclusion.
Then $r\leq\rho(X)$. If $r=\rho(X)$, there is an isomorphism
\[
 X\simeq\prod_{i=1}^r\PP^{n_i},\qquad
 \cE_i=\pr_i^*T_{\PP^{n_i}}\subseteq T_X,
\]
where $\pr_i$ is the $i^{th}$ projection under this isomorphism.
Conversely, every product of positive-dimensional projective
spaces admits such a collection of regular foliations of length $r=\rho(X)$.
\end{theorem}

\begin{proof}
By \cite[Theorem~1.3]{Druel17}, each $\cE_i$ produces a smooth projective morphism $q_i:X\to Y_i$ with connected fibers, where $Y_i$ is smooth and projective and $\cE_i=T_{X/Y_i}$. We obtain $H^1(X,\OO_X)=0$ by Lemmas~\ref{lem:wf-targetfiber}
and~\ref{lem:NF-RC}\textup{(iii)}. Thus, Blanchard's theorem \cite{Blanchard} in the form recalled in \cite[Section \emph{Toric varieties}]{OW02}, descends the torus $\mathbb{T}$ action on $X$ to $Y_i$ such that $q_i:X\rightarrow Y_i$ is equivariant with respect to these actions. Taking differentials of the $\mathbb{T}$-equivariant map $q_i$ shows that
$\cE_i=T_{X/Y_i}=\ker(dq_i)$ is $\mathbb{T}$-invariant. By \cite[Proposition~3.1]{ChangChen25}, for each $1\leq i\leq r$ there is a linear subspace $W_i\subseteq\mathfrak t:=\operatorname{Lie}(\mathbb T)$ such that, under the translation trivialization,
\[
 \cE_i|_{\mathbb T}
 =\OO_{\mathbb T}\otimes_{\CC}W_i\subseteq
 T_{\mathbb T}\simeq \OO_{\mathbb T}\otimes_{\CC}\mathfrak t.
\]
Consequently, for every $I\subseteq\{1,\ldots,r\}$, $\cE_I|_{\mathbb T}=\OO_{\mathbb T}\otimes_{\CC} \bigl(\bigoplus_{i\in I}W_i\bigr)$. Applying the converse construction of proof of \cite[Proposition~3.1]{ChangChen25} to $\bigoplus_{i\in I}W_i\subseteq \OO_{\mathbb T}\otimes_{\CC}\mathfrak{t}$ (obtained by restricting the subbundle inclusion $\cE_I\subseteq T_X$ to the dense torus $\mathbb{T}\subseteq X$), we see that every partial sum $\cE_I$ is a $\mathbb{T}$-invariant regular foliation on $X$.

We first prove the inequality $r\leq \rho(X)$. Fix a fiber $F_i$ of $q_i$. For $j\ne i$, integrability of
$\cE_i\oplus\cE_j$ gives
$[T_{X/Y_i},\cE_j]\subseteq T_{X/Y_i}\oplus\cE_j$.
Lemma~\ref{lem:wf-transverse-descent} therefore identifies
$\cE_j$ with a bundle pulled back from $Y_i$. Consequently,
\[
 \cE_i|_{F_i}=T_{F_i},\qquad
 \cE_j|_{F_i}\simeq\OO_{F_i}^{\oplus n_j}\quad(j\ne i).
\]
By Lemma~\ref{lem:wf-targetfiber}, $F_i$ is weak Fano.
Choose a general ample complete-intersection curve
$C_i\subset F_i$ if $n_i>1$, and $C_i=F_i$ when $n_i=1$.
Since $-K_{F_i}$ is nef and big, it is numerically nonzero, whereby
Lemma~\ref{lem:nefzero} yields
$a_i:=(-K_{F_i})\cdot C_i>0$. Setting $\cL_j:=\det\cE_j$,
we obtain
\begin{equation}\label{eq:toric-diagonal}
 \cL_j\cdot C_i=\delta_{ij}a_i,\qquad a_i>0.
\end{equation}
Thus $[\cL_1],\ldots,[\cL_r]$ are linearly independent in
$N^1(X)_{\RR}$, proving $r\leq\rho(X)$. For the rest of the proof, we focus on the equality $r=\rho(X)$ case.

Suppose that $r=\rho(X)$, so the preceding classes $[\cL_1],\ldots,[\cL_r]$ form a basis of $N^1(X)_{\RR}$. Applying \cite[Theorem~1.3]{Druel17} to the regular foliation $\cF:=\bigoplus_i\cE_i$, we obtain a smooth morphism $q:X\rightarrow B$ with connected fibers and $\cF=T_{X/B}$. Let $\iota_i:F_i\hookrightarrow X$ be the inclusion of any fiber of the preceding maps $q_i:X\rightarrow Y_i$.
Since $d\iota_i(T_{F_i})\subseteq T_{X/B}|_{F_i}=\ker (dq)|_{F_i}$, we have
\[
 d(q|_{F_i})=(dq)|_{F_i}\circ d\iota_i=0.
\]
Hence $q|_{F_i}$ is locally constant, and therefore
constant because $F_i$ is connected. In particular, $q$ contracts
every $C_i$. If $\dim B>0$ and $H$ is ample on $B$,
surjectivity of $q$ gives $[q^*H]\ne0$, whereas
$q^*H\cdot C_i=0$ for all $i$. Since $[\cL_1],\ldots,[\cL_{\rho(X)}]$ form a basis of $N^1(X)_{\RR}$, writing $[q^\ast H]=\sum_{i}b_i[\cL_i]$, we see $0=q^\ast H\cdot C_i=b_ia_i$ by \eqref{eq:toric-diagonal}, which implies $b_i=0$ for all $1\leq i\leq \rho(X)$. This contradicts $[q^\ast H]\neq 0$. Therefore $B$ must be a point and $T_X=\bigoplus_{i=1}^{\rho(X)}\cE_i$. We can now prove that $X$ is a product of projective spaces.

For this, we use the connection argument from the proof of
Theorem~\ref{thm:productcharacterization2}. For $r=\rho(X)>1$,
Lemma~\ref{lem:wf-transverse-descent} applied along $q_1:X\rightarrow Y_1$ with $\cE_1=T_{X/Y_1}$, descends the remaining summands $\cE_j$'s ($2\leq j\leq r$) to a splitting
$T_{Y_1}=\bigoplus_{j=2}^r\cG_j$ such that all the corresponding partial sums $\cG_J:=\bigoplus_{j\in J}\cG_j$ on $Y_1$ are integrable for $J\subseteq \{2,\dots, n\}$. By
Lemma~\ref{lem:wf-targetfiber}, $Y_1$ is weak Fano and
rationally connected. Applying \cite[ Lemma~3.18]{Horing07} to the connection
\[
 T_X=T_{X/Y_1}\oplus\bigoplus_{j=2}^r\cE_j
\]
yields $X\simeq F_1\times Y_1$, such that $\cE_j\simeq q_1^\ast \cG_j$, for all $j\geq2$. Iterating on $Y_1$ yields
\[
 X\simeq\prod_{i=1}^rZ_i,\qquad
 \cE_i=\pr_i^*T_{Z_i}\subseteq T_X,
\]
where the $Z_i$ are smooth projective varieties of dimensions
$n_i\geq1$. Pullbacks of divisor classes from different factors are
numerically independent, we have
\[
 r=\rho(X)\geq\sum_{i=1}^r\rho(Z_i)\geq r,
\]
whereby $\rho(Z_i)=1$ for every $1\leq i\leq r$. Each projection $\pr_i: X\rightarrow Z_i$
is surjective from a complete toric variety and thus $Z_i\simeq\PP^{n_i}$ for all $1\leq i\leq r$ by \cite[Theorem~1]{OW02}.
The converse follows from the factor tangent bundles and the
Picard number of a product of projective spaces.

\end{proof}

\section{Special characterizations}\label{sec:special}

In this section we prove the product analogues of \cite{KO73, Mori79, Ye94} as special instances of the general characterization Theorems~\ref{thm:Main2} and \ref{thm:Main3}.

\begin{theorem}\label{specialcase1proof}
 Let $X$ be a smooth projective complex variety of dimension $n>0$ and let $r>0$ be an integer. Furthermore let $n_1,\dots,n_r$ be positive integers such that $T_X$ admits a holomorphic splitting $T_X=\bigoplus_{i=1}^r\cE_i$, where each $\cE_i$ is nef vector bundle of  $\rk\cE_i=n_i$ and $\bigotimes_{i=1}^r\cE_i$ is ample. Then, possibly after reordering, there is an isomorphism 
 \[X\simeq \prod_{i=1}^{r}\PP^{n_i} \text{ such that } \cE_i\simeq \pr_i^*T_{\PP^{n_i}},\] 
 where $\pr_i:X\rightarrow \PP^{n_i}$ is the natural projection. 
\end{theorem}

\begin{proof}
Let $\cA:=\bigotimes_{i=1}^r\cE_i$ and $R:=\prod_{i=1}^r n_i=\rk\cA$. By Lemma~\ref{lem:schur-positivity}, the line bundle $\det\cA$ is ample and each $\det\cE_i$ is nef. Since
\[
 \det\cA\simeq
 \bigotimes_{i=1}^r(\det\cE_i)^{\otimes R/n_i},
 \qquad
 \OO_X(-K_X)\simeq\bigotimes_{i=1}^r\det\cE_i,
\]
we obtain
\[
 \OO_X(-K_X)^{\otimes R}\simeq
 \det\cA\otimes
 \bigotimes_{i=1}^r
 (\det\cE_i)^{\otimes(R-R/n_i)}
\]
is ample and thus $X$ is Fano. For every $1\leq i\leq r$, let
\[
 \cV_i=\cE_i,\qquad r_i=n_i,\qquad
 \lambda_i=(1),\qquad \theta_i=\id_{\cE_i}.
\]
Since $\bS_{(1)}\cE_i=\cE_i$ and
$\cF:=\bigoplus_{i=1}^r\cE_i=T_X$ is integrable, all the
hypotheses of Theorem~\ref{thm:Main3} are satisfied. Consequently, possibly after reordering, we obtain the desired isomorphism
\[
 X\simeq\prod_{i=1}^r Z_i,
 \qquad
 \cE_i\simeq\pr_i^*T_{Z_i},
\]
where each $Z_i$ is $\PP^{n_i}$ or a smooth hyperquadric
$Q_{n_i}$ with $n_i\geq2$.
The quadric alternative, by
Theorem~\ref{thm:Main2}\textup{(d)}, requires
$\lambda_i=(k_i^{r_i})$ for some $k_i\geq1$, and hence
$\ell(\lambda_i)=r_i=n_i\geq2$. This contradicts
$\lambda_i=(1)$. Therefore $Z_i\simeq\PP^{n_i}$ for every $i$,
which proves the assertion.
\end{proof}

We now prove the product analogue of \cite{KO73}.

\begin{theorem}\label{specialcase2proof}
Let $X$ be a smooth projective complex variety of dimension $n>0$ and let $r>0$ be an integer. Let $n_1,\dots,n_r$ be positive integers such that $T_X$ admits a holomorphic splitting $T_X=\bigoplus_{i=1}^r\cE_i$ with $\rk\cE_i=n_i$. Furthermore, assume that for each $1\leq i\leq r$, there exists a nef line bundle $\cL_i$ such that $\bigotimes_{i=1}^r\cL_i$ is ample. Then the following statements hold.
\begin{enumerate}[(i)]
\item If $\cM_i:=\det\cE_i\otimes\cL_i^{-n_i}$ is nef for every $1\leq i\leq r$, then, possibly after reordering, there is an isomorphism
\begin{equation}\label{Eq:productspecialcaseproof}
 X\simeq\prod_{i=1}^r Z_i
 \text{ such that }
 \cE_i\simeq\pr_i^*T_{Z_i}
 \text{ and }
 \cL_i\simeq\pr_i^*H_i,
\end{equation}
where $\pr_i:X\rightarrow Z_i$ is the natural projection and for all $1\leq i\leq r$:
\[
 (Z_i,H_i, \cM_i)\simeq(\PP^{n_i},\OO_{\PP^{n_i}}(1), \cL_i)
 \quad\text{or}\quad
 (Z_i,H_i, \cM_i)\simeq(Q_{n_i},\OO_{Q_{n_i}}(1), \OO_X).
\]
\item If $\det\cE_i\otimes\cL_i^{-\otimes(n_i+1)}$ is nef for every $1\leq i\leq r$, then \eqref{Eq:productspecialcaseproof} holds such that for all $1\leq i\leq r$:  \[(Z_i,H_i, \det\cE_i)\simeq(\PP^{n_i},\OO_{\PP^{n_i}}(1), \cL_i^{\otimes (n_i+1)}).\]
\item If $\det\cE_i\equiv\cL_i^{n_i}$ for every $1\leq i\leq r$, then \eqref{Eq:productspecialcaseproof} holds such that for all $1\leq i\leq r$:  \[(Z_i,H_i, \det\cE_i)\simeq(Q_{n_i},\OO_{Q_{n_i}}(1), \cL_i^{\otimes n_i}).\] Here $\equiv$ denotes numerical equivalence of line bundles.
\end{enumerate}

\end{theorem}

\begin{proof}
First we prove (i). The splitting $T_X=\bigoplus_{i=1}^r\cE_i$ and definition of $\cM_i$ yields
\[
 \OO_X(-K_X)\simeq
 \bigotimes_{i=1}^r\det\cE_i
 \simeq
 \left(\bigotimes_{i=1}^r\cL_i\right)\otimes
 \bigotimes_{i=1}^r
 \left(\cL_i^{\otimes(n_i-1)}\otimes\cM_i\right).
\]
Since $\bigotimes_{i=1}^r\cL_i$ is ample and $\cL_i$ and $\cM_i$ are nef, it follows that $X$ is Fano. Thus the base point free theorem \cite[Theorem~3.3]{KM98} implies that each $\cM_i$ is semiample since it is nef and $\cM_i\otimes\OO_X(-K_X)$ is ample. Choose an integer $k>0$ such that $\cM_i^{\otimes k}$ is globally generated for all $1\leq i\le r$, and choose $0\neq \sigma_i\in H^0(X, \cM_i^{\otimes k})$. For $1\leq i\leq r$, set $\cV_i=\cL_i^{\oplus n_i}$ and $\lambda_i=(k^{n_i})$. By Example~\ref{Ex:1}(iii) we have
\[
 \bS_{\lambda_i}\cV_i\simeq\cL_i^{kn_i},
 \qquad
 \bS_{\lambda_i}\cE_i\simeq(\det\cE_i)^k.
\]
The chosen non-zero sections $\sigma_i\in H^0(X, \cM_i^{\otimes k})= H^0(X, (\det\cE_i)^{\otimes k}\otimes\cL_i^{-\otimes kn_i})$ 
therefore gives an injective morphism
$\bS_{\lambda_i}\cV_i\hookrightarrow\bS_{\lambda_i}\cE_i$
of coherent $\OO_X$-modules. Each $\cV_i$ is nef and,
\[
 \bigotimes_{i=1}^r\cV_i
 \simeq
 \left(\bigotimes_{i=1}^r\cL_i\right)^{\oplus R}, \text{ where } R=\prod_{i=1}^r n_i.
\]
Thus $\bigotimes_{i=1}^r\cV_i$ is ample by hypothesis on $\bigotimes_{i=1}^r\cL_i$. Since $\bigoplus_{i=1}^{r}\cE_i=T_X$ is integrable as well, we see that Theorem~\ref{thm:Main3} applies whereby we obtain
\[
 X\simeq\prod_{i=1}^rZ_i,\qquad
 \cE_i\simeq\pr_i^*T_{Z_i},\qquad
 \cL_i^{\oplus n_i}\simeq\pr_i^*\cT_i,
\]
where each $\cT_i$ is an ample bundle on $Z_i$ as classified in Theorem~\ref{thm:Main3}. We can now determine the line bundles $\cL_i$. Fix a slice
$\iota_i:Z_i\hookrightarrow\prod_{j=1}^rZ_j$ obtained by fixing a point in every other factor, and let $H_i=\iota_i^*\cL_i$. Then
$\cT_i\simeq H_i^{\oplus n_i}$, so $H_i$ is ample. If
$Z_i\simeq\PP^{n_i}$ with $n_i\geq2$, write $H_i\simeq\OO_{\PP^{n_i}}(d_i)$. The determinant bound in Theorem~\ref{thm:Main1} imposes $n_i\leq n_id_i\leq n_i+1$, thereby forcing $d_i=1$. If $Z_i=Q_{n_i}$ with $n_i\geq2$, Theorem~\ref{thm:Main3} gives
$\cT_i\simeq\OO_{Q_{n_i}}(1)^{\oplus n_i}$ and so
$H_i\simeq\OO_{Q_{n_i}}(1)$.  Since $\cT_i\simeq H_i^{\oplus n_i}$, we therefore obtain $\cL_i^{\oplus n_i}\simeq(\pr_i^*H_i)^{\oplus n_i}$, yielding $\cL_i\simeq\pr_i^*H_i$ on $X$. Finally, the identities
\[
 \det T_{\PP^{n_i}}\simeq\OO_{\PP^{n_i}}(n_i+1),
 \qquad
 \det T_{Q_{n_i}}\simeq\OO_{Q_{n_i}}(n_i)
\]
show that $\cM_i\simeq\cL_i$ if $Z_i\simeq \PP^{n_i}$ and
$\cM_i\simeq\OO_X$ if $Z_i\simeq Q_{n_i}$. This proves \textup{(i)}.

We now prove (ii). For all $1\leq i\leq r$, if $\det\cE_i\otimes \cL_i^{\otimes-(n_i+1)}$ is nef, then the bundle $\cM_i=(\det\cE_i\otimes\cL_i^{-\otimes(n_i+1)})\otimes\cL_i$ is nef, whereby part (i) applies. If $Z_i\simeq Q_{n_i}$, then $\det\cE_i\otimes\cL_i^{-(n_i+1)}$ would restrict to
$H_i^{-1}$ on $Z_i$, which is not nef and contradicts the assumption of (ii). Thus the only possibility is $(Z_i,H_i, \det\cE_i)\simeq(\PP^{n_i},\OO_{\PP^{n_i}}(1), \cL_i^{\otimes (n_i+1)})$ for all $1\leq i\leq r$.

Finally we prove (iii). By assumption, every $\cM_i$ is numerically trivial and hence nef. Thus, we are in the setting of part (i). But in (i), if the factor $Z_i\simeq \PP^{n_i}$, then $\cM_i\simeq\cL_i$, which contradicts numerical triviality of $\cM_i$, since $\cL_i$ restricts to the ample line bundle $H_i$ on $Z_i$. Consequently, by the remaining possibility of (i), we have $(Z_i,H_i, \det\cE_i)\simeq(Q_{n_i},\OO_{Q_{n_i}}(1), \cL_i^{\otimes n_i})$ for all $1\leq i\leq r$, which completes the proof.
\end{proof}

We note that Theorem~\ref{specialcase2proof}(ii) can also be obtained using \cite[Proposition~4.5]{Watanabe26}. We finish with a product analogue of \cite{Ye94}.

\begin{theorem}\label{specialcase3proof}
Let $X$ be a smooth projective complex variety of dimension $n>0$ and let $r>0$ be an integer. Let $n_1,\dots,n_r$ be positive integers such that $T_X$ admits a holomorphic splitting $T_X=\bigoplus_{i=1}^r\cE_i$ with $\rk\cE_i=n_i$. Furthermore, assume that for each $1\leq i\leq r$, there exists a $\cL_i^{\otimes 2}$-valued everywhere nondegenerate holomorphic conformal form $q_i\in H^0\left(X,\Sym^2\cE_i\otimes\cL_i^{-2}\right)$ on $\cE_i$ for some nef line bundle $\cL_i$ such that $\bigotimes_{i=1}^r\cL_i$ is ample. Then, possibly after reordering, there is an isomorphism
\[
 X\simeq\prod_{i=1}^rQ_{n_i}
 \text{ such that }
 \cE_i\simeq\pr_i^*T_{Q_{n_i}}
 \text{ and }
 \cL_i\simeq\pr_i^*\OO_{Q_{n_i}}(1),
\]
where $\pr_i:X\rightarrow Q_{n_i}$ is the natural projection. 
If every $n_i$ is odd, the same conclusion holds if each $\cE_i$ admits an everywhere nondegenerate holomorphic conformal form $g_i\in H^0\left(X,\Sym^2\cE_i^*\otimes\cM_i\right)$ for any choice of nef line bundles $\cM_i$ such that $\bigotimes_{i=1}^r\cM_i$ is ample. In this case,
$\cM_i\simeq\pr_i^*\OO_{Q_{n_i}}(2)$.
\end{theorem}

\begin{proof}
The existence of a $\cL_i^{\otimes 2}$-valued everywhere nondegenerate holomorphic conformal form $q_i\in H^0\left(X,\Sym^2\cE_i\otimes\cL_i^{-2}\right)$ on $\cE_i$ is equivalent to an isomorphism
$\cE_i^*\otimes\cL_i^2\simeq\cE_i$, whereby taking determinants yields $(\det\cE_i)^2\simeq\cL_i^{2n_i}$. Consequently we obtain the numerical equivalence $\det\cE_i\equiv\cL_i^{n_i}$, which along with the other hypotheses of Theorem~\ref{specialcase3proof} gives the asserted product decomposition and the identification of the line bundles by Theorem~\ref{specialcase2proof}\textup{(iii)}.

We now show that if $\rk\cE_i=n_i$ are odd for all $1\leq i\leq r$, then it suffices for each $\cE_i$ to admit an everywhere nondegenerate holomorphic conformal form $g_i\in H^0\left(X,\Sym^2\cE_i^*\otimes\cM_i\right)$ for any choice of nef line bundles $\cM_i$. This is because taking determinants of the isomorphism $\cE_i\simeq\cE_i^*\otimes\cM_i$ yields
$ (\det\cE_i)^2\simeq\cM_i^{n_i}$. For all $1\leq i\leq r$, define
\[
 \cL_i:=\det\cE_i\otimes\cM_i^{-(n_i-1)/2},
\]
whereby the determinant identity $ (\det\cE_i)^2\simeq\cM_i^{n_i}$ implies $\cL_i^2\simeq\cM_i$. It follows that for all $1\leq i\leq r$, the line bundle $\cL_i$ is nef, since $\cM_i$ is. Furthermore,
 \[\left(\bigotimes_{i=1}^r\cL_i\right)^{\otimes 2}
 \simeq\bigotimes_{i=1}^r\cM_i
\]
is ample, so $\bigotimes_{i=1}^r\cL_i$ is ample.
The inverse of $g_i$ defines an everywhere nondegenerate section of
$\Sym^2\cE_i\otimes\cM_i^{-1}
\simeq\Sym^2\cE_i\otimes\cL_i^{-\otimes 2}$. Thus, the first assertion of Theorem~\ref{specialcase3proof} applies and yields the product of quadrics conclusion along with
$\cM_i\simeq\cL_i^2\simeq\pr_i^*\OO_{Q_{n_i}}(2)$ for all $1\leq i\leq r$.
\end{proof}

\printbibliography
\end{document}